\documentclass[3p]{elsarticle}
\journal{}
\usepackage{graphicx}
\usepackage{amsmath}
\usepackage{amssymb}
\usepackage{amsthm}
\usepackage{latexsym}
\usepackage{bm}
\usepackage{color}
\usepackage{stmaryrd}
\usepackage{mathrsfs}
\usepackage{array}
\usepackage{algorithm}
\usepackage{algpseudocode}
\usepackage{multirow}
\usepackage{url}
\SetSymbolFont{stmry}{bold}{U}{stmry}{m}{n}
\newtheorem{remark}{Remark}
\numberwithin{equation}{section}
\newcolumntype{P}[1]{>{\centering\arraybackslash}p{#1}}

\begin{document}

\begin{frontmatter}
\title{The nested block preconditioning technique for the incompressible Navier-Stokes equations with emphasis on hemodynamic simulations}

\author{Ju Liu}
\ead{liuju@stanford.edu}

\author{Weiguang Yang}
\ead{wgyang@stanford.edu}

\author{Melody Dong}
\ead{mldong@stanford.edu}

\author{Alison L. Marsden}
\ead{amarsden@stanford.edu}

\address{Department of Pediatrics (Cardiology), Department of Bioengineering, and Institute for Computational and Mathematical Engineering, Stanford University, Clark Center E1.3, 318 Campus Drive, Stanford, CA 94305, USA}

\begin{abstract}
We develop a novel iterative solution method for the incompressible Navier-Stokes equations with boundary conditions coupled with reduced models. The iterative algorithm is designed based on the variational multiscale formulation and the generalized-$\alpha$ scheme. The spatiotemporal discretization leads to a block structure of the resulting consistent tangent matrix in the Newton-Raphson procedure. As a generalization of the conventional block preconditioners, a three-level nested block preconditioner is introduced to attain a better representation of the Schur complement, which plays a key role in the overall algorithm robustness and efficiency. This approach provides a flexible, algorithmic way to handle the Schur complement for problems involving multiscale and multiphysics coupling. The solution method is implemented and benchmarked against experimental data from the nozzle challenge problem issued by the US Food and Drug Administration. The robustness, efficiency, and parallel scalability of the proposed technique are then examined in several settings, including moderately high Reynolds number flows and physiological flows with strong resistance effect due to coupled downstream vasculature models. Two patient-specific hemodynamic simulations, covering systemic and pulmonary flows, are performed to further corroborate the efficacy of the proposed methodology.
\end{abstract}


\end{frontmatter}


\section{Introduction}
\label{sec:introduction}

\subsection{Motivation and literature survey}
\label{sec:motivation}
Fast-growing interest in cardiovascular modeling \cite{Bao2014,Figueroa2017,Marsden2013a} and ever-increasing computing power \cite{top-500-list-site} has created a pressing need to efficiently simulate flow physics with high resolution. Indeed, running hemodynamic simulations with millions of unknowns has now become a routine part of scientific and clinical research and is now being used routinely in clinical decision-making \cite{Taylor2013,Zarins2013}. More and more, there is a growing consensus that proper design of the preconditioning techniques has become a critical issue for flow simulations, and in particular, hemodynamic simulations. In our opinion, under the setting of parallel computing, a desirable preconditioning technique for flow problems, and especially for hemodynamics, should possess the following attributes.
\begin{enumerate}
\item The algorithm should be robust with respect to physical parameters as well as spatiotemporal discretization methodologies.
\item The algorithm should be scalable in terms of fixed-size (strong) scalability as well as isogranular (weak) scalability.
\item The algorithm should require minimum user interventions. Ideally, users should not need to implement new matrices for the definition of the preconditioner.
\end{enumerate}
To achieve the above-mentioned goals, different approaches have been devised to precondition the matrix problem in the setting of the Krylov subspace method. There is a school that favors the idea of domain decomposition (DD) methods because they can be conveniently implemented in the parallel setting with the existing algebraic-based solver reused \cite{Hwang2005,Kong2019,Wu2014}. Although the DD method can be used as a black-box technique with almost no user intervention, its locality nature limits its parallel scalability \cite{Elman2008,Cyr2012}, and employing incomplete factorization in subdomains renders it non-robust for saddle-point problems \cite{Washio2005}. The recently introduced multi-level DD method borrows concepts from the multigrid method and seems to be a promising direction to overcome the scalability issue \cite{Tuminaro2002,Shadid2005}.

The physics-based block preconditioning technique, as an alternative approach, has attracted concentrated research for flow problems \cite{Benzi2005}. Consider solving a fully implicit scheme for the incompressible Navier-Stokes equations using the consistent Newton-Raphson method, the problem boils down to repeatedly solving a linear system of equations with a $2\times 2$ block structure,
\begin{align*}
\boldsymbol{\mathcal A} :=
\begin{bmatrix}
\boldsymbol{\mathrm A} & \boldsymbol{\mathrm B} \\[0.3mm]
\boldsymbol{\mathrm C} & \boldsymbol{\mathrm D}
\end{bmatrix}.
\end{align*}
In the above, $\boldsymbol{\mathrm A}$, $\boldsymbol{\mathrm B}$, and $\boldsymbol{\mathrm C}$ can be regarded as a discrete convection-diffusion-reaction operator, a discrete gradient operator, and a discrete divergence operator, with additional numerical modeling terms, respectively. The matrix $\boldsymbol{\mathrm D}$ arises purely due to the subgrid-scale modeling. The matrix $\boldsymbol{\mathcal A}$ can be factored into a lower triangular, a diagonal, and an upper triangular matrices as follows,
\begin{align*}
\boldsymbol{\mathcal A} = \boldsymbol{\mathcal L} \boldsymbol{\mathcal D} \boldsymbol{\mathcal U} =
\begin{bmatrix}
\boldsymbol{\mathrm I} & \boldsymbol{\mathrm O} \\[0.3em]
\boldsymbol{\mathrm C} \boldsymbol{\mathrm A}^{-1} & \boldsymbol{\mathrm I}
\end{bmatrix}
\begin{bmatrix}
\boldsymbol{\mathrm A} & \boldsymbol{\mathrm O} \\[0.3em]
\boldsymbol{\mathrm O} & \boldsymbol{\mathrm S}
\end{bmatrix}
\begin{bmatrix}
\boldsymbol{\mathrm I} & \boldsymbol{\mathrm A}^{-1} \boldsymbol{\mathrm B} \\[0.3em]
\boldsymbol{\mathrm O} & \boldsymbol{\mathrm I}
\end{bmatrix},
\end{align*}
with $\boldsymbol{\mathrm S} := \boldsymbol{\mathrm D} - \boldsymbol{\mathrm C} \boldsymbol{\mathrm A}^{-1} \boldsymbol{\mathrm B}$ being the Schur complement. Therefore, the design of the preconditioner for $\boldsymbol{\mathcal A}$ reduces to solving smaller systems associated with $\boldsymbol{\mathrm A}$ and $\boldsymbol{\mathrm S}$. This technique is attractive because it combines the merits of the Chorin-Teman projection method \cite{Chorin1968,Kim1985} used in the finite volume community and the fully implicit method \cite{Gresho1998,Jansen2000} used in the finite element community. Roughly speaking, the physics-based block preconditioner can be regarded as an algebraic procedure that wraps the projection method inside the Krylov iteration. The Schur complement $\boldsymbol{\mathrm S}$ can be viewed as an algebraic manifestation of the pressure Poisson equation in the fully implicit scheme. Unlike the classical projection method, the separation of the physical fields does not take place in the temporal scheme, thus avoiding considerations of artificial pressure boundary conditions and time step size control \cite{Guermond2006}.

A solution method is deemed to be robust with respect to physical parameters if there is no significant impact on its convergence rate with varying physical parameters. For flow problems, this parameter is typically the Reynolds number, which, in hemodynamic simulations, ranges from $\mathcal O(10^{-3})$ in the capillary vessels to $\mathcal O(10^3)$ in the aorta. The classical SIMPLE preconditioner extracts the diagonal of $\boldsymbol{\mathrm A}$ to define a sparse approximation of the Schur complement. This approach ignores the convection information and is therefore non-robust with respect to the Reynolds number. In fact, it has been a major research thrust to search for a preconditioner that is insensitive to the Reynolds number and at the same time remains scalable in the computational fluid dynamics (CFD) community using fully implicit schemes \cite{Elman2008,Cyr2012,Elman2003,Deparis2014,Elman2006,Kay2002,Shadid2016,Silvester2001,Turek1999}. Representative examples include the pressure convection-diffusion (PCD) preconditioner \cite{Kay2002} and the least-squares commutator (LSC) preconditioner \cite{Elman2006}. Both are known to be scalable with respect to discretization resolution and are mildly affected by the Reynolds number. However, they are not without shortcomings. The major drawback of the PCD preconditioner is that it requires the assembly of a new matrix, which is implementationally inconvenient and computationally inefficient. The LSC preconditioner was proposed to remedy that issue by only using the existing blocks in $\boldsymbol{\mathcal A}$. However, the LSC preconditioner was proposed and tested based on a mixed formulation using inf-sup stable elements. It has been recently noticed that the LSC preconditioner is quite sensitive to the spatial discretization method, as it cannot converge for the stabilized formulation using equal-order interpolations \cite{Cyr2012}. The drawbacks of the PCD and LSC preconditioners indicate that the satisfaction of the three attributes listed at the beginning of this section remains a challenging task for the incompressible Navier-Stokes equations.

\subsection{A three-level nested block preconditioner}
\begin{figure}
	\begin{center}
	\begin{tabular}{c}
\includegraphics[angle=0, trim=0 0 0 0, clip=true, scale = 0.45]{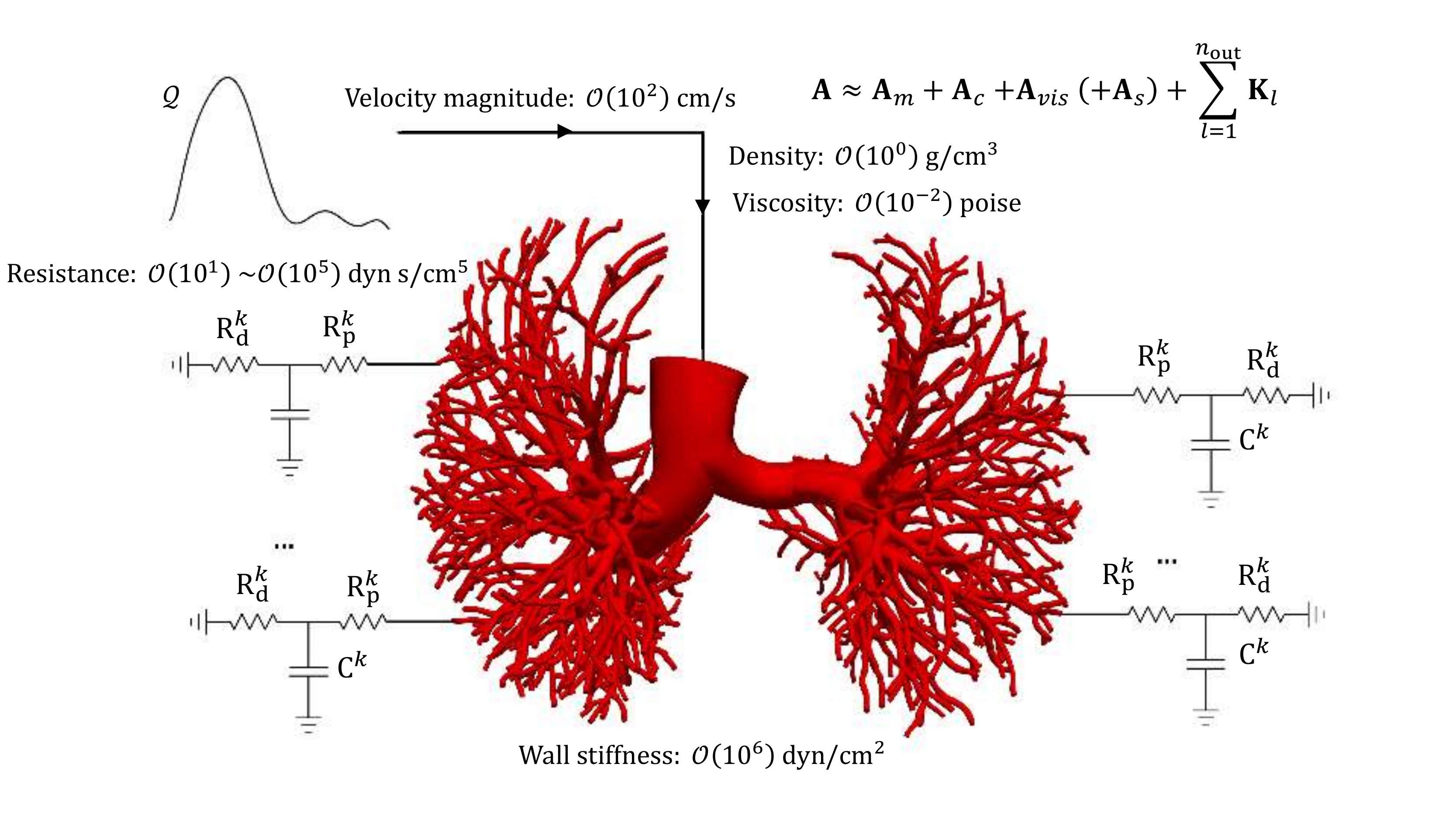}
\end{tabular}
\end{center}
\caption{Illustration of the hemodynamic simulation setting and the structure of the sub-matrix $\boldsymbol{\mathrm A}$. This sub-matrix can be approximately decomposed into matrices arise from the transient term ($\boldsymbol{\mathrm A}_m$), the convection term ($\boldsymbol{\mathrm A}_c$), the viscous term ($\boldsymbol{\mathrm A}_{vis}$), the reduced models on the outflow boundaries ($\boldsymbol{\mathrm K}_l$), and the wall stiffness term if one is considering FSI ($\boldsymbol{\mathrm A}_s$). These terms involve physical parameters and the typical magnitude of the physical parameters in cardiovascular simulations are given in the centimeter-gram-second (CGS) units. } 
\label{fig:hemodynamics_and_algebra}
\end{figure}

In the setting of blood flow simulations, the situation is more complicated. In addition to the convection term, the downstream vasculature is often modeled as a reduced model \cite{Peiro2009}, which is coupled to the three-dimensional problem via the outflow boundary conditions. This coupling strategy is also referred to as the geometric multiscale modeling \cite{Quarteroni2016}, and it leads to a modification of the original matrix problem with rank-one matrices multiplied by a factor proportional to the resistance value on the outlet surfaces (see Sec. \ref{sec:predictor multicorrector algorithm}). For realistic problems, the resistance value is not small, and correspondingly the reduced models have a significant impact on the matrix properties. If we look one step further, more challenges arise in the setting of fluid-structure interaction (FSI). Considering the coupled momentum formulation as an example, an additional wall stiffness matrix enters into the tangent matrix \cite{Figueroa2006}. Again, this modification cannot be neglected because the physiologically realistic Young's modulus is by no means small. The contributions to the tangent matrix from different physical sources are illustrated in Fig. \ref{fig:hemodynamics_and_algebra}. \textit{Consequently, over and above the long-lasting efforts for preconditioning the incompressible Navier-Stokes equations, special consideration needs to be further exercised so that the geometric multiscale models and multiphysics coupling, such as FSI, can be properly taken into account in preconditioner design.}

If one examines the resulting algebraic system carefully, it can be noticed that the rank-one matrices and the wall stiffness matrix enter into the block matrix $\boldsymbol{\mathrm A}$ only, without affecting the other three block matrices. This means that the matrix $\boldsymbol{\mathrm A}$ contains multiple contributions, including the transient term, the viscous term, the convection term, the subgrid scale modeling terms, the rank-one modifications defined on the outlet surfaces, and the stiffness matrix defined on the wall surface if one is solving an FSI problem. The true difficulty comes from designing an effective sparse approximation of the Schur complement that incorporates the information from the above terms. Due to the special algebraic structure of the rank-one matrices, the Sherman-Morrison formula has been utilized to design a Schur complement approximation for hemodynamic simulations \cite{Moghadam2013a}. That approach, together with the bi-partitioned iterative algorithm (BIPN) \cite{Moghadam2015}, constitute the backbone solver technology in the current \texttt{svSolver} \cite{simvascular-simulation-guide}, the CFD solver in the SimVascular package \cite{Updegrove2017}. Yet, it is known that the same solver technology performs poorly when the vessel wall is deformable, and the reason is apparently due to the ignorance of the wall stiffness matrix in the sparse approximation of the Schur complement. Also, one may reasonably expect performance degradation when simulating strong convection problems with \texttt{svSolver} since it does not account for the convection terms \cite{Moghadam2013a}. As mentioned above, designing a robust and scalable preconditioner for the convection term alone remains an unsettled issue. \textit{Therefore, it will be quite challenging, if not impossible, to explicitly construct an algebraic form of the preconditioning technique accounting for all aforementioned terms simultaneously.} This leads us to alternatively consider constructing the Schur complement approximation via an algorithm, rather than an algebraic form.

Recently, it has been realized that one may apply the so-called Schur complement reduction procedure (SCR) to serve as a preconditioning technique. In the original SCR procedure, one solves the sub-matrices $\boldsymbol{\mathrm A}$ and $\boldsymbol{\mathrm S}$ to a prescribed tolerance and thereby solves the linear problem associated with $\boldsymbol{\mathcal A}$ in one pass \cite{Benzi2005,May2008}. Although the excessive memory cost for storing the Schur complement can be resolved by using a matrix-free algorithm to define the action of $\boldsymbol{\mathrm S}$ on a vector, solving a linear system associated with $\boldsymbol{\mathrm S}$ to a high precision is still prohibitively expensive. To alleviate this issue, one may use the SCR procedure as a preconditioning technique by wrapping it inside an iterative solution method, which leads to a three-level nested algorithm structure (see Sec. \ref{sec:SCR}). In the outer level, the iterative method strives to solve $\boldsymbol{\mathcal A}$ either via a static iterative method or a Krylov subspace method; in the intermediate level, the sub-matrices $\boldsymbol{\mathrm A}$ and $\boldsymbol{\mathrm S}$ are solved (not necessarily to high precision) as a preconditioning technique to accelerate the outer iteration; in the inner level, the matrix-free algorithm for $\boldsymbol{\mathrm S}$ necessitates solving with $\boldsymbol{\mathrm A}$. The inner level solver can be a key ingredient when the matrix $\boldsymbol{\mathrm A}$ involves complex contributions from non-traditional sources. Although the introduction of an inner solver may ostensibly increase the computational burden, it can in fact dramatically enhance the solver robustness without losing efficiency if the setting for the inner solver is properly tuned (see Sec. \ref{sec:patient_specific_examples}). In the conventional physics-based block preconditioner, the sparse approximation of $\boldsymbol{\mathrm S}$ is explicitly constructed without invoking the inner level solver, and thence exhibits a two-level structure \cite{Cyr2016}. In theory, the proposed three-level nested preconditioning technique can be viewed as a generalization of the conventional block preconditioners.

The SCR procedure has been used as a preconditioning technique for CFD problems within the Richardson iteration scheme \cite{Manguoglu2008,Manguoglu2011} and  the biconjugate gradient stabilized method \cite{Manguoglu2009}. Those results clearly justified the advantage of using SCR as a preconditioning technique over several standard preconditioners. In our prior work, we investigated the use of FGMRES \cite{Saad1993} preconditioned by the SCR procedure for hyperelasticity \cite{Liu2019}, based on a unified continuum modeling framework \cite{Liu2018,Liu2019a}. New ingredients were added to further enhance the three-level nested algorithm previously introduced in \cite{Manguoglu2008,Manguoglu2011,Manguoglu2009}. In addition to using the FGMRES algorithm as the outer solver, we applied the GMRES algorithm preconditioned by the algebraic multigrid (AMG) preconditioner at the intermediate and inner levels. This combined the merits of both the conventional block preconditioner \cite{Cyr2012,Deparis2014,Elman2008} and the nested algorithm proposed in \cite{Manguoglu2008,Manguoglu2011,Manguoglu2009}. Also, a sparse matrix was constructed to precondition and accelerate the matrix-free solution procedure for $\boldsymbol{\mathrm S}$. In this work, we further investigate the efficacy of this solution method for hemodynamic problems, by examining its performance when convection and the geometric multiscale coupling contributions are significant. This study serves as a stepping stone towards a novel iterative solution method for FSI problems, based on our recently proposed FSI framework \cite{Liu2018} and related iterative solution method for hyperelasticity \cite{Liu2019}.

\subsection{Remarks on the numerical formulation}
We want to point out that the numerical formulation adopted in this work is different from those reported in the existing literature in several aspects. First, in several prior works \cite{Figueroa2006,Taylor1998}, an integration-by-parts was performed for the divergence operator in the continuity equation. We do not favor this approach because it contradicts with the setting of the pressure function space in the Galerkin formulation \cite{Gresho1998}. Second, the pressure variable has traditionally been evaluated by the backward Euler method in the generalized-$\alpha$ scheme \cite{Bazilevs2007a}. This choice will degrade the temporal accuracy, and this issue has recently been rectified by evaluating the pressure at the intermediate time step \cite{Liu2020a}. Third, in the same vein, the traction forces used to couple the three-dimensional model with the reduced model on the outlet surfaces are evaluated at the intermediate step as well, in contrast to the prior approach \cite{Moghadam2013}. Fourth, the definition of the stabilization parameter in the variational multiscale formulation is modified for simplex elements. This modification makes the stabilization parameter remain invariant under node renumbering and is recently introduced in \cite{Pauli2017,Danwitz2019}.

\section{Governing equations and the spatiotemporal discretizations}
\label{sec:continuum_problem_and_discretization}
In this section, we introduce the strong-form problem and the fully discrete problem generated by the variational multiscale formulation and the generalized-$\alpha$ method.

\subsection{Strong-form problem}
Let $\Omega \subset \mathbb R^3$ be a fixed bounded open set with sufficiently smooth (e.g. Lipschitz) boundary $\Gamma := \partial \Omega$. The time interval is denoted $(0,T) \subset \mathbb R$ with $T>0$. The governing equations for the incompressible flow of a Newtonian fluid can be stated as follows.
\begin{align}
\label{eq:ns_mom}
\bm 0 &= \rho \frac{\partial \bm v}{\partial t}  + \rho  \bm v \cdot \nabla \bm v - \nabla \cdot \bm \sigma - \rho \bm f, && \mbox{ in } \Omega \times (0,T), \\
\label{eq:ns_mass}
0 &= \nabla \cdot \bm v, && \mbox{ in } \Omega \times (0, T),
\end{align}
wherein
\begin{align}
\label{eq:ns_sigma}
\bm \sigma := 2\mu \bm \varepsilon(\bm v) - p \bm I, \quad \bm \varepsilon(\bm v) := \frac12 \left( \nabla \bm v + \nabla \bm v^T \right).
\end{align}
In the above, $\rho$ is the fluid density, $\bm v$ is the velocity field, $\bm \sigma$ is the Cauchy stress, $\bm f$ is the body force, $\mu$ is the dynamic viscosity, $\bm \varepsilon$ is the rate-of-strain tensor, $p$ is the pressure, and $\bm I$ is the second-order identity tensor. The initial condition is given by a divergence-free velocity field $\bm v_0$ as
\begin{align}
\label{eq:initial_condition}
\bm v(\cdot, 0) = \bm v_0(\cdot), && \mbox{ in } \bar{\Omega}.
\end{align}
The boundary $\Gamma$ can be partitioned into two non-overlapping subdivisions, that is,
\begin{align}
\label{eq:gamma_partition}
\overline{\Gamma} = \overline{\Gamma_g \cup \Gamma_h}, \mbox{ and } \emptyset = \Gamma_g \cap \Gamma_h.
\end{align}
In the above, the subscripts $g$ and $h$ indicate the Dirichlet and Neumann partitions, respectively. The unit outward normal vector to $\Gamma$ is denoted as $\bm n$. Given the Dirichlet data $\bm g$ and the boundary traction $\bm h$, the boundary conditions can be stated as
\begin{align}
\label{eq:dirichlet_bc}
\bm v = \bm g && \mbox{ on } \Gamma_{g} \times (0,T), \\
\label{eq:neumann_bc}
\bm \sigma \bm n = \bm h  && \mbox{ on } \Gamma_{h} \times (0,T).
\end{align}

\begin{figure}
	\begin{center}
	\begin{tabular}{c}
\includegraphics[angle=0, trim=30 0 30 20, clip=true, scale = 0.25]{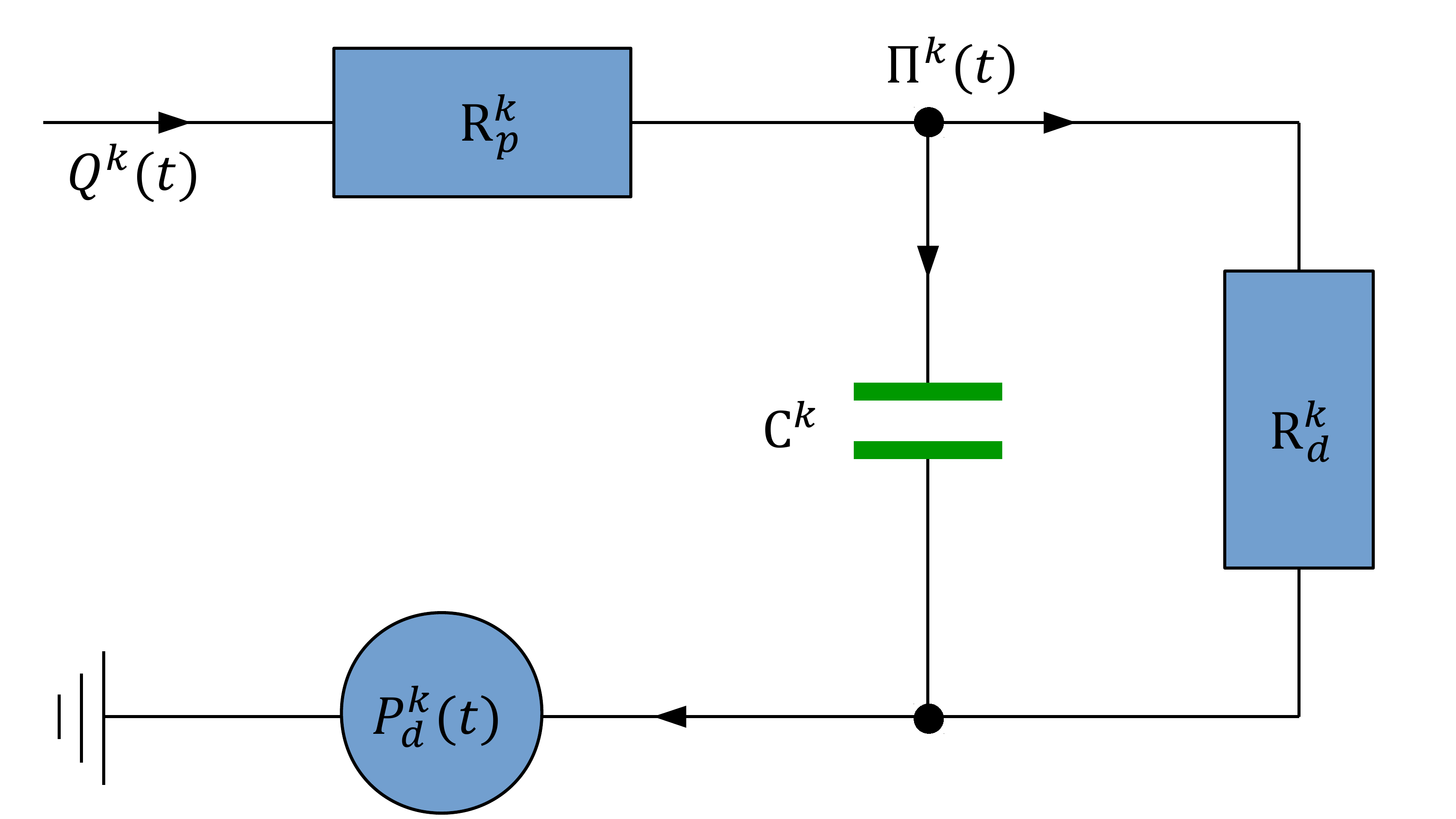}
\end{tabular}
\end{center}
\caption{Schematic representation of the three element Windkessel electric model.} 
\label{fig:RCR_circuit}
\end{figure}

\subsection{The boundary conditions and geometric multiscale modeling}
In the mathematical modeling of cardiovascular biofluids, the description of the boundaries needs to be further refined. We consider that $\Gamma_g$ can be further partitioned into two non-overlapping subdivisions: $\overline{\Gamma_g} = \overline{\Gamma_{\mathrm{in}} \cup \Gamma_{\mathrm{wall}} }$. On $\Gamma_{\mathrm{in}}$, the velocity is specified by a prescribed velocity profile $\bm v_{\mathrm{in}}$, while on $\Gamma_{\mathrm{wall}}$, we impose no-slip boundary condition. This suggests that the Dirichlet data $\bm g$ can be explicitly defined as,
\begin{align*}
\bm g = 
\begin{cases}
\bm v_{\mathrm{in}} & \mbox{ on } \Gamma_{\mathrm{in}}, \\
\bm 0 & \mbox{ on } \Gamma_{\mathrm{wall}}.
\end{cases}
\end{align*}
The velocity profiles on the inlet surfaces are typically Poiseuille or Womersley. For patient-specific geometries, the inlet surfaces do not necessarily have a circular shape, and we invoke a special mapping technique introduced in \cite{Takizawa2010} to generate the desired inflow profile on $\Gamma_{\mathrm{in}}$.

Specialized downstream boundary conditions are essential for capturing the physiologically realistic flow and pressure conditions, and here we use lumped parameter models as the geometric multiscale model to describe the downstream circulation in this work, though other reduced models (e.g. one-dimensional models) can also be employed. Consider the Neumann boundary $\Gamma_h$ that assumes $\mathrm n_{\mathrm{out}}$ non-overlapping planar surfaces,
\begin{align*}
\Gamma_h = \bigcup_{k=1}^{\mathrm n_{\mathrm{out}}} \Gamma^{k}_{\mathrm{out}}, \qquad \overline{\Gamma^i_{\mathrm{out}}} \cap \overline{\Gamma^j_{\mathrm{out}}} = \emptyset, \mbox{ for } 1 \leq i,j \leq \mathrm n_{\mathrm{out}}, i \neq j.
\end{align*}
On each outlet surface, the traction boundary condition is defined as
\begin{align*}
\bm h(t) = -P^k(t) \bm n \qquad \mbox{ on } \Gamma_{\mathrm{out}}^{k}.
\end{align*}
wherein $P^k(\cdot) : \mathbb R_{+} \rightarrow \mathbb R$ is a scalar function of time. We define the flow rate on an outlet surface as
\begin{align}
Q^k(t) := \int_{\Gamma_{\mathrm{out}}^k} \bm v \cdot \bm n d\Gamma.
\end{align}
By invoking the Dirichlet-to-Neumann method, the values of $P^k$ at each time instant are determined by solving a set of algebraic-differential equations with the flow rates as inputs. As an example, we consider the following equations known as the three-element Windkessel model,
\begin{align}
\label{eq:RCR_circuit_1}
& \frac{d\Pi^k(t)}{dt} =\mathcal{F}(\Pi^k(t), Q^k(t)) := -\frac{\Pi^k(t)}{\mathrm R_{\mathrm d}^k\mathrm C^k } + \frac{Q^k(t)}{\mathrm C^k }, \\
\label{eq:RCR_circuit_2}
& P^k(t) = \mathrm R_{\mathrm p}^k Q^k(t) + \Pi^k(t) + P^k_{\mathrm d}(t).
\end{align}
In the above, $\mathrm C^k$, $\mathrm R_{\mathrm d}^k$, and $\mathrm R_{\mathrm p}^k$ are constants and represent the compliance, distal resistance, and proximal resistance of the downstream vasculature to each outlet $\Gamma^k_{\mathrm{out}}$, $\Pi^k(t)$ represents the pressure drop across the distal resistor, and $P^k_{\mathrm d}(t)$ represents the distal reference pressure. An electric circuit analogue of this model is illustrated in Fig. \ref{fig:RCR_circuit}. For this Windkessel model, one may obtain an analytic representation of $P^{k}(t)$ in terms of $Q^k(t)$ as follows,
\begin{align}
\label{eq:RCR_analytic_formula}
P^k(t) = \int_0^t \left( \exp(-\frac{t-s}{\mathrm R_{\mathrm d}^k \mathrm C^k}) \frac{Q^k(s)}{\mathrm C^k} \right) ds + \mathrm R_{\mathrm p}^k Q^k(t) + P^k_{\mathrm d}(t) + \exp(-\frac{t}{\mathrm R_{\mathrm d}^k \mathrm C^k})\left( P^k(0) - \mathrm R_{\mathrm p}^kQ^k(0) -P^k_{\mathrm d}(0) \right).
\end{align}
As a special case, when $\mathrm C^k \rightarrow 0^+$, the three-element Windkessel model \eqref{eq:RCR_circuit_1}-\eqref{eq:RCR_circuit_2} reduces to a simpler resistance boundary condition,
\begin{align}
\label{eq:resistance_bc}
P^k(t) = \mathrm R^k  Q^k(t) + P^k_{\mathrm d}(t),
\end{align}
with $\mathrm R^k := \mathrm R_{\mathrm p}^k + \mathrm R_{\mathrm d}^k$.

\begin{remark}
In this work, we focus on using the Dirichlet-to-Neumann method to model the physics of the downstream domain. It is worth pointing out that there exists another coupling approach, in which flow rates are solved via algebraic-differential equations and are passed to the three-dimensional model by imposing velocity profiles on boundary surfaces. For example, this method is used when coupling a heart model with the inlet of an aorta.
\end{remark}

\subsection{Variational multiscale formulation}
The semi-discrete formulation is constructed based on the residual-based variational multiscale method \cite{Bazilevs2007a}. Let $\mathcal S_{\bm v}$ and $\mathcal S_{p}$ denote the trial solution spaces of the fluid velocity and pressure, and let $\mathcal V_{p}$ and $\mathcal V_{\bm v}$ be the corresponding test function spaces. These spaces are spanned by finite element basis functions, and in this work, linear polynomials. The Dirichlet boundary condition \eqref{eq:dirichlet_bc} is built into the definition of $\mathcal S_{\bm v}$. The union of element interiors is denoted as $\Omega^{\prime}$. The semi-discrete formulation can be stated as follows. Find $\bm y_h(t):= \left\lbrace \bm v_h(t), p_h(t) \right\rbrace^T \in \mathcal S_{\bm v} \times \mathcal S_{p}$ such that
\begin{align}
\label{eq:semi-discrete-v}
& \mathbf B_m \left( \bm w_h ;  \dot{\bm y}_h, \bm y_h \right) = 0, && \forall \bm w_h \in \mathcal V_{\bm v}, \\
\label{eq:semi-discrete-p}
& \mathbf B_p\left( q_h; \dot{\bm y}_h, \bm y_h \right) = 0, && \forall q_h \in \mathcal V_{p}, 
\end{align}
where
\begin{align}
\label{eq:fsi_fluid_residual_based_vms_momentum}
& \mathbf B_m \left( \bm w_h ;  \dot{\bm y}_h, \bm y_h \right) := \mathbf B_m^{\textup{vol}} \left( \bm w_h ;  \dot{\bm y}_h, \bm y_h \right) + \mathbf B_m^{\textup{bc}} \left( \bm w_h ;  \dot{\bm y}_h, \bm y_h \right) + \mathbf B_m^{\textup{bf}} \left( \bm w_h ;  \dot{\bm y}_h, \bm y_h \right), \displaybreak[2] \\
& \mathbf B_m^{\textup{vol}} \left( \bm w_h ;  \dot{\bm y}_h, \bm y_h \right) := \int_{\Omega} \bm w_h \cdot \rho \left( \frac{\partial \bm v_h}{\partial t} + \bm v_h \cdot \nabla \bm v_h - \bm b \right) d\Omega  - \int_{\Omega} \nabla \cdot \bm w_h p_h d\Omega_{\bm x} + \int_{\Omega} 2\mu  \bm \varepsilon(\bm w_h) : \bm \varepsilon(\bm v_h)  d\Omega \nonumber \displaybreak[2] \\
& \hspace{1mm}  - \int_{\Omega^{\prime}} \nabla \bm w_h : \left( \rho \bm v^{\prime} \otimes \bm v_h \right) d\Omega + \int_{\Omega^{\prime}} \nabla \bm v_h : \left( \rho \bm w_h \otimes \bm v^{\prime} \right) d\Omega  - \int_{\Omega^{\prime}} \nabla \bm w_h : \left( \rho \bm v^{\prime} \otimes \bm v^{\prime} \right) d\Omega_{\bm x} - \int_{\Omega^{\prime}} \nabla \cdot \bm w_h p^{\prime} d\Omega, \displaybreak[2] \\
& \mathbf B_m^{\textup{bc}} \left( \bm w_h ;  \dot{\bm y}_h, \bm y_h \right) := - \int_{\Gamma_h} \bm w_h \cdot \bm h d\Gamma, \displaybreak[2] \\
\label{eq:fsi_back_flow_stabilization}
& \mathbf B_m^{\textup{bf}} \left( \bm w_h ;  \dot{\bm y}_h, \bm y_h \right) := - \int_{\Gamma_h}  \rho \beta \left(\bm v_h \cdot \bm n \right)_{-} \bm w_h \cdot \bm v_h d\Gamma, \displaybreak[2] \\
\label{eq:fsi_fluid_residual_based_vms_mass}
& \mathbf B_p\left( q_h ; \dot{\bm y}_h, \bm y_h \right) := \int_{\Omega} q_h \nabla \cdot \bm v_h d\Omega  - \int_{\Omega^{\prime}} \nabla q_h \cdot  \bm v^{\prime} d\Omega, \displaybreak[2] \\
& \bm v^{\prime} := -\bm \tau_{M} \left( \rho \frac{\partial \bm v_h}{\partial t} + \rho \bm v_h \cdot \nabla \bm v_h  + \nabla p_h - \mu \Delta \bm v^h - \rho \bm b \right), \displaybreak[2] \\
& p^{\prime} := -\tau_C \nabla \cdot \bm v_h, \displaybreak[2] \\
& \bm \tau_M := \tau_M \bm I, \displaybreak[2] \\
\label{eq:fsi_fluid_def_tau_m}
& \tau_M := \frac{1}{\rho}\left( \frac{\mathrm C_{\mathrm T}}{\Delta t^2} + \bm v_h \cdot \bm G \bm v_h + \mathrm C_{\mathrm I} \left( \frac{\mu}{\rho} \right)^2 \bm G : \bm G \right)^{-\frac12}, \displaybreak[2] \\
& \tau_C := \frac{1}{\tau_M \textup{tr}\bm G}, \displaybreak[2] \\
& G_{ij} := \sum_{k=1}^{3} \frac{\partial \xi_k}{\partial x_i} M_{kl} \frac{\partial \xi_l}{\partial x_j}, \displaybreak[2] \\
\label{eq:def_K_for_scale_G}
& \bm M = [ M_{kl} ] = \frac{\sqrt[3]{2}}{2}\begin{bmatrix}
2 & 1 & 1 \\
1 & 2 & 1 \\
1 & 1 & 2
\end{bmatrix}, \displaybreak[2] \\
& \bm G : \bm G := \sum_{i,j=1}^{3} G_{ij} G_{ij}, \displaybreak[2] \\
& \textup{tr}\bm G := \sum_{i=1}^{3} G_{ii}, \displaybreak[2] \\
& \left( \bm v_h \cdot \bm n \right)_{-} := \frac{\bm v_h \cdot \bm n - |\bm v_h \cdot \bm n|}{2} = 
\begin{cases}
\bm v_h \cdot \bm n & \quad \mbox{ if } \bm v_h \cdot \bm n < 0, \\
0 & \quad \mbox{ if } \bm v_h \cdot \bm n \geq 0.
\end{cases}
\end{align}
In the above, $\bm \xi = \left\lbrace \xi_i \right\rbrace_{i=1}^{3}$ are the coordinates of an element in the parent domain; the value of $\mathrm C_{\mathrm I}$ relies on the polynomial order of the interpolation basis functions and takes the value of $36$ for linear interpolations \cite{Franca1992,Figuero2006}; the value of $\mathrm C_{\mathrm T}$ is taken to be $4$. The term $\mathbf B_m^{\textup{bf}}$ is an additional term added to enhance the overall numerical robustness in the presence of locally reversed flows near the outlet surfaces \cite{Bazilevs2009a,Esmaily-Moghadam2011}. The parameter $\beta$ is non-dimensional, and its value is fixed to be $0.2$ in this work, following the practices in \texttt{svSolver} \cite{simvascular-simulation-guide}.

\begin{remark}
In \eqref{eq:def_K_for_scale_G}, $\bm M$ is introduced for simplex elements because their standard reference elements are not symmetric, and nodal permutations may lead to changes in the definition $\bm G$ without using $\bm M$. The entries of $\bm M$ are obtained by mapping the reference element to a regular simplex without changing the element volume. For a detailed derivation of $\bm M$ and its variants in two- and four-dimensions, readers are referred to \cite{Pauli2017,Danwitz2019}.
\end{remark}

\begin{remark}
The term \eqref{eq:fsi_back_flow_stabilization} is introduced as a modification of the Neumann boundary condition to compensate for the incoming energy due to locally reversed flows. It can be shown that choosing $\beta = 1.0$ leads to an energy stable formulation \cite{Bertoglio2016}. In practice, it is found numerically that choosing $\beta$ smaller than $1.0$ is often sufficient for robust performances \cite{Esmaily-Moghadam2011} and is apparently less intrusive for flow physics. For a recent summary of the treatment for backflow instabilities, readers are referred to \cite{Bertoglio2018}.
\end{remark}

\subsection{Temporal discretization}
To derive the fully discrete formulation for problems with geometric multiscale coupled boundary conditions, we apply two time-stepping schemes for the three-dimensional and zero-dimensional domains separately. In the three-dimensional domain, the generalized-$\alpha$ method is applied; in the zero-dimensional domain, the boundary traction is obtained by an explicit fourth-order Runge-Kutta method with subdivided time intervals \cite{Moghadam2013}.

\subsubsection{Temporal discretization of the three-dimensional problem}
\label{sec:temporal_scheme_3d_problem}
Let the time interval $(0,T)$ be divided into a set of $N_{\mathrm{ts}}$ subintervals of size $\Delta t_n := t_{n+1} - t_n$, which is delimited by a discrete time vector $\left\lbrace t_n \right\rbrace_{n=0}^{N_{\mathrm{ts}}}$. The approximations to the solution vector and its first time derivative evaluated at the time step $t_n$ are denoted as $\bm y_n := \left\lbrace \bm v_n, p_n \right\rbrace^T$ and $\dot{\bm y}_n:= \left\lbrace \dot{\bm v}_n, \dot{p}_n \right\rbrace^T$. The approximation to $\bm h(t)$ at time $t_n$ is denoted as $\bm h_n$. Let $\bm e_i$ be the Cartesian basis vector with $i=1,2,3$ and $N_A$ be the basis function of the discrete function spaces for the velocity component as well as pressure. We may then define the residual vectors as 
\begin{align*}
\boldsymbol{\mathrm R}^{\textup{vol}}_m\left( \dot{\bm y}_{n}, \bm y_{n} \right) :=& \left\lbrace \mathbf B_m^{\textup{vol}}\left( N_A \bm e_i ;  \dot{\bm y}_{n}, \bm y_{n} \right) \right\rbrace,  && 
\boldsymbol{\mathrm R}_m^{\textup{bc}}\left( \dot{\bm y}_{n}, \bm y_{n} \right) := \left\lbrace \mathbf B_m^{\textup{bc}}\left( N_A \bm e_i ;  \dot{\bm y}_{n}, \bm y_{n} \right) \right\rbrace, \\
\boldsymbol{\mathrm R}^{\textup{bf}}_m\left( \dot{\bm y}_{n}, \bm y_{n} \right) :=& \left\lbrace \mathbf B_m^{\textup{bf}}\left( N_A \bm e_i ;  \dot{\bm y}_{n}, \bm y_{n} \right) \right\rbrace,  &&
\boldsymbol{\mathrm R}_m\left( \dot{\bm y}_{n}, \bm y_{n} \right) := \left\lbrace \mathbf B_m\left( N_A \bm e_i ;  \dot{\bm y}_{n}, \bm y_{n} \right) \right\rbrace, \\
\boldsymbol{\mathrm R}_p\left( \dot{\bm y}_{n}, \bm y_{n} \right) :=& \left\lbrace \mathbf B_p\left( N_A ;  \dot{\bm y}_{n}, \bm y_{n} \right) \right\rbrace.
\end{align*}
With the above notations, we have $\boldsymbol{\mathrm R}_m\left( \dot{\bm y}_{n}, \bm y_{n} \right) = \boldsymbol{\mathrm R}^{\textup{vol}}_m\left( \dot{\bm y}_{n}, \bm y_{n} \right) + \boldsymbol{\mathrm R}^{\textup{bc}}_m\left( \dot{\bm y}_{n}, \bm y_{n} \right) + \boldsymbol{\mathrm R}^{\textup{bf}}_m\left( \dot{\bm y}_{n}, \bm y_{n} \right)$ due to their definitions. The fully discrete scheme can be stated as follows. At time step $t_n$, given $\dot{\bm y}_n$, $\bm y_n$, and the time step size $\Delta t_n$, find $\dot{\bm y}_{n+1}$ and $\bm y_{n+1}$ such that
\begin{align}
\label{eq:gen_alpha_fully_discrete_linear_momentum}
& \boldsymbol{\mathrm R}_m(\dot{\bm y}_{n+\alpha_m}, \bm y_{n+\alpha_f}) = \bm 0, \displaybreak[2] \\
\label{eq:gen_alpha_fully_discrete_pressure_rate}
& \boldsymbol{\mathrm R}_p(\dot{\bm y}_{n+\alpha_m}, \bm y_{n+\alpha_f}) = \bm 0, \displaybreak[2] \\
\label{eq:gen_alpha_def_y_n_plus_1}
& \bm y_{n+1} = \bm y_{n} + \Delta t_n \dot{\bm y}_n, + \gamma \Delta t_n \left( \dot{\bm y}_{n+1} - \dot{\bm y}_{n}\right), \displaybreak[2] \\
\label{eq:gen_alpha_def_y_n_alpha_m}
& \dot{\bm y}_{n+\alpha_m} = \dot{\bm y}_{n} + \alpha_m \left(\dot{\bm y}_{n+1} - \dot{\bm y}_{n} \right), \displaybreak[2] \\
\label{eq:gen_alpha_def_y_n_alpha_f}
& \bm y_{n+\alpha_f} = \bm y_{n} + \alpha_f \left( \bm y_{n+1} - \bm y_{n} \right).
\end{align}
In the above system, there are three parameters $\alpha_m$, $\alpha_f$ and $\gamma$, whose values determine critical numerical properties of the discrete dynamic system. To ensure second-order accuracy and unconditional stability (for linear problems), and optimal high frequency dissipation, they are parametrized as
\begin{align*}
\alpha_m = \frac12 \left( \frac{3-\varrho_{\infty}}{1+\varrho_{\infty}} \right), \quad \alpha_f = \frac{1}{1+\varrho_{\infty}}, \quad \gamma = \frac{1}{1+\varrho_{\infty}},
\end{align*}
wherein $\varrho_{\infty} \in [0,1]$ denotes the spectral radius of the amplification matrix at the highest mode \cite{Jansen2000,Chung1993}. Except the simulations presented in Section \ref{subsec:fda_benchmark}, we choose $\varrho_{\infty} = 0.5$ in the generalized-$\alpha$ method.

\begin{remark}
In the literature, the pressure is typically evaluated at the time step $n+1$ rather than $n+\alpha_f$ within the generalized-$\alpha$ scheme\cite{Figueroa2006,Bazilevs2007a,Taylor1998,Moghadam2013}. In a recent investigation, we found that that choice leads to a first-order temporal convergence rate for the pressure \cite{Liu2020a}. Evaluating the pressure at the time step $n+\alpha_f$ not only recovers the second order accuracy but also simplifies the implementation.
\end{remark}

\subsubsection{Temporal discretization of the zero-dimensional problem}
Considering the implicitly coupled outflow boundary condition given by the Windkessel model \eqref{eq:gen_alpha_fully_discrete_linear_momentum}, one needs to evaluate the values of the proximal pressure at time $t_{n+\alpha_f}$. Although an analytic solution  \eqref{eq:RCR_analytic_formula} exists for the specific model considered in this work, it is typically unavailable for general lumped parameter models. Therefore, we consider performing a time integration for the algebraic-differential equations \eqref{eq:RCR_circuit_1}-\eqref{eq:RCR_circuit_2} for its evaluation. Following the framework proposed in \cite{Moghadam2013}, we solve the ordinary differential equations \eqref{eq:RCR_circuit_1} via an explicit fourth-order Runge-Kutta method by subdividing the time interval $(t_n, t_{n+1})$ into $n_{\mathrm{ts}}$ equal-sized subintervals $(t_{n,m}, t_{n,m+1})$ that satisfies 
\begin{align*}
t_{n,m} := t_n + m h, \quad \mbox{ for } m=0, \cdots, n_{\mathrm{ts}}, \mbox{ with } h := \frac{\Delta t_n}{n_{\mathrm{ts}}}.
\end{align*}
We denote by $\Pi^k_{n,m}$ the approximations to $\Pi^k(t_{n,m})$, and we may define the approximation to the flow rate at the intermediate time step $t_{n,m}$ by linear interpolations,
\begin{align*}
Q^k_{n,m}:= \left(1 -  \frac{m}{n_{\mathrm{ts}}} \right) Q^k_{n,0} + \frac{m}{n_{\mathrm{ts}}} Q^k_{n,n_{\mathrm{ts}}} = \left(1 -  \frac{m}{n_{\mathrm{ts}}} \right) Q^k_n + \frac{m}{n_{\mathrm{ts}}} Q^k_{n+1}, \quad \mbox{ for } m=0, \cdots, n_{\mathrm{ts}}.
\end{align*}
In the above, the values of $Q^k_{n}$ and $Q^k_{n+1}$ can be explicitly calculated by 
\begin{align*}
Q^k_{n,0} = Q^k_n := \int_{\Gamma_{\mathrm{out}}^k} \bm v_n \cdot \bm n d\Gamma, \quad Q^k_{n,n_{\mathrm{ts}}} = Q^k_{n+1} := \int_{\Gamma_{\mathrm{out}}^k} \bm v_{n+1} \cdot \bm n d\Gamma.
\end{align*}
Given the values of $\Pi^k_n$, $Q^k_{n}$, and $Q^k_{n+1}$ as the input data, the algorithm for obtaining $P^k_{n+1}$, the approximation to $P^k\left(t_{n+1}\right)$, is stated in Algorithm \ref{algorithm:runge_kutta_4}. To simplify notations, we denote the dependency of $P^k_{n+1}$ on the input data $\Pi^k_{n}$, $Q^k_{n}$, and $Q^k_{n+1}$ through the Algorithm \ref{algorithm:runge_kutta_4} as $P^k_{n+1} = \mathcal G(\Pi^k_{n}, Q^k_{n}, Q^k_{n+1})$ in the following text.
\begin{algorithm}[H]
\caption{The fourth-order Runge-Kutta method for solving \eqref{eq:RCR_circuit_1}-\eqref{eq:RCR_circuit_2}.}
\label{algorithm:runge_kutta_4}
\begin{algorithmic}[1]
\State \texttt{Set $\Pi^k_{n,0} \gets \Pi^k_{n}$, $Q^k_{n,0} \gets Q^k_{n}$, and $Q^k_{n,n_{\mathrm{ts}}} \gets Q^k_{n+1}$}

\For{ \texttt{$m =0$ to $n_{\mathrm{ts}}-1$ } }
	\State \texttt{Calculate $K_1 \gets \mathcal F\left(\Pi^k_{n,m}, Q^k_{n,m}\right)$}
	\State \texttt{Calculate $K_2 \gets \mathcal F\left(\Pi^k_{n,m}+ \frac13 K_1 h, \frac23 Q^k_{n,m} + \frac13 Q^k_{n,m+1} \right)$}
	\State \texttt{Calculate $K_3 \gets \mathcal F\left(\Pi^k_{n,m}- \frac13 K_1 h + K_2 h, \frac13 Q^k_{n,m} + \frac23 Q^k_{n,m+1}\right)$}
	\State \texttt{Calculate $K_4 \gets \mathcal F\left(\Pi^k_{n,m} + K_1 h - K_2 h + K_3 h, Q^k_{n,m+1}\right)$}
	\State \texttt{Calculate $\Pi^k_{n,m+1} \gets \Pi^k_{n,m} + \frac18 K_1 h + \frac38 K_2 h + \frac38 K_3 h + \frac18 K_4 h$}
\EndFor

\State \Return $P^k_{n+1} \gets \mathrm R^k_{\mathrm p} Q^k_{n+1} + \Pi^k_{n,n_{\mathrm{ts}}} + P^k_{d}(t_{n+1})$
\end{algorithmic}
\end{algorithm}

\subsection{Predictor multi-corrector algorithm}
\label{sec:predictor multicorrector algorithm}
The system of equations \eqref{eq:gen_alpha_fully_discrete_linear_momentum}-\eqref{eq:gen_alpha_def_y_n_alpha_f} are solved using the Newton-Raphson method with consistent linearization. At the time step $t_{n+1}$, the solution vector $\bm y_{n+1}$ is solved by the following predictor multi-corrector algorithm. We denote $\bm y_{n+1,(l)} := \left\lbrace \bm v_{n+1,(l)}, p_{n+1,(l)} \right\rbrace^T$ as the solution vector for the three-dimensional problem evaluated at the Newton-Raphson iteration step $l=0,\cdots, l_{\textup{max}}$. The residual vectors evaluated at the iteration stage $l$ are denoted as
\begin{align*}
\boldsymbol{\mathrm R}_{(l)} :=& \left\lbrace \boldsymbol{\mathrm R}_{m,(l)}, \boldsymbol{\mathrm R}_{p,(l)} \right\rbrace^T, 
\end{align*}
with
\begin{align*}
\boldsymbol{\mathrm R}_{m,(l)}^{\textup{vol}} :=& \boldsymbol{\mathrm R}_m^{\textup{vol}}\left( \dot{\bm y}_{n+\alpha_m, (l)}, \bm y_{n+\alpha_f, (l)} \right), &&
\boldsymbol{\mathrm R}_{m,(l)}^{\textup{bc}} := \boldsymbol{\mathrm R}_m^{\textup{bc}}\left( \dot{\bm y}_{n+\alpha_m, (l)}, \bm y_{n+\alpha_f, (l)} \right), \displaybreak[2] \\
\boldsymbol{\mathrm R}_{m,(l)}^{\textup{bf}} :=& \boldsymbol{\mathrm R}_m^{\textup{bf}}\left( \dot{\bm y}_{n+\alpha_m, (l)}, \bm y_{n+\alpha_f, (l)} \right), &&
\boldsymbol{\mathrm R}_{m,(l)} := \boldsymbol{\mathrm R}_m\left( \dot{\bm y}_{n+\alpha_m, (l)}, \bm y_{n+\alpha_f, (l)} \right) = \boldsymbol{\mathrm R}_{m,(l)}^{\textup{vol}} + \boldsymbol{\mathrm R}_{m,(l)}^{\textup{bc}} + \boldsymbol{\mathrm R}_{m,(l)}^{\textup{bf}}, \displaybreak[2] \\
\boldsymbol{\mathrm R}_{p,(l)} :=& \boldsymbol{\mathrm R}_p\left( \dot{\bm y}_{n+\alpha_m, (l)}, \bm y_{n+\alpha_f, (l)} \right). &&
\end{align*}
In the above, the term $\boldsymbol{\mathrm R}_{m,(l)}^{\textup{bc}}$ can be explicitly written as
\begin{align*}
\boldsymbol{\mathrm R}_{m,(l)}^{\textup{bc}} = \sum_{k=1}^{\mathrm n_{\mathrm{out}}} \int_{\Gamma_{\mathrm{out}}^k} P^k_{n+\alpha_f, (l)} N_A n_i d\Gamma = \sum_{k=1}^{\mathrm n_{\mathrm{out}}} \int_{\Gamma_{\mathrm{out}}^k} \left( (1-\alpha_f)P^k_{n} + \alpha_f P^k_{n+1,(l)} \right) N_A n_i d\Gamma.
\end{align*} 
With $\bm v_{n+1,(l)}$, the corresponding flow rate at the $\Gamma_{\mathrm{out}}^k$ surface is denoted by $Q^k_{n+1,(l)}$. We use $P^k_{n+1,(l)}$ to represent the solution of the zero-dimensional problem solved with $\Pi^k_{n}$, $Q^k_{n}$, and $Q^k_{n+1,(l)}$ as the input data for Algorithm \ref{algorithm:runge_kutta_4}. The consistent tangent matrix associated with the above residual vectors is
\begin{align}
\label{eq:def_matrix_K}
\boldsymbol{\mathcal A}_{(l)} =
\begin{bmatrix}
\boldsymbol{\mathrm A}_{(l)} & \boldsymbol{\mathrm B}_{(l)} \\[0.3mm]
\boldsymbol{\mathrm C}_{(l)} & \boldsymbol{\mathrm D}_{(l)}
\end{bmatrix},
\end{align}
wherein
\begin{align}
\label{eq:def_matrix_A}
& \boldsymbol{\mathrm A}_{(l)} = \boldsymbol{\mathrm F}_{(l)} + \boldsymbol{\mathrm K}_{(l)}, \displaybreak[2] \\
& \boldsymbol{\mathrm F}_{(l)} = \alpha_m \frac{\partial \boldsymbol{\mathrm R}_m^{\textup{vol}}\left( \dot{\bm y}_{n+\alpha_m, (l)}, \bm y_{n+\alpha_f, (l)} \right)}{\partial \dot{\bm v}_{n+\alpha_m, (l)}} + \alpha_f \gamma \Delta t_n \frac{\partial \boldsymbol{\mathrm R}_m^{\textup{vol}}\left( \dot{\bm y}_{n+\alpha_m, (l)}, \bm y_{n+\alpha_f, (l)} \right)}{\partial \bm v_{n+\alpha_f, (l)} } \displaybreak[2] \\
& \hspace{10mm} + \alpha_f \gamma \Delta t_n \frac{\partial \boldsymbol{\mathrm R}_m^{\textup{bf}}\left( \dot{\bm y}_{n+\alpha_m, (l)}, \bm y_{n+\alpha_f, (l)} \right)}{\partial \bm v_{n+\alpha_f, (l)} }, \displaybreak[2] \\
& \boldsymbol{\mathrm K}_{(l)} = \alpha_f \gamma \Delta t_n \frac{\partial \boldsymbol{\mathrm R}_m^{\textup{bc}}\left( \dot{\bm v}_{n+\alpha_m, (l)}, \bm y_{n+\alpha_f, (l)} \right)}{\partial \bm v_{n+\alpha_f, (l)} }, \displaybreak[2] \\
& \boldsymbol{\mathrm B}_{(l)} = \alpha_f \gamma \Delta t_n \frac{\partial \boldsymbol{\mathrm R}_m^{\textup{vol}}\left( \dot{\bm y}_{n+\alpha_m, (l)}, \bm y_{n+\alpha_f, (l)} \right)}{\partial \bm p_{n+\alpha_f, (l)} }, \displaybreak[2] \\
& \boldsymbol{\mathrm C}_{(l)} = \alpha_m \frac{\partial \boldsymbol{\mathrm R}_p\left( \dot{\bm y}_{n+\alpha_m, (l)}, \bm y_{n+\alpha_f, (l)} \right)}{\partial \dot{\bm v}_{n+\alpha_m, (l)}} +  \alpha_f \gamma \Delta t_n \frac{\partial \boldsymbol{\mathrm R}_p\left( \dot{\bm y}_{n+\alpha_m, (l)}, \bm y_{n+\alpha_f, (l)} \right)}{\partial \bm v_{n+\alpha_f, (l)} }, \displaybreak[2] \\
\label{eq:def_matrix_D}
& \boldsymbol{\mathrm D}_{(l)} = \alpha_f \gamma \Delta t_n \frac{\partial \boldsymbol{\mathrm R}_p\left( \dot{\bm y}_{n+\alpha_m, (l)}, \bm y_{n+\alpha_f, (l)} \right)}{\partial \bm p_{n+\alpha_f, (l)} }.
\end{align}
Specifically, $\boldsymbol{\mathrm K}_{(l)}$ can be explicitly expressed as
\begin{align}
\boldsymbol{\mathrm K}_{(l)} = \alpha_f \gamma \Delta t_n \sum_{k=1}^{\mathrm n_{\mathrm{out}}} \left( \mathfrak m^k_{(l)}  \bm a^k \bm a^{kT} \right),
\end{align}
wherein $\bm a^k$ is a column vector with its entries defined as 
\begin{align*}
\bm a^k_{Ai} := \left[ \int_{\Gamma_{\mathrm{out}}^k} N_A \bm n_i^k d\Gamma \right],
\end{align*}
and $\mathfrak m^k_{(l)}$ is a scalar variable defined as
\begin{align*}
\mathfrak m^k_{(l)} := \frac{\partial P^k_{n+\alpha_f,(l)}}{\partial Q^k_{n+\alpha_f,(l)}} = \frac{\partial P^k_{n+\alpha_f,(l)}}{\partial P^k_{n+1,(l)}} \frac{\partial P^k_{n+1,(l)}}{\partial Q^k_{n+1,(l)}} \frac{\partial Q^k_{n+1,(l)}}{\partial Q^k_{n+\alpha_f,(l)}} = \alpha_f \frac{\partial P^k_{n+1,(l)}}{\partial Q^k_{n+1,(l)}} \frac{1}{\alpha_f} = \frac{\partial P^k_{n+1,(l)}}{\partial Q^k_{n+1,(l)}}.
\end{align*} 
The matrix $\boldsymbol{\mathrm K}_{(l)}$ is a weighted sum of rank-one matrices $\bm a^k \bm a^{kT}$, with weights being $\alpha_f \gamma \Delta t_n \mathfrak m^k_{(l)}$. When the outflow boundaries are prescribed by the resistance boundary condition \eqref{eq:resistance_bc}, we directly have $\mathfrak m^k_{(l)} = \mathrm R^k$. Yet, for general reduced models, the values of $\mathfrak m^k_{(l)}$ can only be obtained through a difference approximation. One first obtains 
\begin{align*}
\tilde{P}^k_{n+1,(l)} := \mathcal G(\Pi^k_n, Q^k_n, Q^k_{n+1,(l)}+ \epsilon/2) \quad \mbox{ and } \quad \hat{P}^k_{n+1,(l)} := \mathcal G(\Pi^k_n, Q^k_n, Q^k_{n+1,(l)}- \epsilon/2)
\end{align*}
by calling Algorithm \ref{algorithm:runge_kutta_4} and then calculates
\begin{align}
\label{eq:difference_approx_m_k}
m^k_{(l)} \approx \frac{\tilde{P}^k_{n+1,(l)} - \hat{P}^k_{n+1,(l)}}{\epsilon}, \mbox{ with } \epsilon := \textup{max}\left\lbrace \epsilon_{\textup{abs}}, \epsilon_{\textup{rel}}|Q^k_{n+1,(l)}| \right\rbrace.
\end{align}
Following the choice made in the \texttt{svSolver}, we fix $\epsilon_{\textup{abs}} = 10^{-8}$ and $\epsilon_{\textup{rel}} = 10^{-5}$. Based on the above discussion, the predictor multi-corrector algorithm for solving the nonlinear algebraic equations in each time step can be summarized as follows.

\noindent \textbf{Predictor stage}: Set:
\begin{align*}
\bm y_{n+1, (0)} = \bm y_{n}, \quad
\dot{\bm y}_{n+1, (0)} = \frac{\gamma - 1}{\gamma} \dot{\bm y}_{n}.
\end{align*}

\noindent \textbf{Multi-corrector stage}:
Repeat the following steps for \(l = 1, \dots, l_{\textup{max}}\):
\begin{enumerate}
\item Calculate the flow rate 
\begin{align*}
Q^k_{n+1,(l)} := \int_{\Gamma_{\mathrm{out}}^k} \bm v_{n+1,(l)}\cdot \bm n d\Gamma.
\end{align*}

\item Calculate 
\begin{align*}
P^k_{n+1,(l)} &= \mathcal G(\Pi^k_n, Q^k_n, Q^k_{n+1,(l)}), \\ 
\tilde{P}^k_{n+1,(l)} &= \mathcal G(\Pi^k_n, Q^k_n, Q^k_{n+1,(l)}+\frac12 \epsilon), \\ \hat{P}^k_{n+1,(l)} &= \mathcal G(\Pi^k_n, Q^k_n, Q^k_{n+1,(l)}-\frac12 \epsilon),
\end{align*}
by invoking Algorithm \ref{algorithm:runge_kutta_4}.

\item Calculate $m^k_{(l)}$ from the difference quotient \eqref{eq:difference_approx_m_k}.

\item Evaluate the solution vectors at the intermediate stages:
\begin{align*}
\bm y_{n+\alpha_f, (l)} = \bm y_n + \alpha_f \left( \bm y_{n+1,(l-1)} - \bm y_n \right), \quad
\dot{\bm y}_{n+\alpha_m, (l)} = \dot{\bm y}_n + \alpha_m \left( \dot{\bm y}_{n+1,(l-1)} - \dot{\bm y}_n \right).
\end{align*}

\item Evaluate $P^k_{n+\alpha_f,(l)} = (1-\alpha_f)P^k_{n} + \alpha_f P^k_{n+1,(l)}$.

\item Assemble the residual vectors $\boldsymbol{\mathrm R}_{m,(l)}$ and $\boldsymbol{\mathrm R}_{p,(l)}$ using $\bm y_{n+\alpha_f, (l)}$, $\dot{\bm y}_{n+\alpha_m, (l)}$, and $P^k_{n+\alpha_f, (l)}$.

\item Let $\|\boldsymbol{\mathrm R}_{(l)}\|_{\mathfrak{l}^2}$ denote the $\mathfrak l^2$-norm of the residual vector. If either one of the following stopping criteria 
\begin{align*}
& \frac{\|\boldsymbol{\mathrm R}_{(l)}\|_{\mathfrak l^2}}{\|\boldsymbol{\mathrm R}_{(0)}\|_{\mathfrak l^2}} \leq \textup{tol}_{\textup{R}}, \qquad \|\boldsymbol{\mathrm R}_{(l)}\|_{\mathfrak l^2} \leq \textup{tol}_{\textup{A}},
\end{align*}
is satisfied for two prescribed tolerances $\textup{tol}_{\textup{R}}$, $\textup{tol}_{\textup{A}}$, set the solution vector at time step $t_{n+1}$ as $\dot{\bm y}_{n+1} = \dot{\bm y}_{n+1, (l-1)}$ and $\bm y_{n+1} = \bm y_{n+1, (l-1)}$, and exit the multi-corrector stage; otherwise, continue to step 8.

\item Assemble the tangent matrices \eqref{eq:def_matrix_K}-\eqref{eq:def_matrix_D}.

\item Solve the following linear system of equations for $\Delta \dot{p}_{n+1,(l)}$ and $\Delta \dot{\bm v}_{n+1,(l)}$,
\begin{align}
\label{eq:pre-multi-corrector-stage-5-matrix-problems}
\begin{bmatrix}
\boldsymbol{\mathrm A}_{(l)} & \boldsymbol{\mathrm B}_{(l)} \\[0.3mm]
\boldsymbol{\mathrm C}_{(l)} & \boldsymbol{\mathrm D}_{(l)}
\end{bmatrix}
\begin{bmatrix}
\Delta \dot{\bm v}_{n+1,(l)} \\[0.3mm]
\Delta \dot{p}_{n+1,(l)}
\end{bmatrix}
= -
\begin{bmatrix}
\boldsymbol{\mathrm R}_{m,(l)} \\[0.3mm]
\boldsymbol{\mathrm R}_{p,(l)}
\end{bmatrix}.
\end{align}

\item Let $\Delta \dot{\bm y}_{n+1,(l)} := \left\lbrace \Delta \dot{\bm v}_{n+1,(l)}, \Delta \dot{p}_{n+1,(l)} \right\rbrace^T$. Update the solution vectors as
\begin{align*}
\bm y_{n+1,(l)} = \bm y_{n+1,(l)} + \gamma \Delta t_n \Delta \dot{\bm y}_{n+1,(l)}, \quad
\dot{\bm y}_{n+1,(l)} = \dot{\bm y}_{n+1,(l)} + \Delta \dot{\bm y}_{n+1,(l)},
\end{align*}
and return to step 1.

\end{enumerate}
For the numerical simulations presented in this work, unless otherwise specified, we adopt the tolerances for the nonlinear iteration as $\textup{tol}_{\textup{R}} = \textup{tol}_{\textup{A}} = 10^{-6}$ and the maximum number of iterations as $l_{\textup{max}}=20$. Notice that in step 5, we evaluated the pressure from the reduced model at the intermediate time step $n+\alpha_f$ to make it consistent with the generalized-$\alpha$ scheme given in Sec. \ref{sec:temporal_scheme_3d_problem}. This is different from the approach adopted in prior studies \cite{Moghadam2013}, in which the pressure from the reduced model and the pressure in the three-dimensional problem are both evaluated at the time step $n+1$. 

\section{Iterative solution method with the nested block preconditioning}
In this section, we focus on developing the iterative solution method for \eqref{eq:pre-multi-corrector-stage-5-matrix-problems}, which is the most time consuming part in the predictor multi-corrector algorithm. For notational simplicity, we neglect the subscript $(l)$ used to represent the iteration steps in the multi-corrector stage. In this section, we consider the linear system $\boldsymbol{\mathcal A} \bm x = \bm r$, where the matrix and vectors adopt the following block structure,
\begin{align*}
\boldsymbol{\mathcal A} :=
\begin{bmatrix}
\boldsymbol{\mathrm A} & \boldsymbol{\mathrm B} \\[0.3mm]
\boldsymbol{\mathrm C} & \boldsymbol{\mathrm D}
\end{bmatrix}, \quad
\bm x :=
\begin{bmatrix}
\bm x_{\bm v} \\[0.3em]
\bm x_{p}
\end{bmatrix},
\quad
\bm r := 
\begin{bmatrix}
\bm r_{\bm v} \\[0.3em]
\bm r_{p}
\end{bmatrix}.
\end{align*}
For a sparse non-symmetric matrix like $\boldsymbol{\mathcal A}$, the GMRES algorithm and its variants \cite{Saad1993,Saad1986}, with a suitable preconditioning technique, are among the most effective general-purpose solution method. In this work, we consider using the FGMRES \cite{Saad1993} as the iterative solution method, since the proposed preconditioner varies over iterations. Let $\boldsymbol{\mathcal P}_{(i)}$ denote the preconditioning matrix in the $i$-th iteration. In the FGMRES algorithm with right preconditioning, one strives to generate a subspace by repeatedly applying $\boldsymbol{\mathcal A} \boldsymbol{\mathcal P}_{(i)}^{-1}$ to a vector, and search for an approximate solution by minimizing the residual over the subspace. This procedure is referred to as the \textit{outer solver} in this work. Similar to the GMRES algorithm, the FGMRES algorithm is restarted every $\mathrm m$ steps and is denoted as FGMRES($\mathrm m$). It is worth mentioning that there are also caveats associated with the FGMRES algorithm. It requires slightly more memory to store the generated subspace. More critically, there is no general convergence theory for it, since it does not generate a standard Krylov subspace. In practice, one might observe stagnation or even divergence due to an improper choice of the preconditioner. Therefore, a suitable design of $\boldsymbol{\mathcal P}_{(i)}$ has become the crux. With the understanding that the preconditioners may vary over the iterations, we neglect the subscript $(i)$ for notational simplicity. In Sec. \ref{sec:SCR}, we focus on designing a preconditioning technique based on the SCR procedure, denoted as $\boldsymbol{\mathcal P}_{\textup{SCR}}$. Following that, we will review the SIMPLE preconditioner $\boldsymbol{\mathcal P}_{\textup{SIMPLE}}$ in Sec. 
\ref{sec:simple_pc} and the iterative algorithm used in the current \texttt{svSolver} \cite{simvascular-simulation-guide} in Sec. \ref{subsec:bipn}. 

\subsection{Schur complement reduction}
\label{sec:SCR}
In this section, we introduce the concept of SCR and design a preconditioning technique based on it. The matrix $\boldsymbol{\mathcal A}$ can be factorized as follows,
\begin{align}
\label{eq:a_ldu_block_factorization}
\boldsymbol{\mathcal A} = \boldsymbol{\mathcal L} \boldsymbol{\mathcal D} \boldsymbol{\mathcal U} =
\begin{bmatrix}
\boldsymbol{\mathrm I} & \boldsymbol{\mathrm O} \\[0.3em]
\boldsymbol{\mathrm C} \boldsymbol{\mathrm A}^{-1} & \boldsymbol{\mathrm I}
\end{bmatrix}
\begin{bmatrix}
\boldsymbol{\mathrm A} & \boldsymbol{\mathrm O} \\[0.3em]
\boldsymbol{\mathrm O} & \boldsymbol{\mathrm S}
\end{bmatrix}
\begin{bmatrix}
\boldsymbol{\mathrm I} & \boldsymbol{\mathrm A}^{-1} \boldsymbol{\mathrm B} \\[0.3em]
\boldsymbol{\mathrm O} & \boldsymbol{\mathrm I}
\end{bmatrix},
\end{align}
with $\boldsymbol{\mathrm S} := \boldsymbol{\mathrm D} - \boldsymbol{\mathrm C} \boldsymbol{\mathrm A}^{-1} \boldsymbol{\mathrm B}$ being the Schur complement of $\boldsymbol{\mathrm A}$. This implies that we may solve $\boldsymbol{\mathcal A} \bm x = \bm r$ by considering
\begin{align*}
\boldsymbol{\mathcal D} \boldsymbol{\mathcal U} \bm x =
\begin{bmatrix}
\boldsymbol{\mathrm A} & \boldsymbol{\mathrm B} \\[0.3em]
\boldsymbol{\mathrm O} & \boldsymbol{\mathrm S}
\end{bmatrix}
\begin{bmatrix}
\bm x_{\bm v} \\[0.3em]
\bm x_{p}
\end{bmatrix}
= \boldsymbol{\mathcal L}^{-1} \bm r =
\begin{bmatrix}
\boldsymbol{\mathrm I} & \boldsymbol{\mathrm O} \\[0.3em]
\boldsymbol{\mathrm C} \boldsymbol{\mathrm A}^{-1} & \boldsymbol{\mathrm I}
\end{bmatrix}^{-1}
\begin{bmatrix}
\bm r_{\bm v} \\[0.3em]
\bm r_{p}
\end{bmatrix}
=
\begin{bmatrix}
\boldsymbol{\mathrm I} & \boldsymbol{\mathrm O} \\[0.3em]
-\boldsymbol{\mathrm C} \boldsymbol{\mathrm A}^{-1} & \boldsymbol{\mathrm I}
\end{bmatrix}
\begin{bmatrix}
\bm r_{\bm v} \\[0.3em]
\bm r_{p}
\end{bmatrix} =
\begin{bmatrix}
\bm r_{\bm v} \\[0.3em]
\bm r_{p} - \boldsymbol{\mathrm C} \boldsymbol{\mathrm A}^{-1} \bm r_{\bm v}
\end{bmatrix}.
\end{align*}
The above system of equations can then be solved by back substitution, which results in the SCR procedure \cite{Benzi2005,May2008}. Therefore, the design of the solution method for the matrix problem $\boldsymbol{\mathcal A}$ reduces to the design of solution methods for matrices $\boldsymbol{\mathrm A}$ and $\boldsymbol{\mathrm S}$. Since the algebraic form of the matrix $\boldsymbol{\mathrm A}$ is available, one may apply GMRES($\mathrm m_{\mathrm A}$) with a suitable preconditioner $\boldsymbol{\mathrm P}_{\mathrm A}$ to solve it. In this work, we consider two options for the preconditioner. We use AMG for robust and scalable performances, and we use the Jacobi preconditioner to achieve efficient performances. The stopping criteria for solving with $\boldsymbol{\mathrm A}$ include the tolerance for the relative error $\delta^r_{\mathrm A}$, the tolerance for the absolute error $\delta^a_{\mathrm A}$, and the maximum number of iterations $\mathrm n^{\textup{max}}_{\mathrm A}$. 

The difficulty of the SCR procedure comes from solving the equations associated with $\boldsymbol{\mathrm S}$ because the Schur complement $\boldsymbol{\mathrm S}$ is defined as a composition of all four block matrices. The appearance of $\boldsymbol{\mathrm A}^{-1}$ inevitably makes $\boldsymbol{\mathrm S}$ a dense matrix. Therefore, forming $\boldsymbol{\mathrm S}$ algebraically is impractical primarily due to the limit of computer memory sizes. Recall that for all iterative methods, one just needs to repeatedly perform matrix-vector multiplications. Therefore, the application of the matrix $\boldsymbol{\mathrm S}$ on a vector can be achieved without an algebraic definition of $\boldsymbol{\mathrm S}$, and this procedure is summarized as the matrix-free algorithm in Algorithm \ref{algorithm:matrix_free_mat_vec_for_S}. Notice that in this algorithm one has to solve an equation associated with $\boldsymbol{\mathrm A}$.  This procedure is referred to as the \textit{inner solver}, since it resides inside the solution procedure for $\boldsymbol{\mathrm S}$. In practice, we invoke the same iterative method and the same preconditioning technique for equations associated with $\boldsymbol{\mathrm A}$ in the inner solver with potentially different tolerances for the relative error $\delta^r_{\mathrm I}$. One may call a GMRES($\mathrm m_{\mathrm S}$) algorithm to solve with $\boldsymbol{\mathrm S}$ using the Algorithm \ref{algorithm:matrix_free_mat_vec_for_S} to define the action of the matrix on a vector. We may further accelerate the algorithm by providing a preconditioner for this GMRES iteration. A sparse matrix $\hat{\boldsymbol{\mathrm S}}:= \boldsymbol{\mathrm D} - \boldsymbol{\mathrm C} \left(\textup{diag}\left(\boldsymbol{\mathrm A}\right)\right)^{-1} \boldsymbol{\mathrm B}$ can be constructed as an approximation to $\boldsymbol{\mathrm S}$, and then an AMG preconditioner based on $\hat{\boldsymbol{\mathrm S}}$ can be constructed and applied as a left preconditioner in the GMRES iterations. This preconditioner is denoted by $\boldsymbol{\mathrm P}_{\mathrm S}$. The stopping criteria for solving with $\boldsymbol{\mathrm S}$ include the tolerance for the relative error $\delta^r_{\mathrm S}$, the tolerance for the absolute error $\delta^a_{\mathrm S}$, and the maximum number of iterations $\mathrm n^{\textup{max}}_{\mathrm S}$.

\begin{algorithm}[H]
\caption{The matrix-free algorithm for the multiplication of $\boldsymbol{\mathrm S}$ with a vector $\bm x_p$.}
\label{algorithm:matrix_free_mat_vec_for_S}
\begin{algorithmic}[1]
\State \texttt{Compute the matrix-vector multiplication $\hat{\bm x}_p \gets \boldsymbol{\mathrm D} \bm x_p$.}
\State \texttt{Compute the matrix-vector multiplication $\bar{\bm x}_p \gets \boldsymbol{\mathrm B} \bm x_p$.}
\State \texttt{Solve for $\tilde{\bm x}_p$ from the linear system} 
\begin{align}
\label{eq:S_inner_A_eqn}
\boldsymbol{\mathrm A} \tilde{\bm x}_p = \bar{\bm x}_p
\end{align}
\texttt{by GMRES($\mathrm m_{\mathrm A}$) preconditioned by $\boldsymbol{\mathrm P}_{\mathrm A}$ with $\delta^r_{\mathrm I}$, $\delta^a_{\mathrm A}$, and $\mathrm n^{\textup{max}}_{\mathrm A}$ prescribed.}
\State \texttt{Compute the matrix-vector multiplication $\bar{\bm x}_p \gets \boldsymbol{\mathrm C} \tilde{\bm x}_p$.}
\State \Return $\hat{\bm x}_p - \bar{\bm x}_p$.
\end{algorithmic}
\end{algorithm}

With the solution method for the matrices $\boldsymbol{\mathrm A}$ and $\boldsymbol{\mathrm S}$ defined, we may define the preconditioner $\boldsymbol{\mathcal P}_{\textup{SCR}}$. Given the solution method GMRES($\mathrm m_{\mathrm A}$), GMRES($\mathrm m_{\mathrm S}$), the preconditioner $\boldsymbol{\mathrm P}_{\mathrm A}$ , $\boldsymbol{\mathrm P}_{\mathrm S}$, the relative tolerances $\delta^r_{\mathrm A}$, $\delta^r_{\mathrm S}$, and $\delta^r_{\mathrm I}$, the absolute tolerances $\delta^a_{\mathrm A}$, $\delta^a_{\mathrm S}$, and the maximum number of iterations $\mathrm n^{\textup{max}}_{\mathrm A}$, $\mathrm n^{\textup{max}}_{\mathrm S}$, the action of $\boldsymbol{\mathcal P}_{\textup{SCR}}^{-1}$ on a vector can be summarized as the Algorithm \ref{algorithm:exact_block_factorization} in below.
\begin{algorithm}[H]
\caption{The action of $\boldsymbol{\mathcal P}_{\textup{SCR}}^{-1}$ on a vector $\bm s := [ \bm s_{\bm v}; \bm s_p]^T$ with the output being $\bm y := [ \bm y_{\bm v}; \bm y_p]^T$.}
\label{algorithm:exact_block_factorization}
\begin{algorithmic}[1]
\State \texttt{Solve for an intermediate velocity $\hat{\bm y}_{\bm v}$ from the equation}
\begin{align}
\label{eq:seg_sol_int_disp}
\boldsymbol{\mathrm A} \hat{\bm y}_{\bm v} = \bm s_{\bm v}
\end{align}
\texttt{by GMRES($\mathrm m_{\mathrm A}$) preconditioned by $\boldsymbol{\mathrm P}_{\mathrm A}$ with $\delta^r_{\mathrm A}$, $\delta^a_{\mathrm A}$, and $\mathrm n^{\textup{max}}_{\mathrm A}$ prescribed.}
\State \texttt{Update the continuity residual by $\bm s_{p} \gets  \bm s_{p} - \boldsymbol{\mathrm C} \hat{\bm y}_{\bm v}$.}
\State \texttt{Solve for $\bm y_p$ from the equation}
\begin{align}
\label{eq:seg_sol_pres}
\boldsymbol{\mathrm S} \bm y_p = \bm s_{p}
\end{align}
\texttt{by Algorithm \ref{algorithm:matrix_free_mat_vec_for_S} and GMRES($\mathrm m_{\mathrm S}$) preconditioned by $\boldsymbol{\mathrm P}_{\mathrm S}$ with $\delta^r_{\mathrm S}$, $\delta^a_{\mathrm S}$, and $\mathrm n^{\textup{max}}_{\mathrm S}$ prescribed.}
\State \texttt{Update the momentum residual by $\bm s_{\bm v} \gets \bm s_{\bm v} - \boldsymbol{\mathrm B} \bm y_{p}$.}
\State \texttt{Solve for $\bm y_{\bm v}$ from the equation}
\begin{align}
\label{eq:seg_sol_disp}
\boldsymbol{\mathrm A} \bm y_{\bm v} = \bm s_{\bm v}
\end{align}
\texttt{by GMRES($\mathrm m_{\mathrm A}$) preconditioned by $\boldsymbol{\mathrm P}_{\mathrm A}$ with $\delta^r_{\mathrm A}$, $\delta^a_{\mathrm A}$, and $\mathrm n^{\textup{max}}_{\mathrm A}$ prescribed.}
\end{algorithmic}
\end{algorithm}
The three solution procedures for \eqref{eq:seg_sol_int_disp}-\eqref{eq:seg_sol_disp} are referred to as \textit{intermediate solvers}. The preconditioner $\boldsymbol{\mathcal P}_{\textup{SCR}}$ is implicitly defined through the intermediate solvers, and therefore it may vary over iterations. In this work, the absolute tolerances are set to be $\delta^a_{\mathrm A} = \delta^a_{\mathrm S} = 10^{-50}$. When the preconditioners $\boldsymbol{\mathrm P}_{\mathrm A}$ and $\boldsymbol{\mathrm P}_{\mathrm S}$ are constructed by AMG, we adopt the implementation of the BoomerAMG preconditioner from the Hypre package, with the setting summarized in Table \ref{table:amg_settings}. Unless otherwise specified, we choose $\mathrm m_{\mathrm A} = \mathrm m_{\mathrm S}= \mathrm n^{\textup{max}}_{\mathrm A} = \mathrm n^{\textup{max}}_{\mathrm S} = 200$.

\begin{table}[htbp]
\begin{center}
\renewcommand{\arraystretch}{1.0}
\begin{tabular}{p{9.5cm} p{5.5cm} }
\hline
Cycle type &  V-cycle \\
Coarsening method & HMIS \\
Interpolation method & Extended method (ext+i) \\
Truncation factor for the interpolation & $0.3$ \\
Threshold for being strongly connected & $0.5$ \\
Maximum number of elements per row for interp. & $5$ \\
The number of levels for aggressive coarsening & $2$ \\
\hline
\end{tabular}
\end{center}
\caption{Settings of the BoomerAMG preconditioner \cite{Falgout2002}. These choices are made based on the balance of robustness and efficiency. Interested readers are referred to \cite{hypre-boomeramg-moose-intro} for a discussion of these options.}
\label{table:amg_settings}
\end{table}

On the outer level, we invoke the FGMRES($\mathrm m$) to minimize the residual for $\boldsymbol{\mathcal A} \bm x = \bm r$, and, unless otherwise specified, we fix $\mathrm m$ as $200$, the absolute tolerance $\delta^a$ as $10^{-50}$, and the maximum number of iterations $\mathrm n^{\textup{max}}$ as $200$. Then the stopping criterion is completely determined by the relative tolerance $\delta^r$.

\begin{remark}
When using higher-order elements, it is often economical to avoid explicitly constructing the tangent matrix. Instead, the action of the tangent matrix on a vector can be achieved through a matrix-free manner using, e.g., the partial assembly approach \cite{Anderson2018} or the Newton-Raphson-Krylov technique \cite{Knoll2004}. In the meantime, a matrix can be created using low-order elements on the same mesh to construct preconditioners. Although this work focuses on the low-order finite element discretization of the incompressible Navier-Stokes equations, it may serve as a competitive candidate for higher-order methods \cite{Pernice2001}. 
\end{remark}

\subsection{The SIMPLE preconditioner}
\label{sec:simple_pc}
The SIMPLE preconditioner is an algebraic analogue of the Semi-Implicit Method for Pressure Linked Equations (SIMPLE) \cite{Patankar1983}. Together with its variants, it is among the most popular choices for CFD \cite{Cyr2012,Deparis2014}, FSI \cite{Deparis2016}, and other multiphysics problems \cite{White2011,Verdugo2016}. It is partly inspired from the block factorization \eqref{eq:a_ldu_block_factorization} and can be represented as follows,
\begin{align*}
\boldsymbol{\mathcal P}_{\textup{SIMPLE}} := 
\begin{bmatrix}
\boldsymbol{\mathrm I} & \boldsymbol{\mathrm O} \\[0.3em]
\boldsymbol{\mathrm C} \boldsymbol{\mathrm A}^{-1} & \boldsymbol{\mathrm I}
\end{bmatrix}
\begin{bmatrix}
\boldsymbol{\mathrm A} & \boldsymbol{\mathrm O} \\[0.3em]
\boldsymbol{\mathrm O} & \hat{ \boldsymbol{\mathrm S} }
\end{bmatrix}
\begin{bmatrix}
\boldsymbol{\mathrm I} & \left(\textup{diag}\left(\boldsymbol{\mathrm A}\right)\right)^{-1} \boldsymbol{\mathrm B} \\[0.3em]
\boldsymbol{\mathrm O} & \boldsymbol{\mathrm I}
\end{bmatrix}
=
\begin{bmatrix}
\boldsymbol{\mathrm A} & \boldsymbol{\mathrm A}\textup{diag}\left(\boldsymbol{\mathrm A}\right)^{-1}\boldsymbol{\mathrm B} \\[0.3em]
\boldsymbol{\mathrm C} & \boldsymbol{\mathrm D}
\end{bmatrix}.
\end{align*}
Remember that $\hat{\boldsymbol{\mathrm S}}:= \boldsymbol{\mathrm D} - \boldsymbol{\mathrm C} \left(\textup{diag}\left(\boldsymbol{\mathrm A}\right)\right)^{-1} \boldsymbol{\mathrm B}$ is a sparse approximation of $\boldsymbol{\mathrm S}$. The purpose of choosing $\left(\textup{diag}\left(\boldsymbol{\mathrm A}\right)\right)^{-1}$ in the upper diagonal matrix is to ensure that the mass equation is not perturbed \cite{Elman2008,Quarteroni2000}. The action of $\boldsymbol{\mathcal P}_{\textup{SIMPLE}}^{-1}$ on a vector can be implemented in a way similar to Algorithm \ref{algorithm:exact_block_factorization} with modifications to the steps 3 and 4. Roughly speaking, one may view the SIMPLE preconditioner as a simplification of the SCR approach without invoking the inner solver. Similar to the SCR preconditoner, one may invoke a Krylov subspace method preconditioned with multigrid \cite{Elman2008,Cyr2012} or domain decomposition preconditioners \cite{Deparis2016} on the intermediate level to define the action of $\boldsymbol{\mathcal P}_{\textup{SIMPLE}}^{-1}$. Since the actual algebraic definition of the SIMPLE preconditioner also varies over successive iterations, one still has to use the FGMRES method as the outer solver. In this work, the settings of the outer and intermediate solvers in $\boldsymbol{\mathcal P}_{\textup{SIMPLE}}$ always follow the settings used in $\boldsymbol{\mathcal P}_{\textup{SCR}}$ in each specific example. 

\subsection{The bi-partitioned iterative algorithm}
\label{subsec:bipn}
In the current \texttt{svSolver}, the linear system of equations is solved by the BIPN algorithm \cite{Esmaily-Moghadam2015} used in conjunction with a specifically designed preconditioner for the resistance-type boundary condition  \cite{Esmaily-Moghadam2013}. The action of this preconditioner on a vector is achieved through a five-step procedure similar to the one given in Algorithm \ref{algorithm:exact_block_factorization}, with the following differences. First, in step 3, the following equation is solved by an unpreconditioned conjugate gradient method,
\begin{align*}
& \tilde{\boldsymbol{\mathrm S}} \bm x_p = \bm r_p, \\
& \tilde{\boldsymbol{\mathrm S}} := \boldsymbol{\mathrm D} - \boldsymbol{\mathrm B}^T \left[ \left(\textup{diag}\left(\boldsymbol{\mathrm A}\right)\right)^{-1}  - \sum_{k=1}^{\mathrm n_{\mathrm{out}}} \alpha_f \gamma \Delta t_n m^k_{(l)} \frac{\bm b^k \bm b^{kT}}{1 + \alpha_f \gamma \Delta t_n m^k_{(l)}\|\bm b^k\|^2_{\mathfrak{l}^2}} \right] \boldsymbol{\mathrm B},\\
& \bm b^k := \left(\textup{diag}\left(\boldsymbol{\mathrm A}\right)\right)^{-1} \bm a^k.
\end{align*}
The definition of $\tilde{\boldsymbol{\mathrm S}}$ is inspired from the Sherman-Morrison formula for inverting a matrix with rank-one modifications \cite{Akguen2001}. Second, on the outer level, the \texttt{svSolver} solves a least-square problem to minimize the residual by inverting a normal matrix. 

It has been demonstrated that this preconditioner is an effective candidate for hemodynamic simulations since the contributions from the outflow boundary conditions are taken into account. However, this method is not without shortcomings. Firstly, the equation associated with $\tilde{\boldsymbol{\mathrm S}}$ is solved without a preconditioner because the rank-one matrices in $\tilde{\boldsymbol{\mathrm S}}$ are not explicitly constructed. The rank-one matrix is a dense block matrix defined on each outlet surface, necessitating a non-trivial sparsity pattern in the matrix allocation stage. It is quite inconvenient to assemble this matrix in a standard finite element code. However, the action of a rank-one matrix on a vector can be efficiently calculated as a vector inner product in the parallel setting \cite[p.~3547]{Bazilevs2009a}, which leads to a matrix-free definition of the matrices $\boldsymbol{\mathrm A}$ and $\tilde{\boldsymbol{\mathrm S}}$. Using a matrix-free algorithm can be convenient and efficient for small- or medium-sized problems. However, without using a preconditioner, we may reasonably expect a loss of scalability in the performance of \texttt{svSolver} as the problem size gets larger. Secondly, similar to $\boldsymbol{\mathcal P}_{\textup{SIMPLE}}$, the convection information is still missing in the definition of this preconditioner, rendering the method ineffective for strong convection problems. Thirdly, the usage of a normal matrix to minimize the residual will inevitably lead to a nearly singular problem with an increased size of the subspace (see Fig. \ref{fig:conv_history_cyl}). Thereby the outer solver is typically limited to at most ten iterations in \texttt{svSolver} \cite{simvascular-simulation-guide}. This signifies a numerical robustness issue and poses a limit on the accuracy that the current \texttt{svSolver} can attain. Lastly, the \texttt{svSolver} underperforms for FSI problems using the coupled momentum method. The reason is apparently that the wall stiffness matrix is not reflected in $\tilde{\boldsymbol{\mathrm S}}$.  Indeed, for more complicated problems, the definition of $\boldsymbol{\mathrm A}$ may involve contributions from different physics. This motivates us to consider using the matrix-free algorithm (i.e. Algorithm \ref{algorithm:matrix_free_mat_vec_for_S}) to attain a better representation of the Schur complement.

\section{Numerical Results}
\label{sec:results}
Numerical simulations reported in this work are all performed on the Stampede2 supercomputer at Texas Advanced Computing Center (TACC), using the Intel Xeon Platinum 8160 ``Skylake" node. Each compute node contains $48$ processors, with $2.1$ GHz nominal clock rate, and $192$ GB RAM. They are interconnected by a $100$ Gb/s Intel Omni-Path network. More technical details about the Stampede2 supercomputer can be found in the user guide \cite{Stampede2-user-guide}.

In this work, we are primarily concerned with blood flow simulations. The fluid density is fixed to be $1.065$ g/cm$^3$ and the dynamic viscosity is fixed to be $0.035$ poise for all simulations. On the inlet surface $\Gamma_{\mathrm{in}}$, we always apply a parabolic velocity profile for $\bm v_{\mathrm{in}}$ based on a prescribed flow rate $\mathcal Q$.

\subsection{FDA idealized medical device benchmark}
\label{subsec:fda_benchmark}
\begin{figure}
	\begin{center}
	\begin{tabular}{c}
\includegraphics[angle=0, trim=0 100 0 100, clip=true, scale = 0.3]{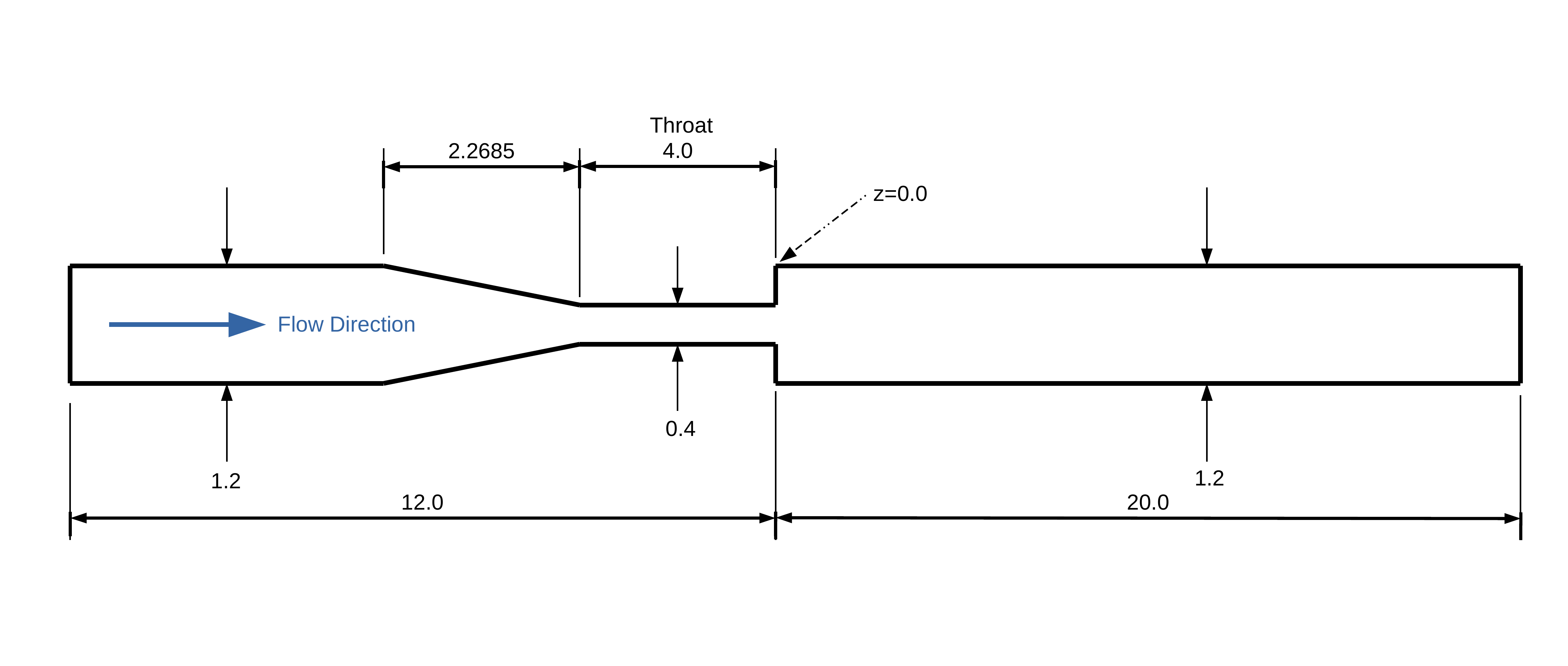}
\end{tabular}
\end{center}
\caption{Schematic representation of the cross section of the computational domain for the idealized medical device problem.} 
\label{fig:fda_benchmark_setting}
\end{figure}

In the first numerical example, we consider a benchmark designed by the US Food and Drug Administration (FDA) to examine the CFD code performance and validate the code by comparing with in vitro experiments. The geometry of the idealized medical device is illustrated in Figure \ref{fig:fda_benchmark_setting}. We follow the practice adopted in \cite{Bergersen2019}, in which the computational domain is chosen to be $32$ cm long. Based on a parabolic velocity profile, we may calculate the maximum velocity on the inlet $\bm v^{\textup{i}}_{\textup{max}}$ as well as at the throat region $\bm v^{\textup{t}}_{\textup{max}}$ from a given value of $\mathcal Q$. Zero traction boundary condition is applied on the outlet surface. The element size $\Delta x$ here is defined as the circumscribing sphere's diameter. In this example, we choose $\varrho_{\infty} = 1.0$ in the generalized-$\alpha$ method to mitigate the numerical dissipation effect, which is critical in this example. With this choice, the time-stepping method reduces to the mid-point rule. The Courant number $\textup{Cr}$ reported in this section is calculated based on the minimum mesh size $\Delta x_{\textup{min}}$ and $\bm v^{\textup{t}}_{\textup{max}}$.

\subsubsection{Solver validation}
\begin{figure}
	\begin{center}
	\begin{tabular}{c}
\includegraphics[angle=0, trim=0 520 0 420, clip=true, scale = 0.2]{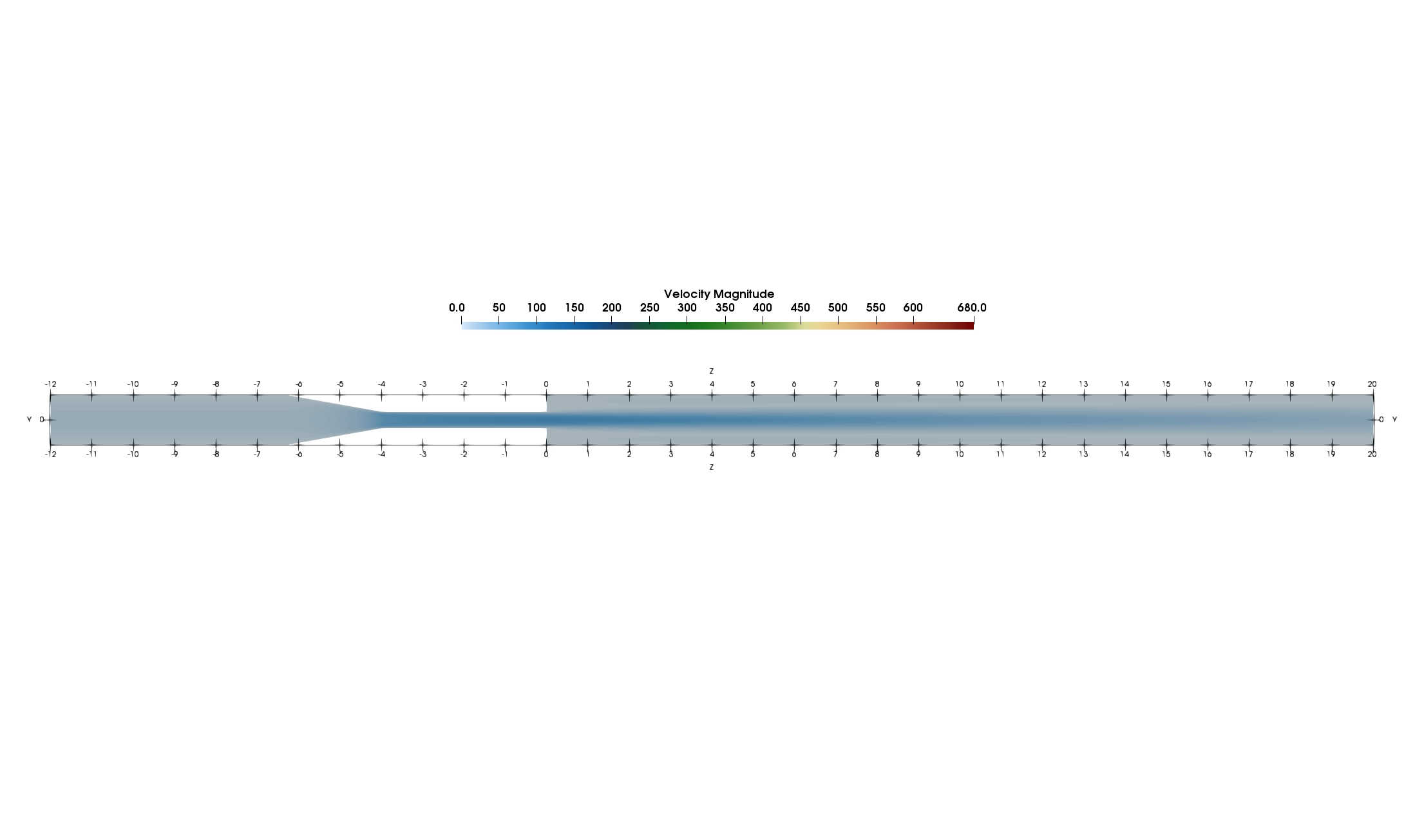} \\
\includegraphics[angle=0, trim=0 520 0 530, clip=true, scale = 0.2]{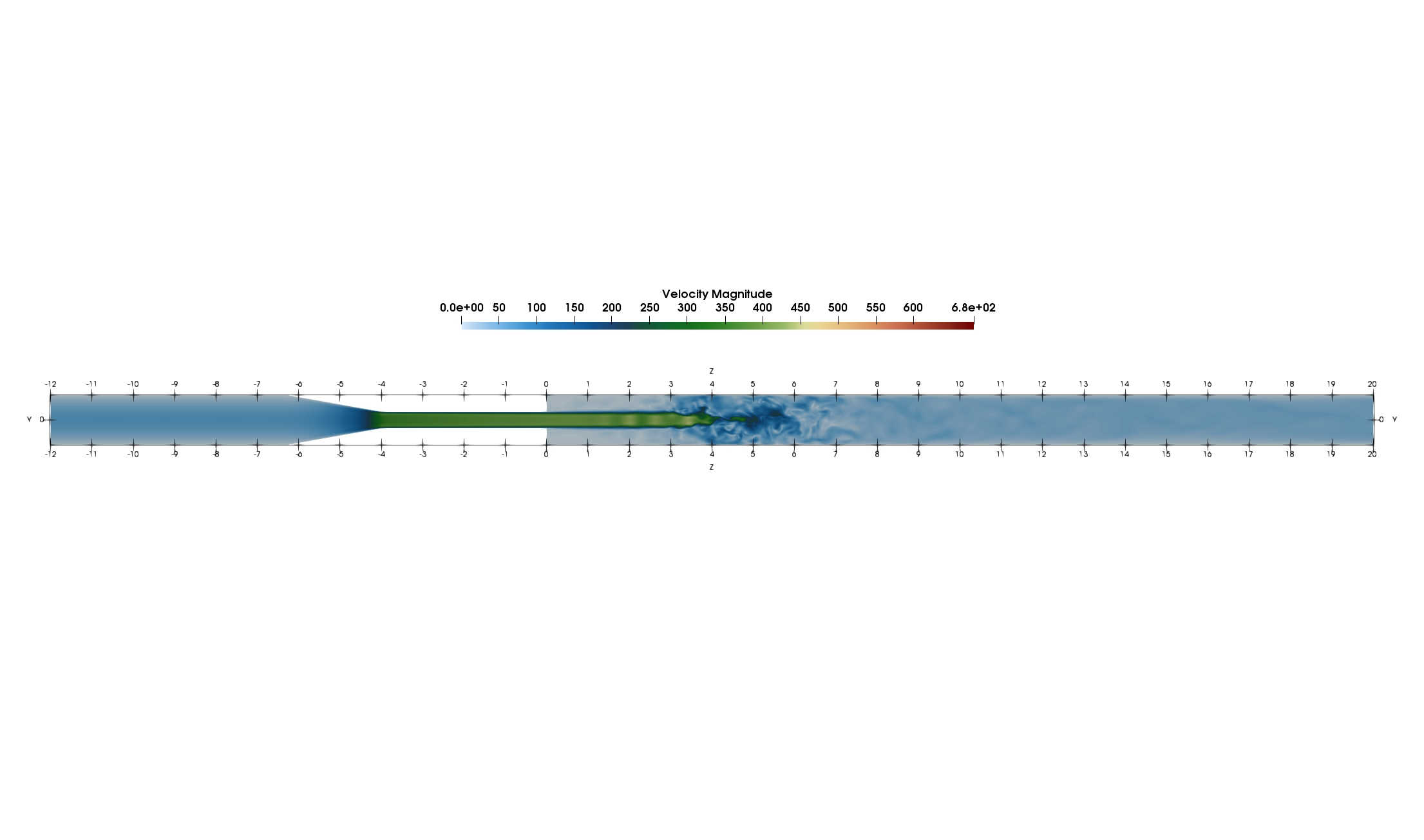}\\
\includegraphics[angle=0, trim=0 520 0 530, clip=true, scale = 0.2]{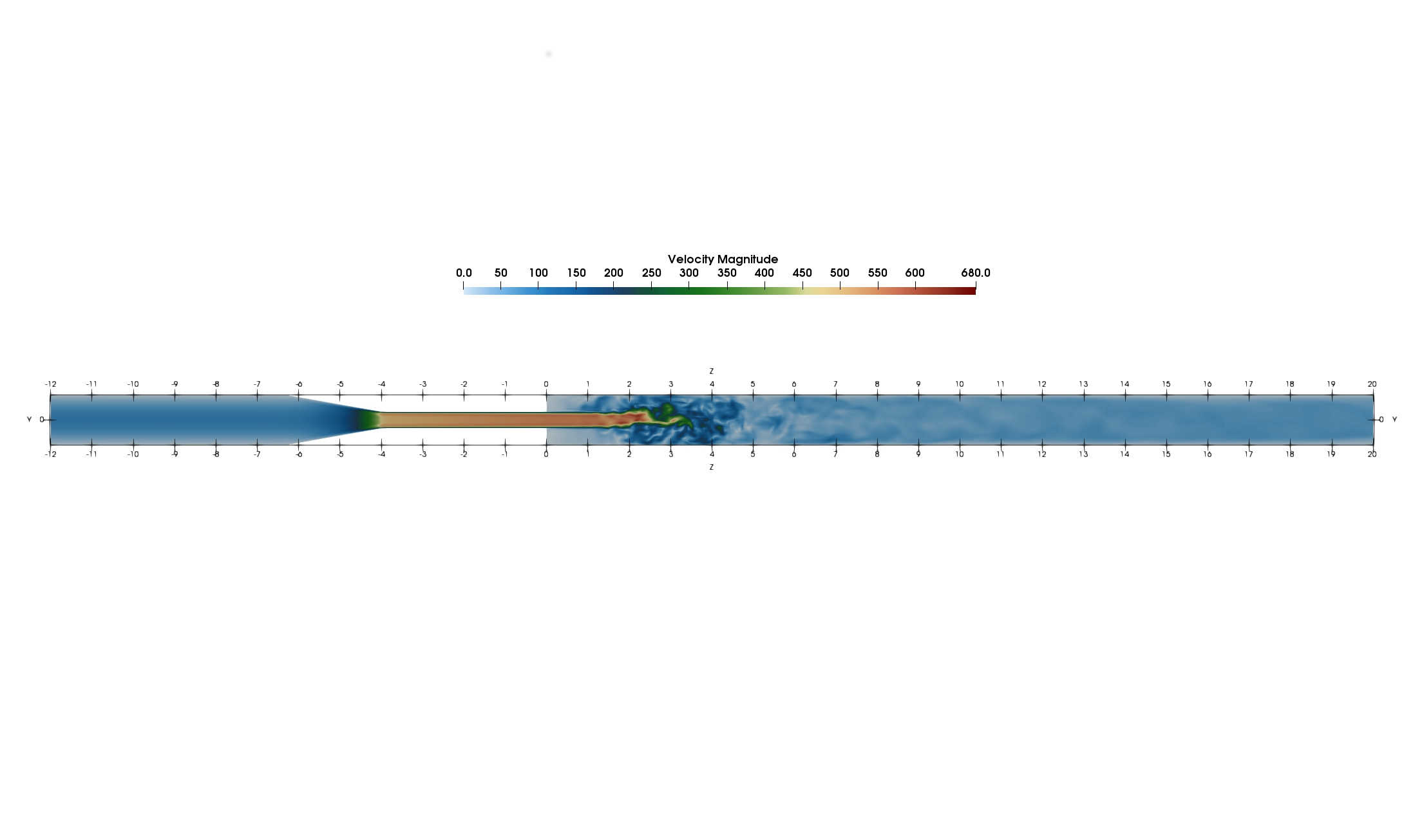}
\end{tabular}
\end{center}
\caption{Instantaneous velocity magnitude in cm/s for $\textup{Re}^{\textup{t}}=500$ (top), $\textup{Re}^{\textup{t}}=3500$ (middle), and $\textup{Re}^{\textup{t}}=5000$ (bottom).} 
\label{fig:fda_benchmark_sanpshots}
\end{figure}

\begin{figure}
	\begin{center}
	\begin{tabular}{cc}
\includegraphics[angle=0, trim=90 50 200 100, clip=true, scale = 0.08]{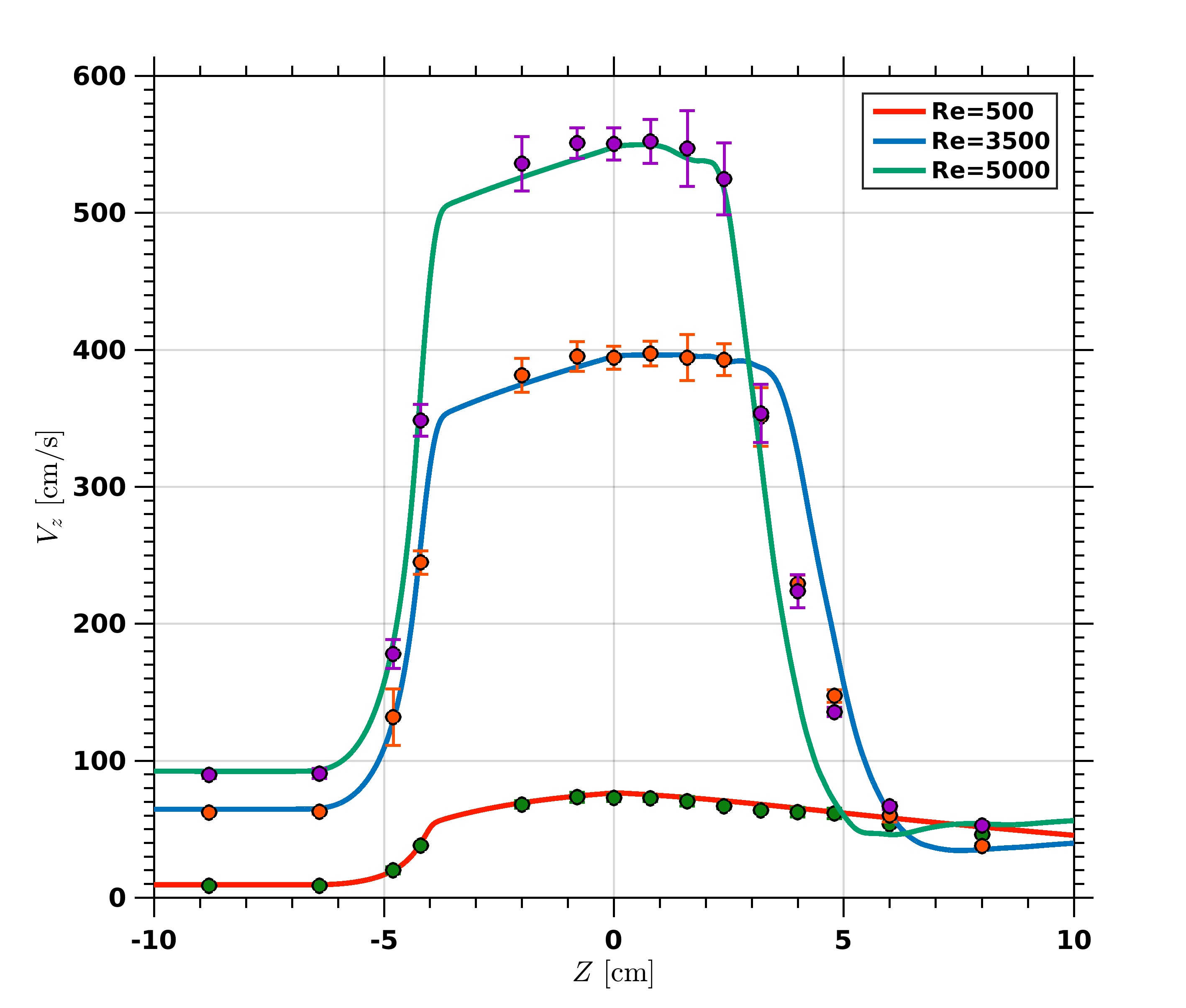} &
\includegraphics[angle=0, trim=90 50 200 100, clip=true, scale = 0.08]{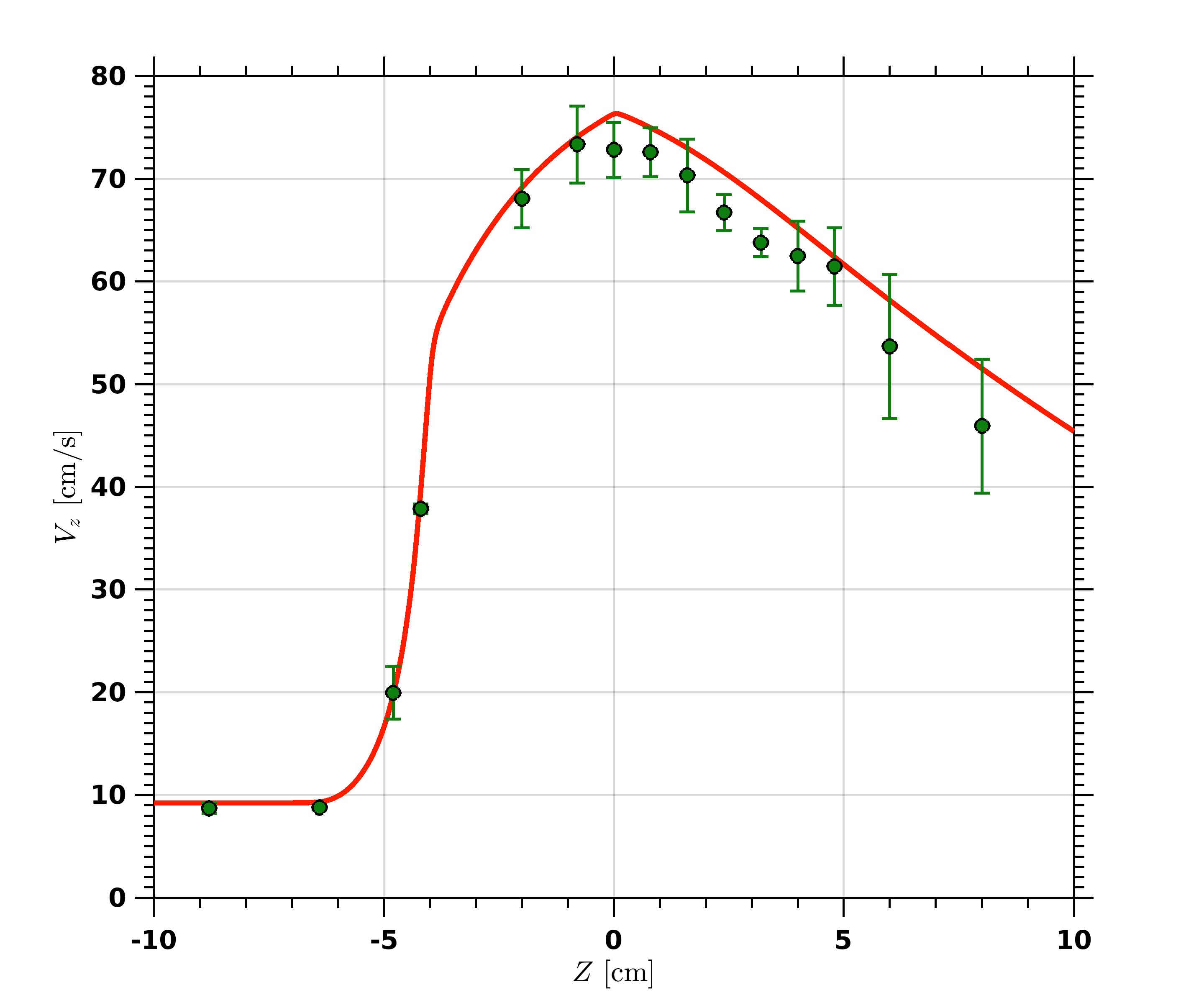} \\
(a) & (b)\\
\includegraphics[angle=0, trim=90 50 200 100, clip=true, scale = 0.08]{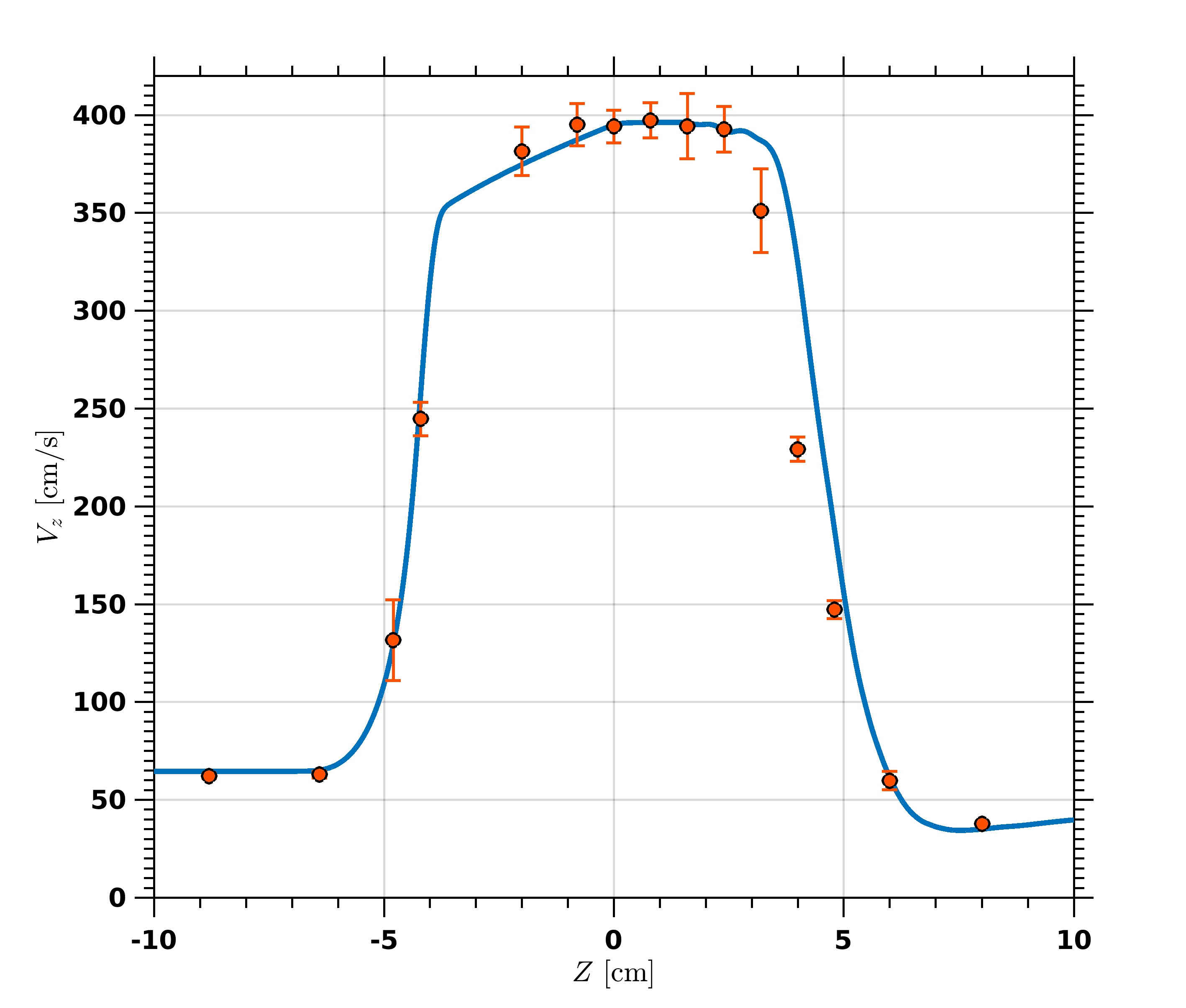} &
\includegraphics[angle=0, trim=90 50 200 100, clip=true, scale = 0.08]{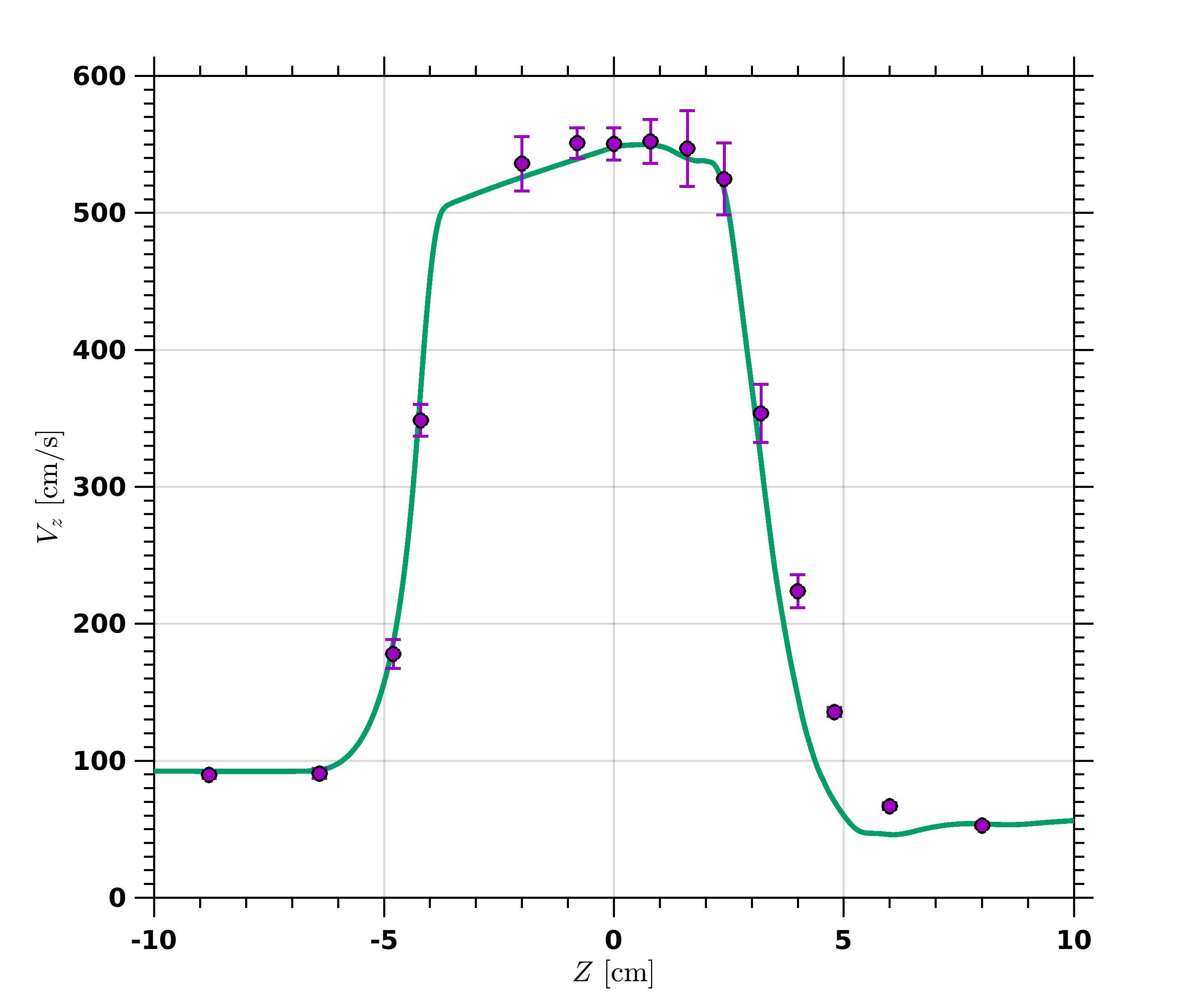} \\
(c) & (d)
\end{tabular}
\end{center}
\caption{The axial velocity along the nozzle centerline for all three Reynolds numbers (a). Detailed views of the velocity are given for $\textup{Re}^{\textup{t}}=500$ (b), $3500$ (c), and $5000$ (d). Solid lines are time-averaged velocities from the CFD simulations, and the symbols represent experimental means with 95\% confidence interval \cite{FDA_Critical_Path_Project_webpage}.} 
\label{fig:fda_benchmark_axial_velo}
\end{figure}

Before investigating this problem, the CFD solver has been verified using the manufactured solution method, and optimal convergence rates for the velocity and pressure have been observed. Here, we perform a suite of simulations to further validate it. An isotropic unstructured mesh is generated by Gmsh \cite{Geuzaine2009} using the frontal algorithm. The mesh consists of $28.8$ million linear tetrahedral elements and $4.94$ million vertices. In particular, we do not perform local mesh refinement near the sudden expansion region to ease the reproducibility of the results and to give a fair evaluation of the predictive capability of our CFD code. The minimum and maximum element sizes of the mesh are $\Delta x_{\textup{min}} = 0.014$ cm and $\Delta x_{\textup{max}} = 0.069$ cm, respectively. In Table \ref{table:fda_setup}, the volumetric flow rates of the inflow, the flow speed, the corresponding Reynolds numbers, and the time step sizes are summarized. With the varying values of $\mathcal Q$, the Reynolds numbers at the throat take the values of $500$, $3500$, and $5000$. These values correspond to the laminar, transitional, and fully turbulent regimes, thereby offering a wide range of flow characteristics for code validation. The inflow velocity profile is perturbed by a random velocity field with zero mean and a standard deviation 1\% of the mean axial velocity \cite{Bergersen2019,Zmijanovic2017}. We start the simulation with zero velocity and the inflow rate is gradually increased to reach the target values. In Figure \ref{fig:fda_benchmark_sanpshots}, the instantaneous velocity magnitudes for the three cases considered are depicted over a cut in the middle of the geometry. For the case of $\textup Re^{\textup{t}}=500$, the flow reaches a steady state and shows an axially symmetric pattern. For the cases of $\textup{Re}^{\textup{t}}=3500$ and $5000$, the flows are highly unsteady, with an anticipated jet break point in the sudden expansion region. In Figure \ref{fig:fda_benchmark_axial_velo}, the time-averaged axial velocities along the nozzle centerline are compared with experimental data obtained by different laboratories \cite{FDA_Critical_Path_Project_webpage}. As can be observed from the figures, the simulation results show a very good agreement with the experimental data.

\begin{table}[htbp]
\begin{center}
\tabcolsep=0.3cm
\renewcommand{\arraystretch}{1.5}
\begin{tabular}{ c c c c c c }
\hline
$\mathcal Q$ [cm$^3$/s]& $\bm v^{\textup{i}}_{\textup{max}}$ [cm/s] & $\textup{Re}^{\textup{i}}$ & $\bm v^{\textup{t}}_{\textup{max}}$ [cm/s]  & $\textup{Re}^{\textup{t}}$ & $\Delta t$ [s] \\
\hline
$5.2062 \times 10^0$ & $9.21 \times 10^0$ & $167$ & $8.29 \times 10^1$  & $500$ & $1.0 \times 10^{-4}$ \\
$3.6444 \times 10^1$ & $6.45 \times 10^1$ & $1167$ & $5.80 \times 10^2$  & $3500$ & $5.0 \times 10^{-5}$ \\
$5.2062 \times 10^1$ & $9.21 \times 10^1$ & $1667$ & $8.29 \times 10^2$  & $5000$ & $2.5 \times 10^{-5}$ \\
\hline 
\end{tabular}
\end{center}
\caption{The flow rates $\mathcal Q$, maximum velocities in the inlet and the throat, the corresponding Reynolds numbers, and the time step sizes used in the simulations. The maximum velocities are calculated based on a parabolic profile.}
\label{table:fda_setup}
\end{table}

\begin{figure}
	\begin{center}
	\begin{tabular}{cc}
\includegraphics[angle=0, trim=30 30 200 120, clip=true, scale = 0.07]{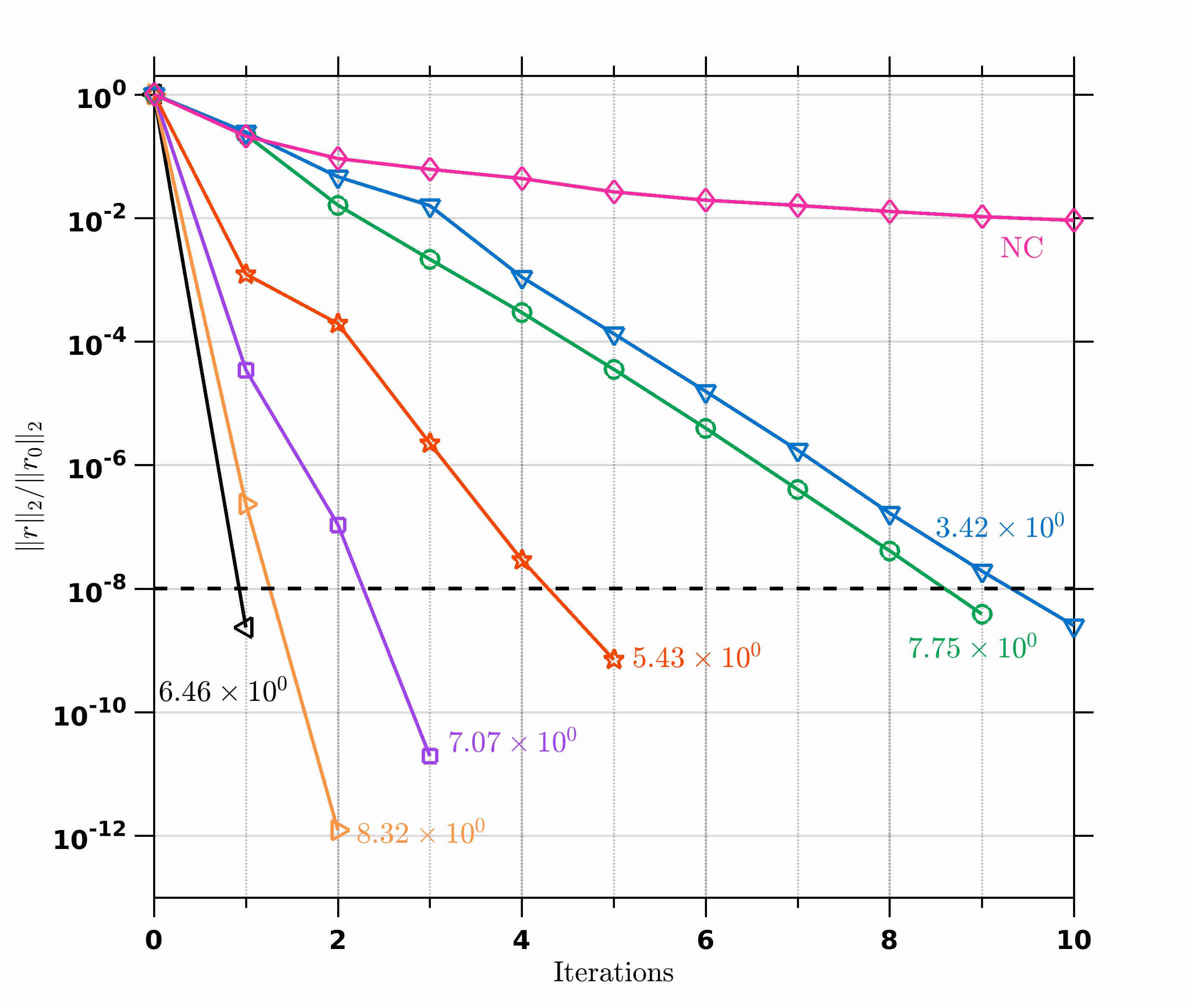} &
\includegraphics[angle=0, trim=30 30 200 120, clip=true, scale = 0.07]{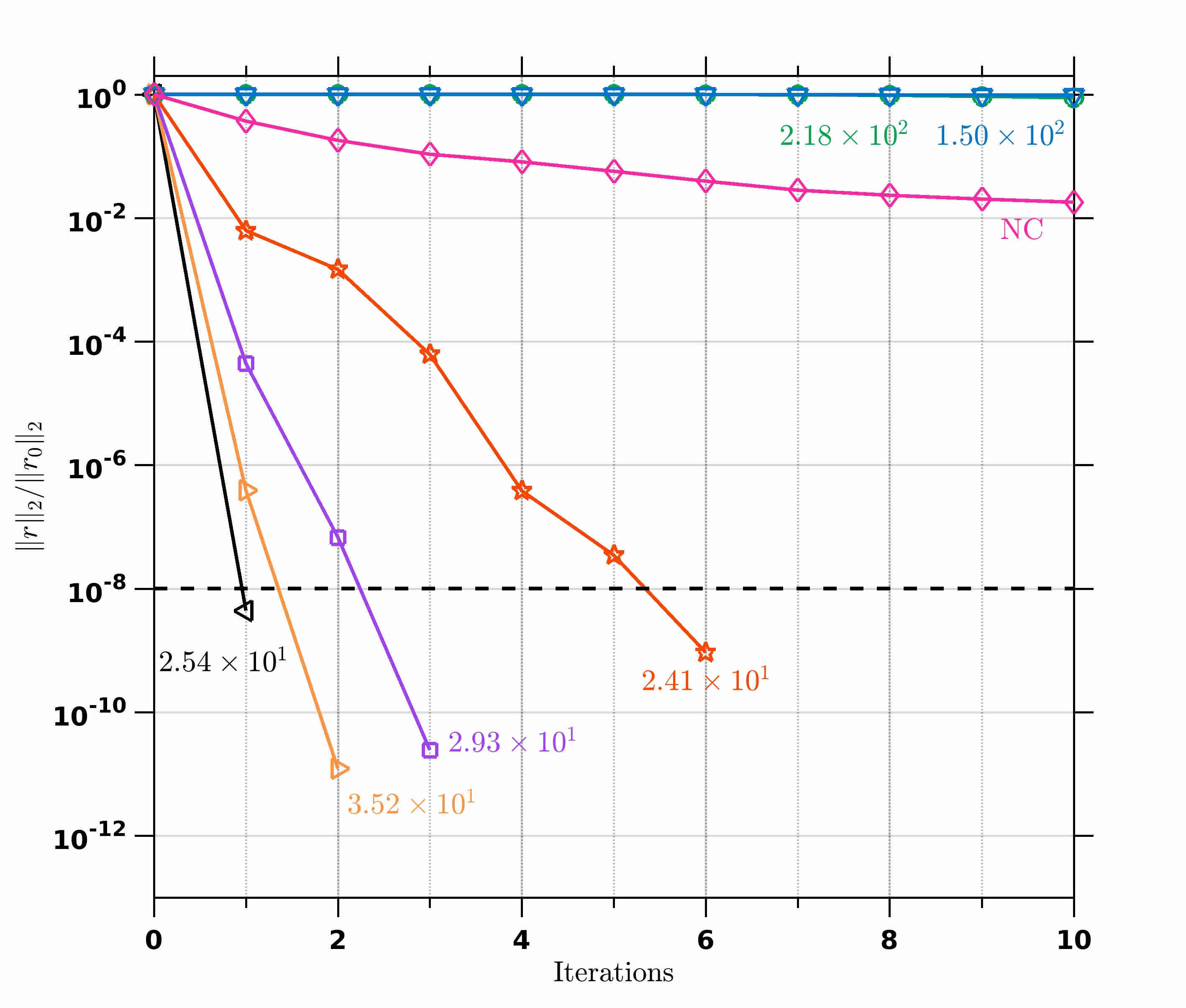} \\
(a) $\textup{Re}^{\textup{t}}=500$, $\textup{Cr}=0.5$ & (b) $\textup{Re}^{\textup{t}}=500$, $\textup{Cr}=2.0$ \\[1mm]
\includegraphics[angle=0, trim=30 30 200 120, clip=true, scale = 0.07]{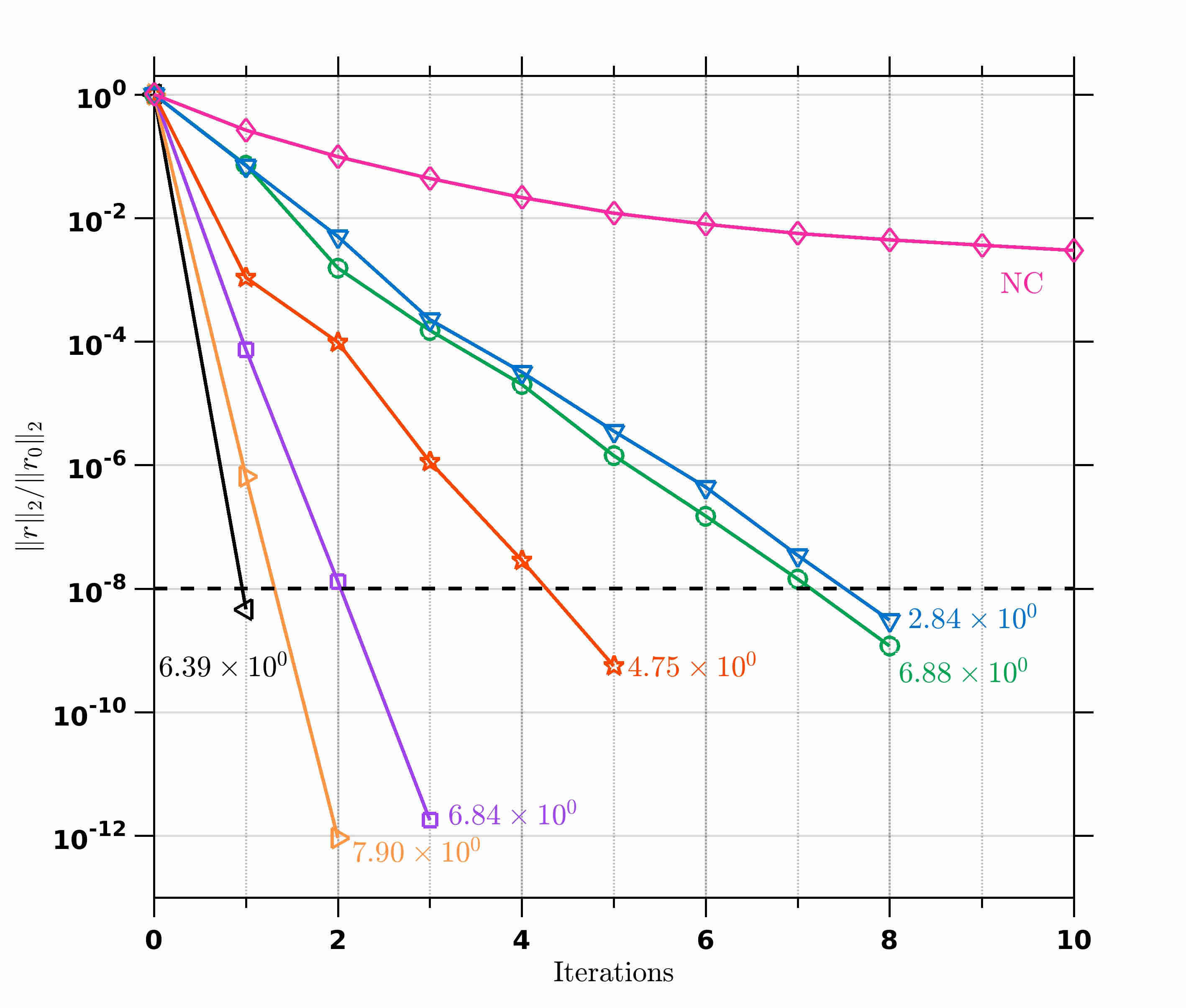} &
\includegraphics[angle=0, trim=30 30 200 120, clip=true, scale = 0.07]{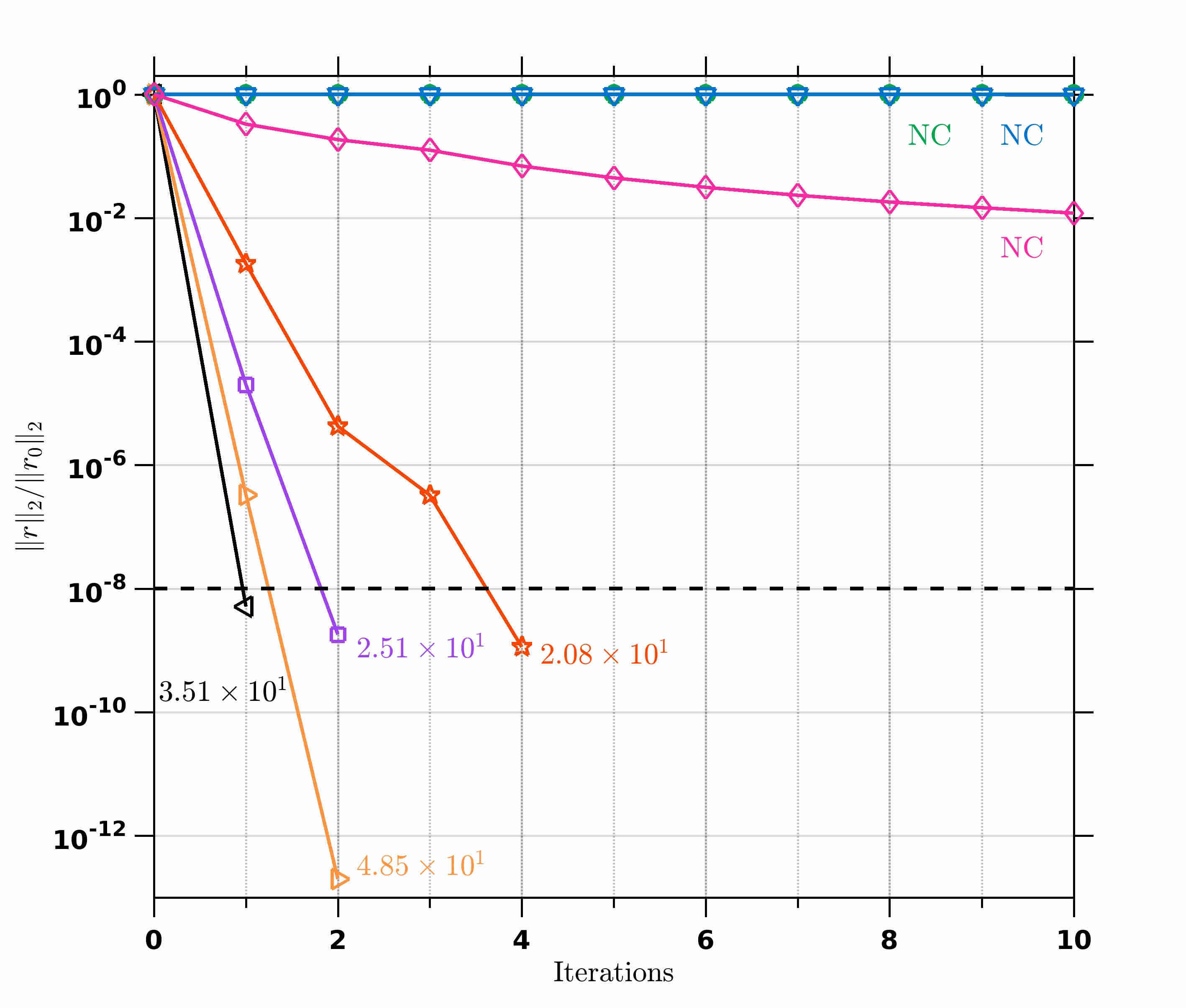} \\
(c) $\textup{Re}^{\textup{t}}=3500$, $\textup{Cr}=0.5$ & (d) $\textup{Re}^{\textup{t}}=3500$, $\textup{Cr}=2.0$ \\[1mm]
\includegraphics[angle=0, trim=30 30 200 120, clip=true, scale = 0.07]{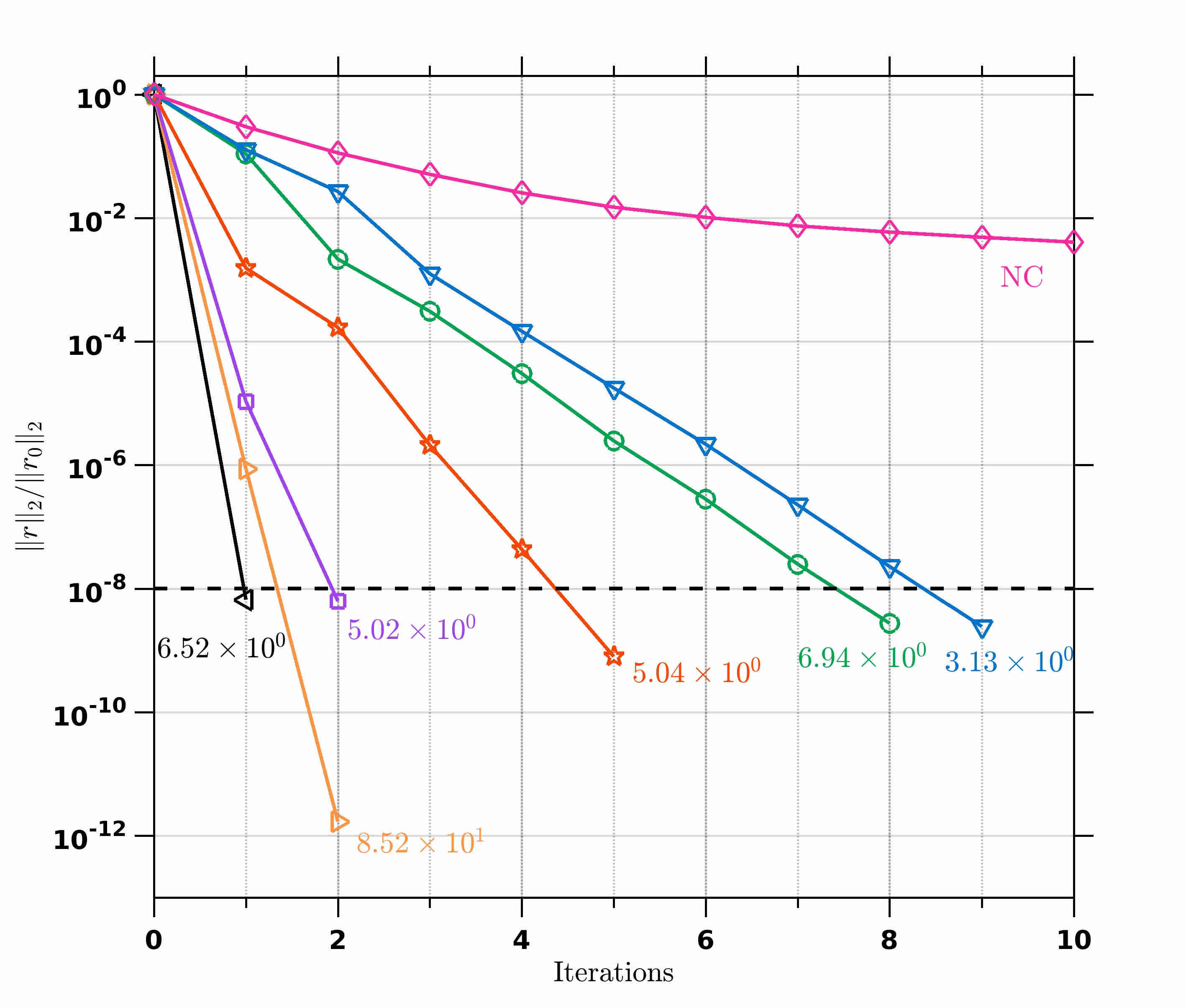} &
\includegraphics[angle=0, trim=30 30 200 120, clip=true, scale = 0.07]{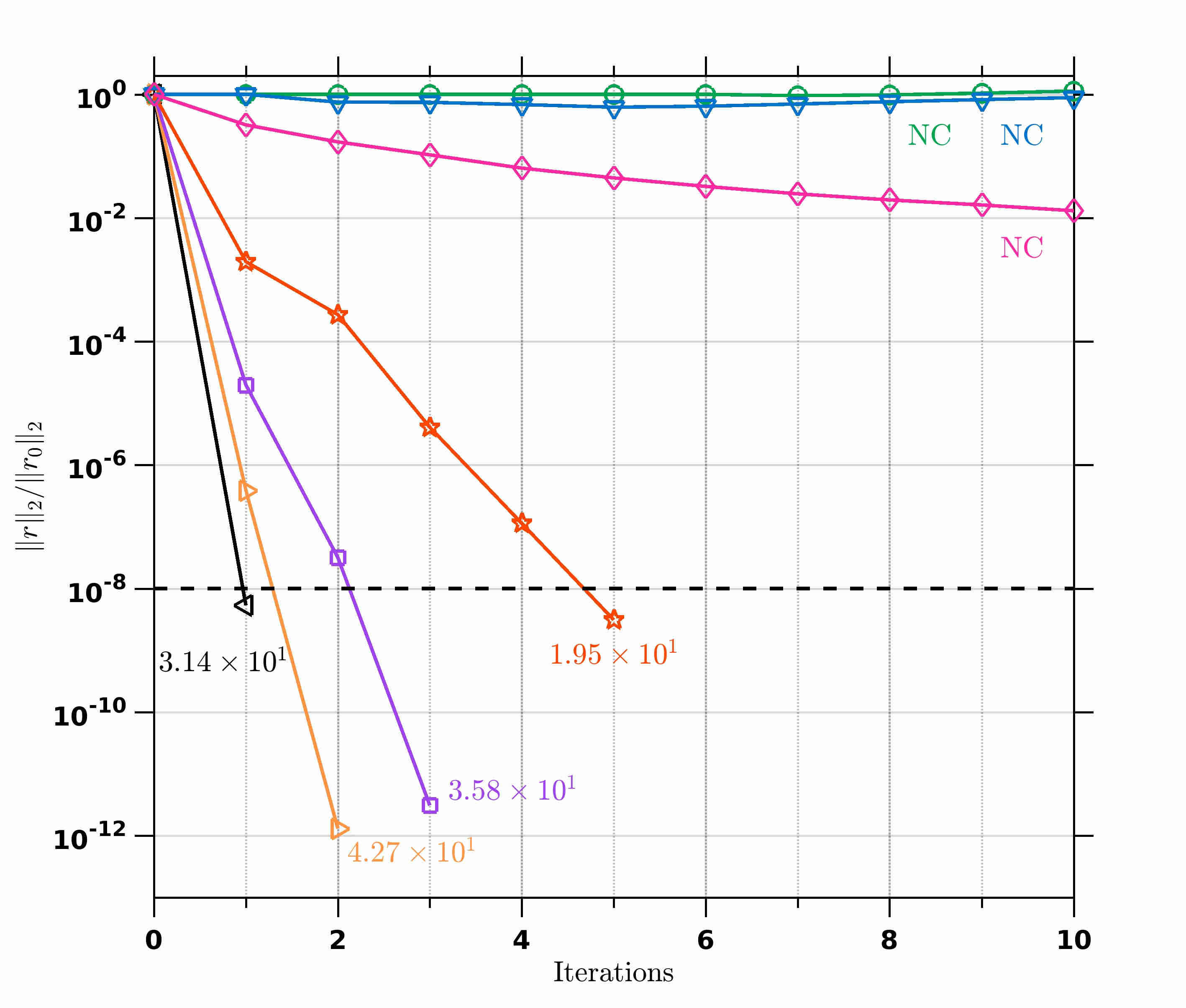} \\
(e) $\textup{Re}^{\textup{t}}=5000$, $\textup{Cr}=0.5$ & (f) $\textup{Re}^{\textup{t}}=5000$, $\textup{Cr}=2.0$ \\[3mm]
\multicolumn{2}{c}{ \includegraphics[angle=0, trim=1020 90 840 310, clip=true, scale = 0.12]{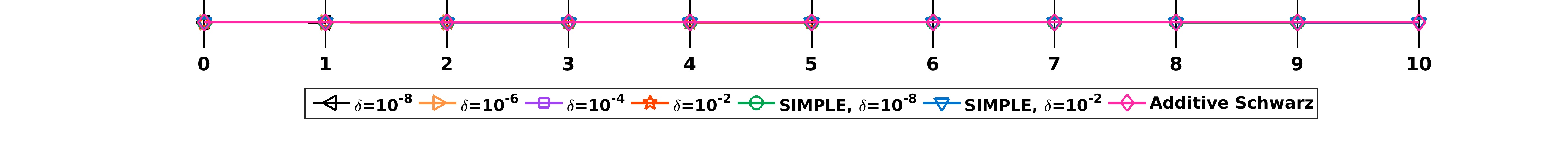}
}
\end{tabular}
\caption{Convergence history for $\textup{Cr}=0.5$ (left column) and $2.0$ (right column) for $\textup{Re}^{\textup{t}}=500$ (top row), $\textup{Re}^{\textup{t}}=3500$ (middle row), and $\textup{Re}^{\textup{t}}=5000$ (bottom row). The horizontal dashed black line indicates the prescribed stopping criterion for the relative error, which is $10^{-8}$ here. The numbers indicate the CPU time for the linear solver in the unit of seconds. NC stands for no convergence within the prescribed number of iterations. In the legend, $\delta$ represents the relative tolerance used in the intermediate and inner solvers.} 
\label{fig:conv_history_fda}
\end{center}
\end{figure}

\subsubsection{Algorithm robustness}
\label{sec:performance_with_varying_reynolds_numbers}
In this test, we examine the solver performance with different Reynolds numbers. The spatial mesh contains $7.82\times 10^5$ tetrahedral elements and $1.43\times 10^5$ vertices. Correspondingly, the system of equations involves $5.71\times 10^5$ unknowns. The minimum and maximum element sizes of this mesh are $0.06$ cm and $0.29$ cm, respectively. We choose to simulate with two different Courant numbers, $\textup{Cr}=0.5$ and $2.0$. The time step sizes based on the Courant number and the Reynolds number are summarized in Table \ref{table:fda_cfl}. The simulations are performed on a single ``Skylake" node with $6$ CPUs assigned. We choose $\delta^r = 10^{-8}$ for the FGMRES iterations. For the nested block preconditioner, we set $\delta^r_A = \delta^r_S = \delta^r_I$ and vary their values from $10^{-8}$ to $10^{-2}$. For the SIMPLE preconditioner, we investigated two choices: $\delta^r_A = \delta^r_S = 10^{-8}$ and $10^{-2}$. For comparison purposes, we also applied the restricted additive Schwarz method implemented in PETSc \cite{petsc-user-ref} with the incomplete LU factorization as the subdomain solver. The convergence histories of this set of simulations are depicted in Figure \ref{fig:conv_history_fda}. 

\begin{table}[htbp]
\begin{center}
\tabcolsep=0.3cm
\renewcommand{\arraystretch}{1.5}
\begin{tabular}{ c c c c }
\hline
$\textup{Cr}$ & $\textup{Re}^{\textup{t}}=500$  & $\textup{Re}^{\textup{t}}=3500$  & $\textup{Re}^{\textup{t}}=5000$ \\
\hline
$0.5$ & $3.69\times 10^{-4}$ & $5.28\times 10^{-5}$ & $3.69\times 10^{-5}$ \\
$2.0$ & $1.48\times 10^{-3}$ & $2.11\times 10^{-4}$ & $1.48\times 10^{-4}$ \\
\hline 
\end{tabular}
\end{center}
\caption{The time step size $\Delta t$ (in the unit of seconds) for different Courant number $\textup{Cr}$ and Reynolds number $\textup{Re}^{\textup{t}}$ in the simulations of Sec. \ref{sec:performance_with_varying_reynolds_numbers}.}
\label{table:fda_cfl}
\end{table}

There are several salient features of the solver performance that can be summarized from Figure \ref{fig:conv_history_fda}. First, the nested block preconditioner always shows a faster rate in residual reduction. For all cases considered, the proposed method requires at most $6$ iterations to achieve convergence. In terms of the CPU time, we found that a looser tolerance for the intermediate and inner solvers always gives more efficient performance. Indeed, for the nested block preconditioner, a significant amount of computing time is spent on the matrix-free approximation of the Schur complement. A looser tolerance can surely reduce the computational burden from that part. This suggests that, in practice, one is advised to flexibly tune the tolerances for $\delta^r_A$, $\delta^r_S$, and $\delta^r_I$ to balance the solver robustness and the cost per iteration (see, e.g., the discussion in Sec. 5.1.1 of \cite{Liu2019}). Also, in our experience, AMG is not the most cost-effective approach for solving with $\boldsymbol{\mathrm A}$, especially in the inner solver. Using a `lightweight' preconditioner may further reduce the overhead.

Second, the SIMPLE preconditioner achieves convergence within ten iterations for all three cases when $\textup{Cr} = 0.5$. The choice of the intermediate solver accuracy does not have a significant impact on the overall convergence rate. The SIMPLE preconditioner with the intermediate solver accuracy $10^{-2}$ gives the fastest performance for all three cases of $\textup{Cr} = 0.5$. However, the drawback of the SIMPLE preconditioner is demonstrated when the Courant number is greater than $1$. When $\textup{Cr} = 2.0$, the SIMPLE preconditioner only achieves convergence for the case of $\textup{Re}^{\textup{t}} = 500$, with around 130 iterations. For larger Courant and Reynolds numbers, the matrix $\boldsymbol{\mathrm A}$ is no more diagonally dominated, and thus the SIMPLE preconditioner shows rather poor performance. The efficacy of the inner solver can be observed from Figure \ref{fig:conv_history_fda} (b), (d), and (e). Even a loose tolerance, such as $\delta^r_I = 10^{-2}$, can lead to a significant improvement of the robustness of the overall method.

Third, the additive Schwarz domain decomposition preconditioner cannot drive the relative error to $10^{-8}$ within the prescribed number of iterations ($10000$ here) for any of the six cases. It is also known that this preconditioner does not scale well. An increase in the problem size will further degrade its performance in terms of the iteration number and CPU time cost.

\subsubsection{Parallel performance}
In this section, we examine the proposed preconditioning technique by investigating its performance under the parallel setting. The stopping criterion for the FGMRES iteration is $\delta^r = 10^{-6}$, and the tolerances for the intermediate and inner solvers are $\delta^r_A = \delta^r_S = \delta^r_I = 10^{-3}$.

\paragraph{Fixed-size scalability}
In this test, we investigate the fixed-size scalability of the algorithm. Rather than starting from a zero solution, we prepared a developed flow profile at $\textup{Re}^{\textup{t}}=3500$. The spatial mesh contains  $5.71\times 10^6$ tetrahedral elements and $1.00\times 10^6$ vertices, and the overall system of equations involves $4\times 10^6$ unknowns. Then the same mesh is re-partitioned based on the given number of processors. The prepared flow profile is mapped to the new partitioned mesh to serve as the initial condition. The time step size is fixed to be $5\times 10^{-5}$ s and the problem is integrated for $40$ time steps. Due to the hierarchical architecture of Stampede2, the simulations are performed with one CPU assigned in each node, which means the number of processors reported in Table \ref{table:strong_scaling} equals the number of nodes. In doing so, the MPI messages are communicated purely through the Intel Omni-Path network. In Table \ref{table:strong_scaling}, we report the timings as well as the speed-up efficiency based on the total time as the number of CPUs increases. It can be observed that the scaling shows nearly optimal speed-up efficiency for a wide range of processor counts for the proposed nested block preconditioner.

\begin{table}[htbp]
\begin{center}
\tabcolsep=0.19cm
\renewcommand{\arraystretch}{1.2}
\begin{tabular}{P{2.0cm} P{2.0cm} P{2.0cm} P{3.0cm} P{2.5cm} P{2.0cm} }
\hline
Proc. & $T_A$ (s) & $T_L$ (s) & Total (s) & Efficiency   \\
\hline
2 &  $2.24\times 10^3$ & $3.33 \times 10^3$ & $5.64\times 10^3$ & $100\%$ \\
4 &  $1.25\times 10^3$ & $1.74 \times 10^3$ & $3.03\times 10^3$ & $93\%$ \\
8 &  $6.15\times 10^2$ & $9.05 \times 10^2$ & $1.54\times 10^3$ & $92\%$ \\
16 &  $2.91\times 10^2$ & $4.78 \times 10^2$ & $7.80\times 10^2$ & $91\%$ \\
32 &  $1.25\times 10^2$ & $2.39 \times 10^2$ & $3.72\times 10^2$ & $95\%$ \\
64 &  $6.38\times 10^1$ & $1.42 \times 10^2$ & $2.11\times 10^2$ & $84\%$ \\
128 &  $3.39\times 10^1$ & $7.80 \times 10^1$ & $1.17\times 10^2$ & $75\%$ \\
256 &  $1.84\times 10^1$ & $5.53 \times 10^1$ & $8.04\times 10^1$ & $55\%$ \\
\hline
\end{tabular}
\end{center}
\caption{The strong scaling performance for the FDA idealized medical device benchmark. In the table, $T_A$ and $T_L$ represent the timings for the matrix assembly and the linear solver, respectively; the efficiency is computed based on the total runtime.}
\label{table:strong_scaling}
\end{table}

\paragraph{Isogranular scalability}
In the next test, we study the isogranular or weak scalability of the proposed preconditioning technique. We generated three sets of isotropic unstructured meshes and partitioned the mesh with the goal of maintaining a constant number of unknowns assigned to each CPU. The initial conditions are prepared by gradually increasing the flow rate to reach the targeted values in one second and maintaining the targeted flow rates for another second. Then the obtained solutions are used as the initial condition for the scaling test. For each different run, the time step size is chosen based on the Courant number, and we select to simulate with $\textup{Cr} = 0.5$, $1.0$, and $2.0$. The statistics of the solver performance are collected for $20$ time steps (Table \ref{table:weak_scaling}). It can be observed that the averaged number of iterations remains between $2$ and $3$ for different problem sizes and the Courant numbers. For a given Courant number, the averaged CPU time grows mildly (three to four times longer with a sixty-four-fold increase of the problem size). 


\begin{table}
\begin{center}
\tabcolsep=0.19cm
\renewcommand{\arraystretch}{1.2}
\begin{tabular}{@{\extracolsep{3pt}}P{1.5cm} P{1.5cm} P{0.9cm} P{1.0cm} P{1.0cm} P{1.0cm} P{1.0cm} P{1.5cm} P{1.0cm}@{}}
\hline
\multirow{2}{*}{$n_{en}$} & \multirow{2}{*}{$n_{eq}$} & \multirow{2}{*}{Proc.} & \multicolumn{2}{c}{$\textup{Cr} = 0.5$} & \multicolumn{2}{c}{$\textup{Cr} = 1.0$} & \multicolumn{2}{c}{$\textup{Cr} = 2.0$}\\
\cline{4-5} \cline{6-7}\cline{8-9}
& & & $\bar{n}$ & $\bar{T}_L$ (s) & $\bar{n}$ & $\bar{T}_L$ (s) & $\bar{n}$ & $\bar{T}_L$ (s)
\\
\hline
\multicolumn{9}{l}{$\textup{Re}^{\textup{t}} = 500$} \\
$4.20\times 10^5$ & $3.10\times 10^5$  &  6 & 2.45 & 1.82 & 2.41 & 2.12 & 2.47 & 4.20 \\
$3.65\times 10^6$ & $2.48\times 10^6$  &  48 & 2.50 & 2.83 & 2.60 & 3.34 & 2.89 & 9.34 \\
$2.88\times 10^7$ & $1.98\times 10^7$  &  384 & 2.90 & 6.97 & 2.95 & 7.78 & 3.02 & 15.05 \\
\hline
\multicolumn{9}{l}{$\textup{Re}^{\textup{t}} = 3500$} \\
$4.20\times 10^5$ & $3.10\times 10^5$  &  6 & 2.38 & 1.99 & 2.36 & 2.23 & 2.13 & 4.23 \\
$3.65\times 10^6$ & $2.48\times 10^6$  &  48 & 2.62 & 3.20 & 2.31 & 3.11 & 2.27 & 6.51 \\
$2.88\times 10^7$ & $1.98\times 10^7$  &  384 & 2.78 & 6.95 & 2.58 & 6.96 & 2.45 & 11.58 \\
\hline
\multicolumn{9}{l}{$\textup{Re}^{\textup{t}} = 5000$} \\
$4.20\times 10^5$ & $3.10\times 10^5$  &  6 & 2.35 & 1.83 & 2.22 & 2.12 & 2.21 & 4.40 \\
$3.65\times 10^6$ & $2.48\times 10^6$  &  48 & 2.76 & 3.24 & 2.52 & 3.31 & 2.33 & 6.64 \\
$2.88\times 10^7$ & $1.98\times 10^7$  &  384 & 2.83 & 6.96 & 2.52 & 7.17 & 2.41 & 13.41 \\
\hline
\end{tabular}
\end{center}
\caption{Comparison of the averaged iteration counts $\bar{n}$ and the CPU time for solving the linear system $\bar{T}_L$ in seconds for the nested block preconditioner. The values of $\Delta x_{\textup{min}}$ for the three meshes are $6.8\times 10^{-2}$, $3.1\times 10^{-2}$, $1.4\times 10^{-2}$ cm, respectively. }
\label{table:weak_scaling}
\end{table}

\subsection{A cylindrical model with resistance boundary condition}
In this example, we investigate the preconditioner performance with a resistance boundary condition \eqref{eq:resistance_bc}. The geometry of this problem can be found as the first example in the simulation guide of SimVascular \cite{simvascular-simulation-guide}, which is a straight cylinder with a radius $2$ cm and length $30$ cm. On the inlet surface, we apply a parabolic velocity profile with the flow rate $\mathcal Q = 100$ cm$^3$/s. On the outlet boundary surface, we apply the resistance boundary condition \eqref{eq:resistance_bc} with varying values of $\mathrm R$. \footnote{Notice that there is only one outlet surface $\Gamma_h = \Gamma_{\mathrm{out}}^1$ in this case, and we neglect the superscript of $\mathrm R^1$ for notational simplicity.} The mesh is generated using TetGen \cite{Si2015} through the SimVascular package \cite{Updegrove2017}, and it consists of $1.74 \times 10^6$ isotropic tetrahedral elements and $2.88 \times 10^5$ vertices. The minimum element size is $\Delta x_{\textup{min}} = 0.1$ cm. For comparison, we also simulate the same problem using \texttt{svSolver} \cite{Updegrove2017,simvascular-simulation-guide}. The solution method of \texttt{svSolver} is based on the BIPN iterative algorithm introduced in Sec. \ref{subsec:bipn}. \textit{We caution the readers that the comparisons made between the proposed method and the \texttt{svSolver} are not apples-to-apples because the numerical formulation and the resulting matrix problem are different.} Comparisons are made merely to give overall evaluations of the two existing solvers under realistic settings.

\begin{figure}
	\begin{center}
	\begin{tabular}{cc}
\includegraphics[angle=0, trim=30 30 200 120, clip=true, scale = 0.072]{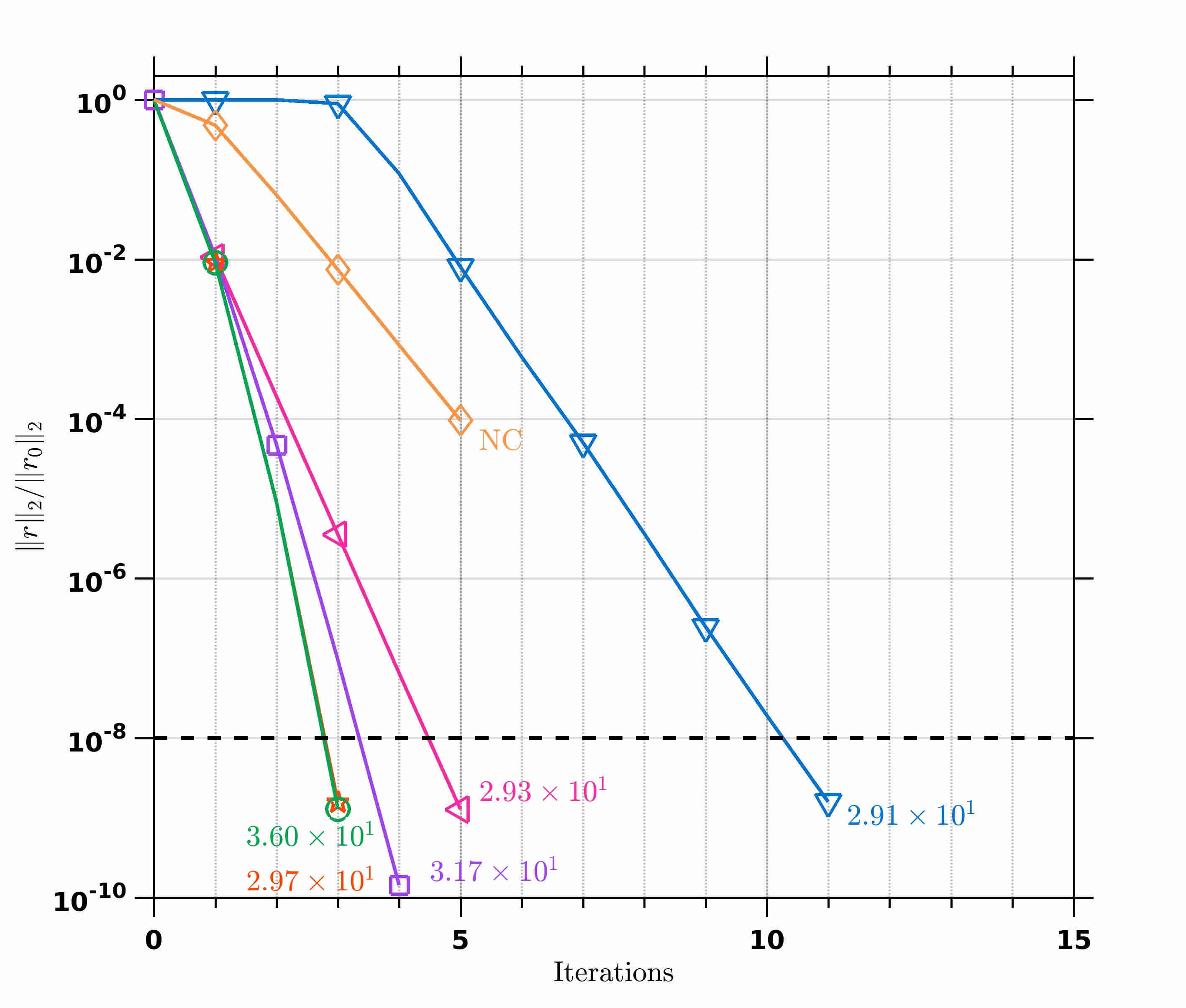} &
\includegraphics[angle=0, trim=30 30 200 120, clip=true, scale = 0.072]{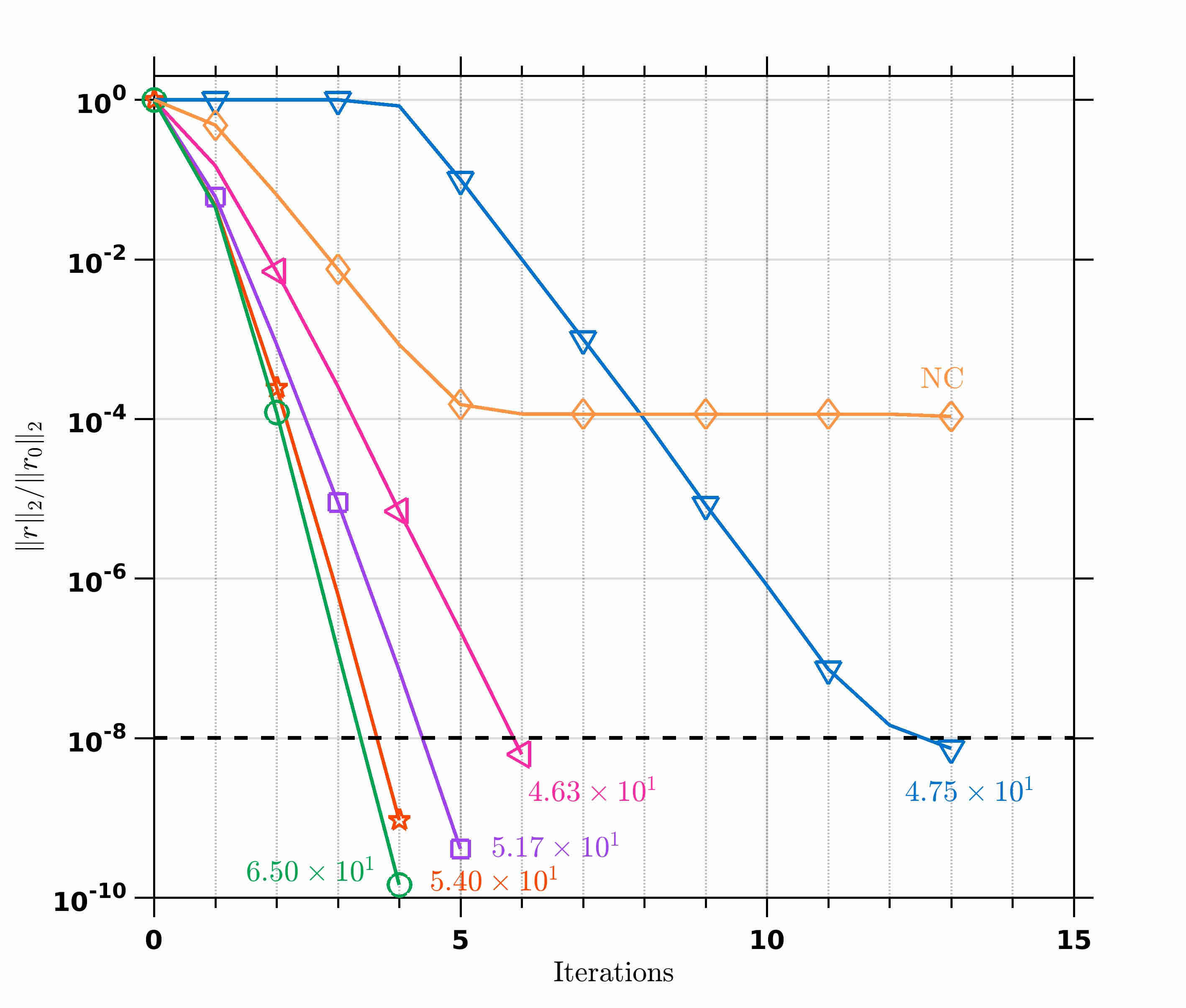} \\
(a) $\mathrm R=10^2$ g/(s cm$^4$) & (b) $\mathrm R=10^3$ g/(s cm$^4$) \\
\includegraphics[angle=0, trim=30 30 200 120, clip=true, scale = 0.072]{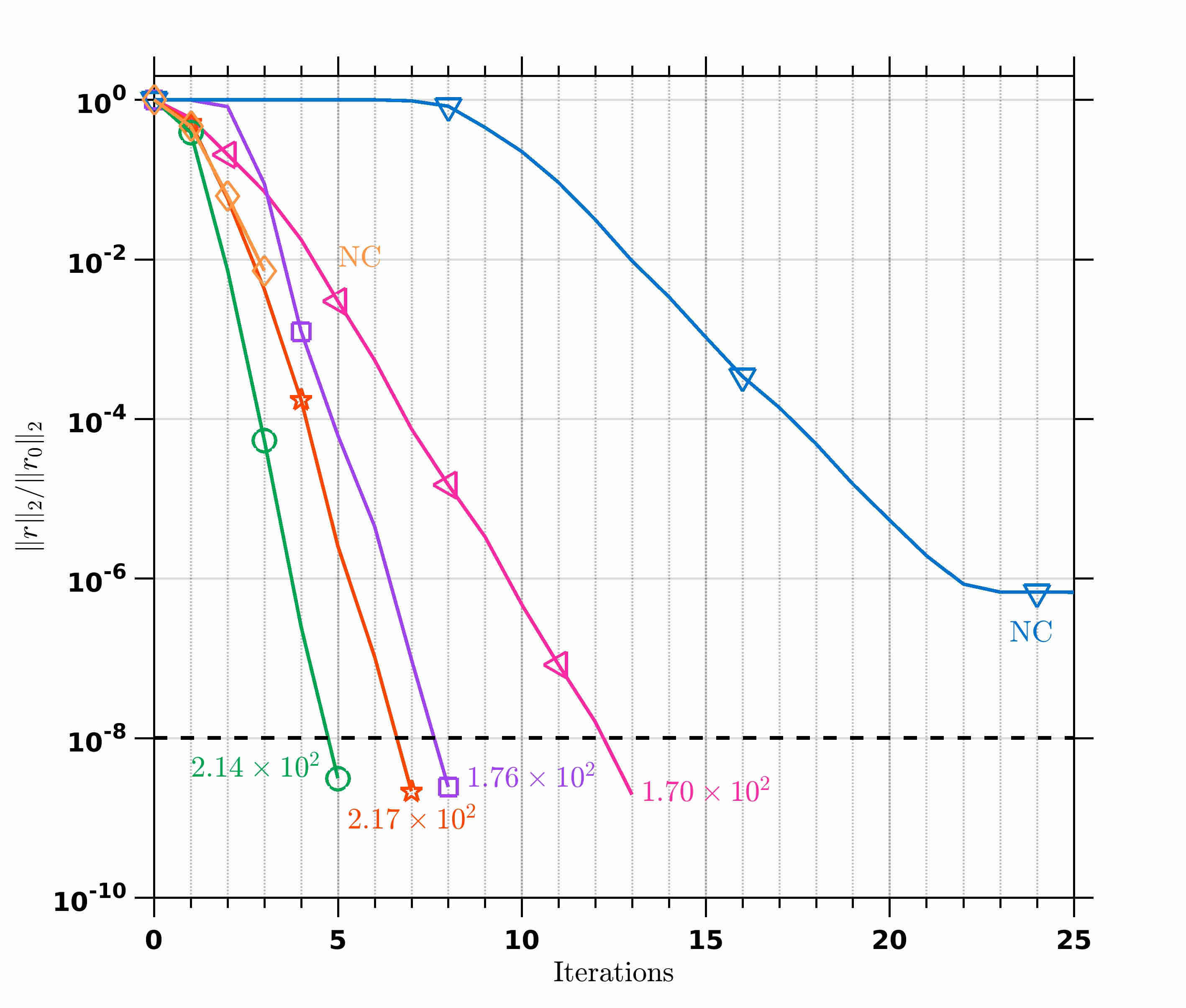} &
\includegraphics[angle=0, trim=30 30 200 120, clip=true, scale = 0.072]{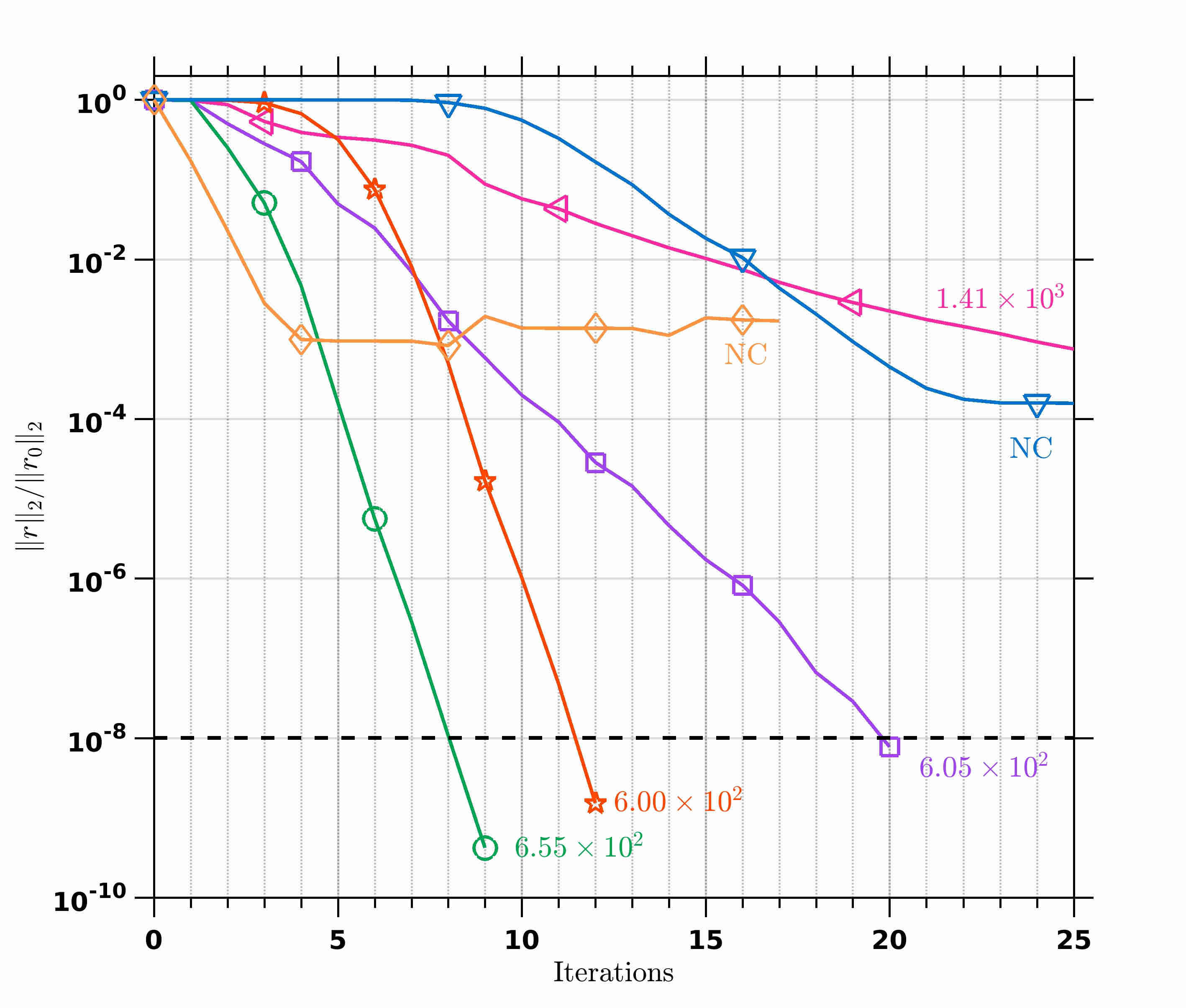} \\
(a) $\mathrm R=10^4$ g/(s cm$^4$) & (b) $\mathrm R=10^5$ g/(s cm$^4$) \\[3mm]
\multicolumn{2}{c}{ \includegraphics[angle=0, trim=1600 80 1420 295, clip=true, scale = 0.12]{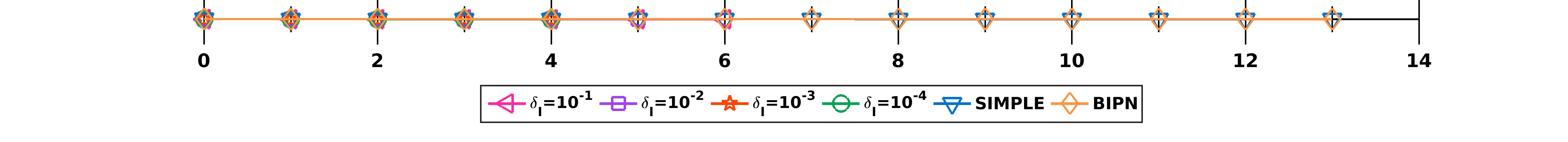} }
\end{tabular}
\caption{Convergence history for different resistance values. The horizontal dashed black line indicates the prescribed stopping criterion for the relative error, which is $10^{-8}$ here. The numbers indicate the CPU time for the linear solver in the unit of seconds. NC stands for no convergence within the prescribed number of iterations.} 
\label{fig:conv_history_cyl}
\end{center}
\end{figure}

In the first study, we choose $\Delta t = 6.3\times 10^{-3}$ s for the simulation so that the Courant number $\textup{Cr} = 1.0$ based on the maximum flow speed and $\Delta x_{\textup{min}}$. The simulations are performed with $24$ CPUs. The resistance value $\mathrm R$ varies from $10^2$ g/(s cm$^4$) to $10^5$ g/(s cm$^4$). Again, we choose $\delta^r = 10^{-8}$ as the stopping criterion. In the SIMPLE and nested block preconditioners, we choose $\delta^r_A = \delta^r_S = 10^{-4}$, and the value of $\delta^r_I$ varies from $10^{-4}$ to $10^{-1}$. We also simulated the same problem using \texttt{svSolver} with tolerances on the momentum and continuity equations fixed to $10^{-4}$. The convergence histories and the CPU time of the linear system solution procedure are depicted in Figure \ref{fig:conv_history_cyl}. In all cases considered, \texttt{svSolver} shows a rapid reduction of the residual in the first a few iterations and stagnates when the relative error is driven to approximately $10^{-4}$. As one proceeds in the linear iteration, the normal matrix (see Eqn. (45) in \cite{Esmaily-Moghadam2015}) in the BIPN algorithm becomes nearly singular and causes divergence eventually. As we increase the resistance value, the SIMPLE preconditioner is also observed to stagnate. With an inner solver invoked, the proposed algorithm robustly drive the residual to the prescribed tolerance. For example, when $\mathrm R=10^5$, the algorithm converges in $64$ iterations with $\delta^r_I = 10^{-1}$. This fact again indicates that a slight improvement of the numerical representation of the Schur complement may significantly improve the convergence rate of the overall iterative algorithm.

\begin{figure}
	\begin{center}
	\begin{tabular}{cc}
\includegraphics[angle=0, trim=30 20 200 120, clip=true, scale = 0.083]{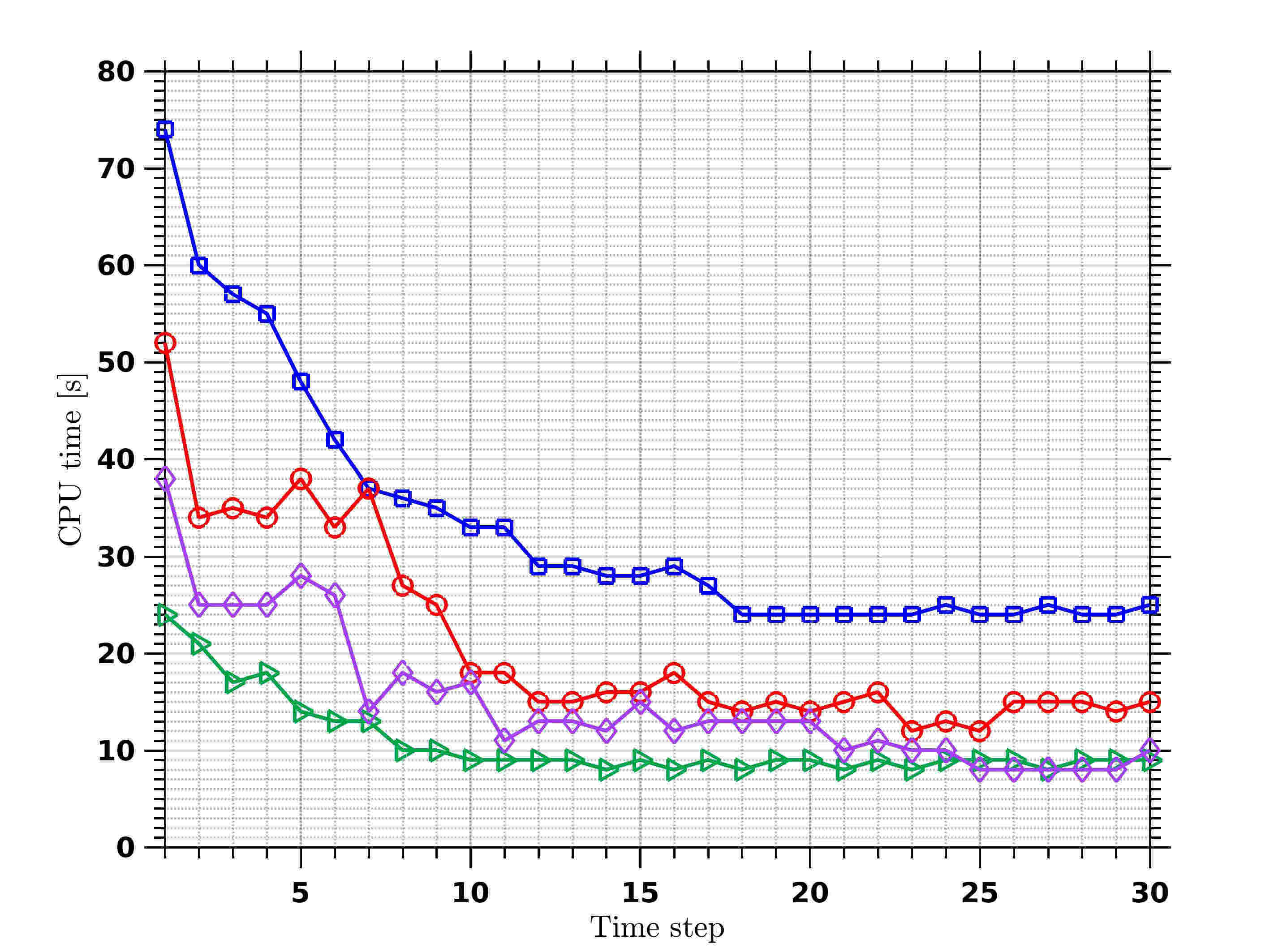} &
\includegraphics[angle=0, trim=30 20 200 120, clip=true, scale = 0.083]{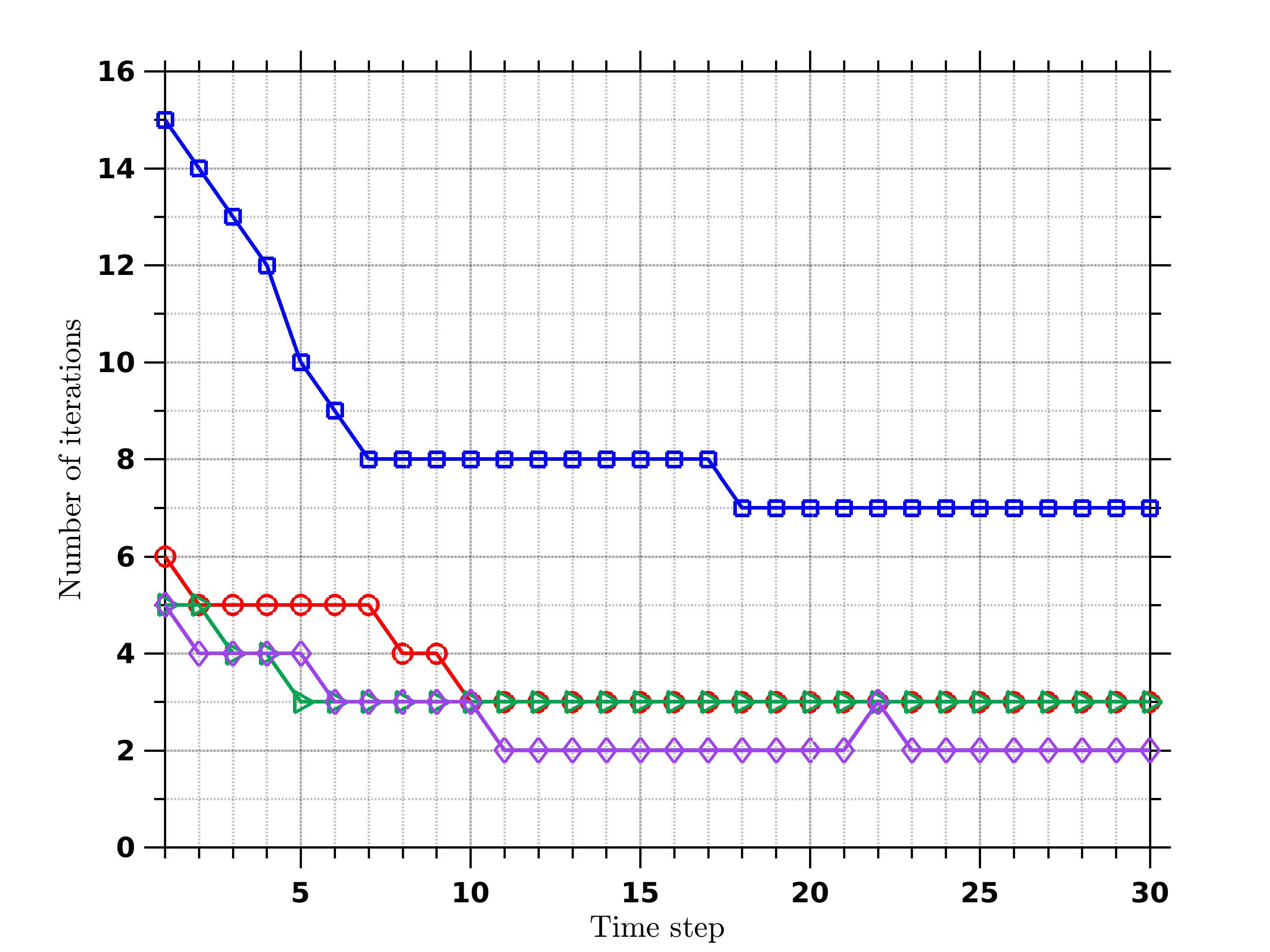} \\
\multicolumn{2}{c}{ \includegraphics[angle=0, trim=1460 80 1280 295, clip=true, scale = 0.12]{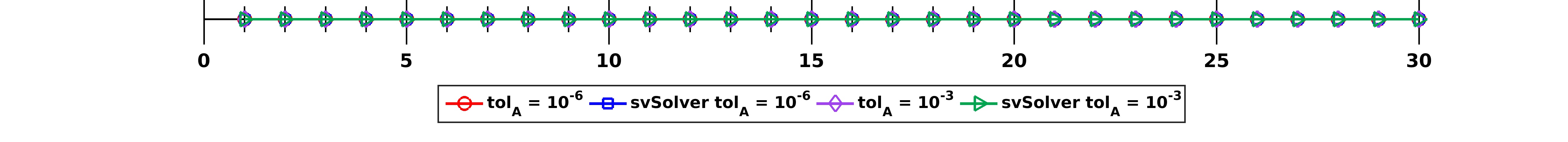} }
\end{tabular}
\caption{Comparison of the CPU time (left) and the number of the multi-corrector iterations (right) in the first $30$ time steps between the proposed algorithm and the \texttt{svSolver}.} 
\label{fig:nest_sv_compare}
\end{center}
\end{figure}

In the second study, we choose to compare the impact of the two solvers on the overall time-stepping algorithms for $30$ time steps. The resistance value is fixed to be $1333$ g/(s cm$^4$), leading to a pressure value of 100 mmHg at the outlet. The time step size is chosen to be $\Delta t = 3.15\times 10^{-3}$ s. The tolerance $\delta^r$ is fixed to be $10^{-3}$, with $\delta^r_A = 10^{-3}$ and $\delta^r_S = \delta^r_I = 10^{-2}$. We also use the Jacobi preconditioner for the inner and intermediate solvers associated with $\boldsymbol{\mathrm A}$, because this option is often more efficient than the AMG preconditioner (see also our discussion in Sec. \ref{sec:additional_discussion}). For \texttt{svSolver}, we limit the maximum number of iterations to $5$ to avoid the numerical instability issue; the tolerances on momentum equations and the continuity equation are fixed to be $10^{-3}$ and $10^{-2}$, respectively. In Figure \ref{fig:nest_sv_compare}, the CPU time and the number of Newton-Raphson iterations in the first $30$ time steps are reported. Notice that in \texttt{svSolver}, only the absolute tolerance is used as the stopping criterion. Therefore two sets of comparisons are made using $\textup{tol}_{\textup{A}} = 10^{-3}$ and $10^{-6}$ as the stopping criterion without invoking the relative tolerance. It takes more iterations and CPU times in the first a few time steps for both solvers. After around the fifteenth time step, both solvers behave steadily over the remaining time steps. When $\textup{tol}_{\textup{A}} = 10^{-3}$, both solvers take approximately the same amount of CPU time for each time step after the twentieth time step. When $\textup{tol}_{\textup{A}} = 10^{-6}$, \texttt{svSolver} requires more iterations and CPU time. The proposed algorithm can robustly solve the consistent tangent matrix to the prescribed tolerance, which leads to a quadratic convergence rate in the Newton-Raphson iteration. For example, after the fifteenth time step, one iteration can drive the relative error below $10^{-3}$ and two iterations drive it below $10^{-6}$. For \texttt{svSolver}, it often cannot give an accurate enough solution for the linear problem within its maximum number of iterations ($5$ here). This results in a higher amount of total computing cost because more multi-corrector steps have to be taken. The efficiency of the proposed algorithm will be further corroborated in the following section.

\begin{figure}
	\begin{center}
	\begin{tabular}{cc}
\includegraphics[angle=0, trim=50 100 50 120, clip=true, scale = 0.17]{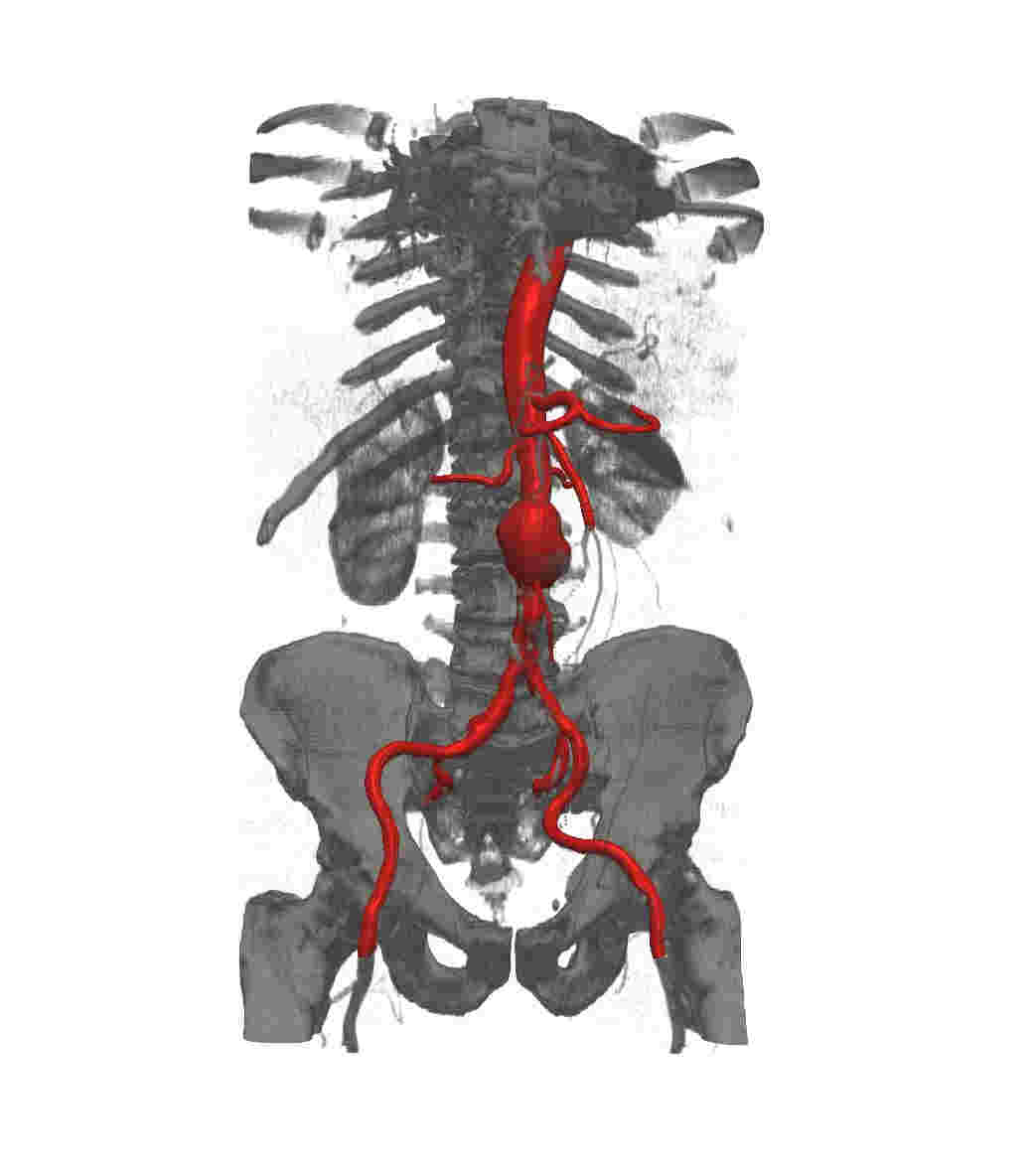} & 
\includegraphics[angle=0, trim=130 0 230 0, clip=true, scale = 0.08]{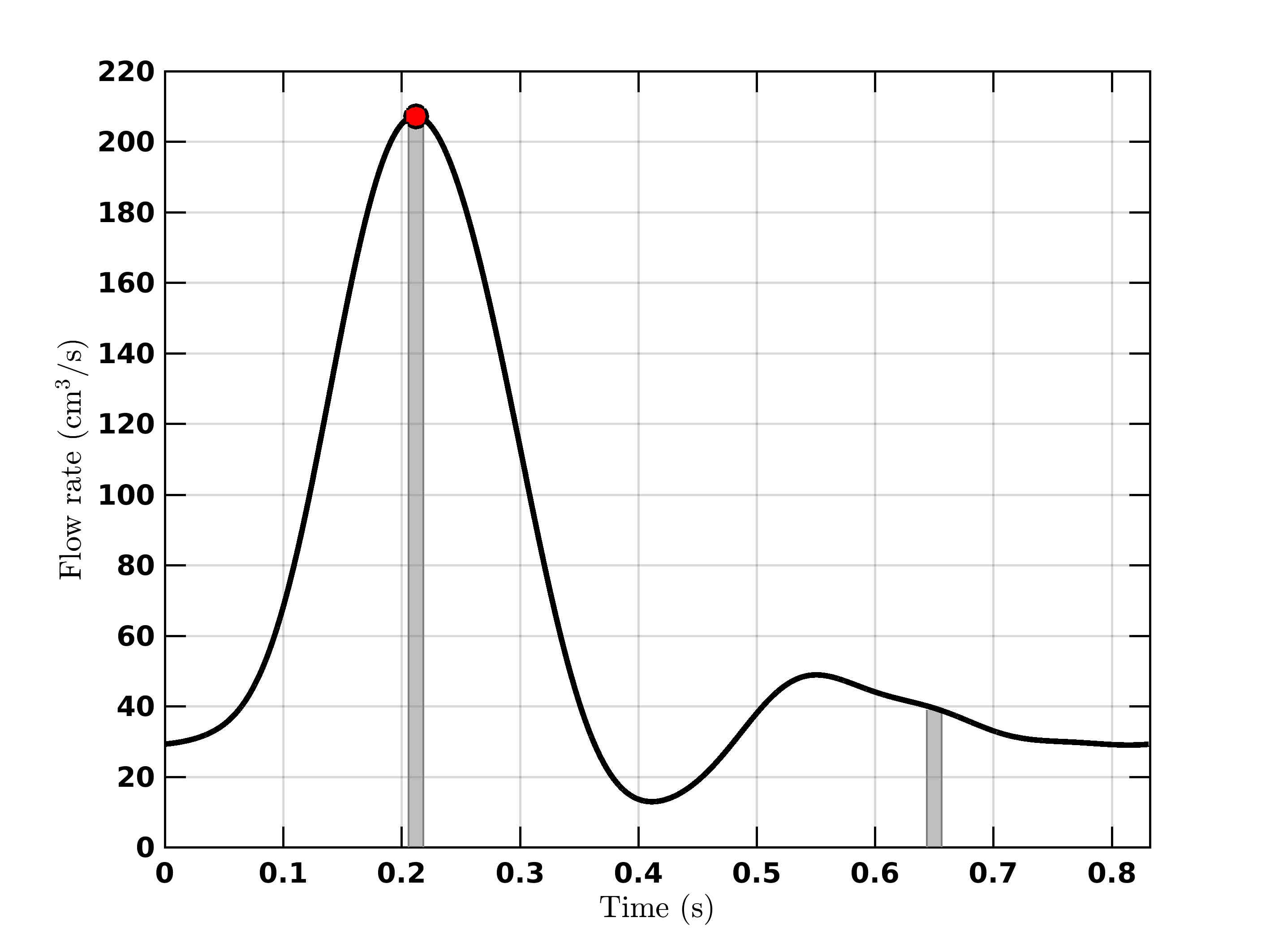}
\\
(a) & (b) \\
\multicolumn{2}{c}{ \includegraphics[angle=0, trim=100 100 100 100, clip=true, scale = 0.04]{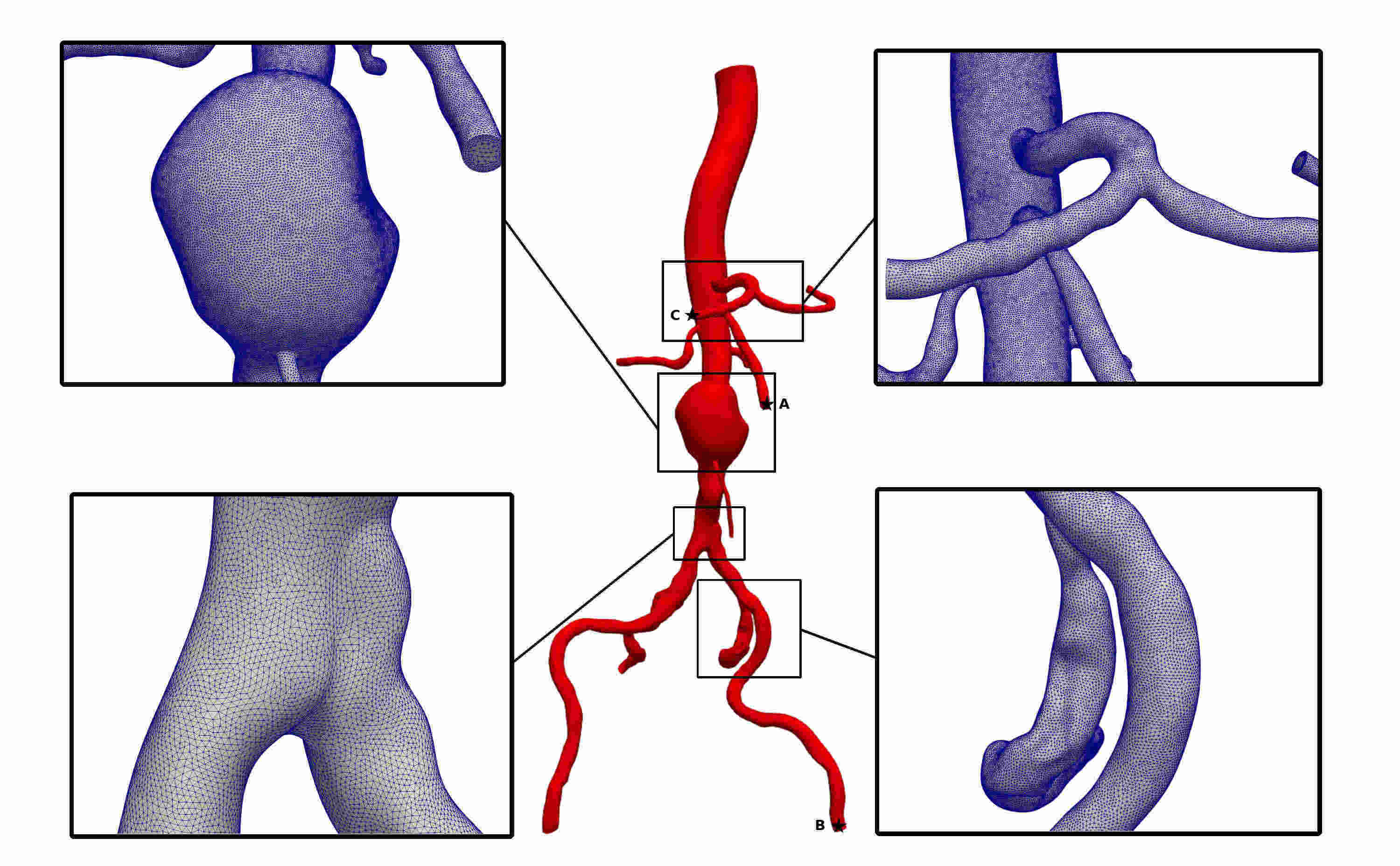} } \\
\multicolumn{2}{c}{(c)}\\
\includegraphics[angle=0, trim=120 60 20 120, clip=true, scale = 0.05]{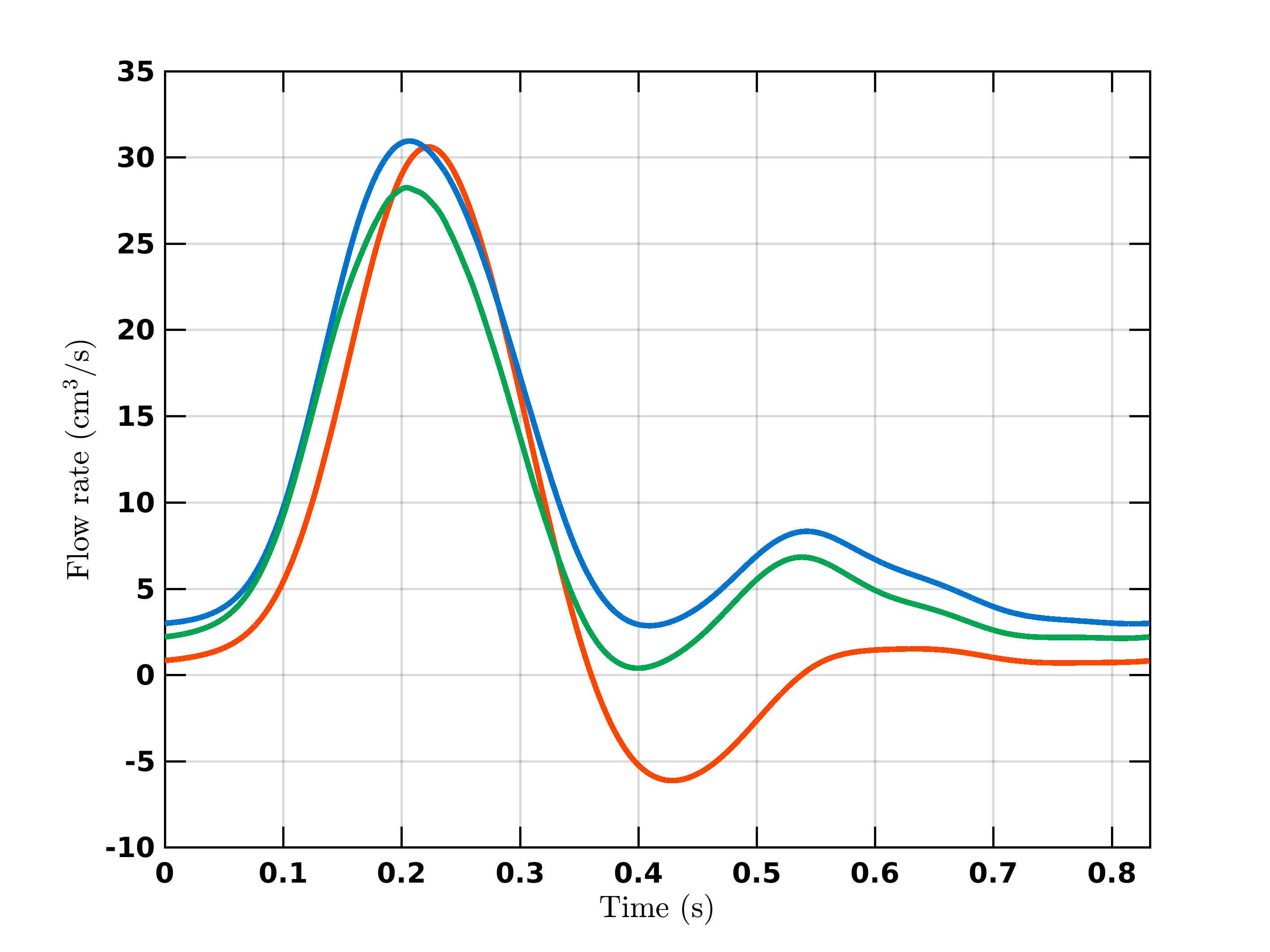} &
\includegraphics[angle=0, trim=20 60 260 120, clip=true, scale = 0.05]{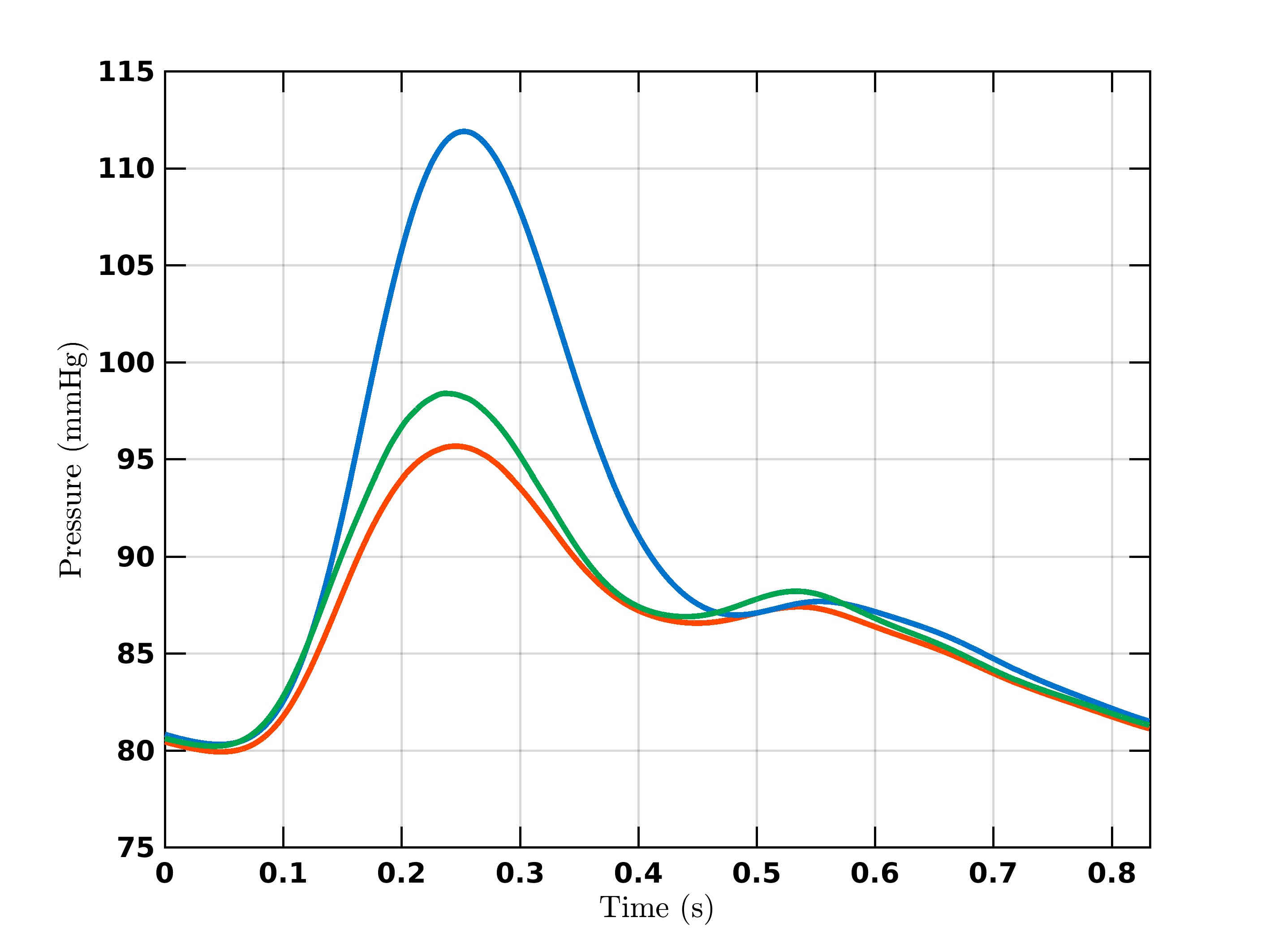} \\
(d) & (e) \\
\multicolumn{2}{c}{ \includegraphics[angle=0, trim=1270 105 1090 300, clip=true, scale = 0.12]{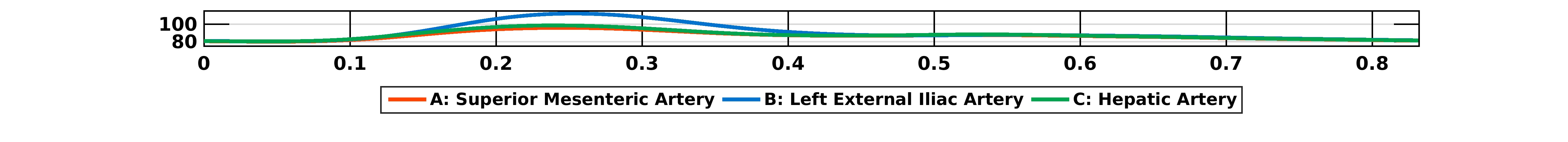}
}
\end{tabular}
\caption{(a) Volume rendering of the medical data. (b) Volumetric flow waveform used to prescribe the velocity on the inlet surface. The red dot represent the peak systole, and the grey regions represent the periods of time where we monitor the CPU time for the two different solvers. (c) The mesh of the aortofemoral model. The flow waveform (d) and pressure (e) on three outlets of the aortofemoral model.} 
\label{fig:aortofemoral_image_and_mesh}
\end{center}
\end{figure}

\begin{figure}
	\begin{center}
	\begin{tabular}{cc}
\includegraphics[angle=0, trim=500 0 600 0, clip=true, scale = 0.22]{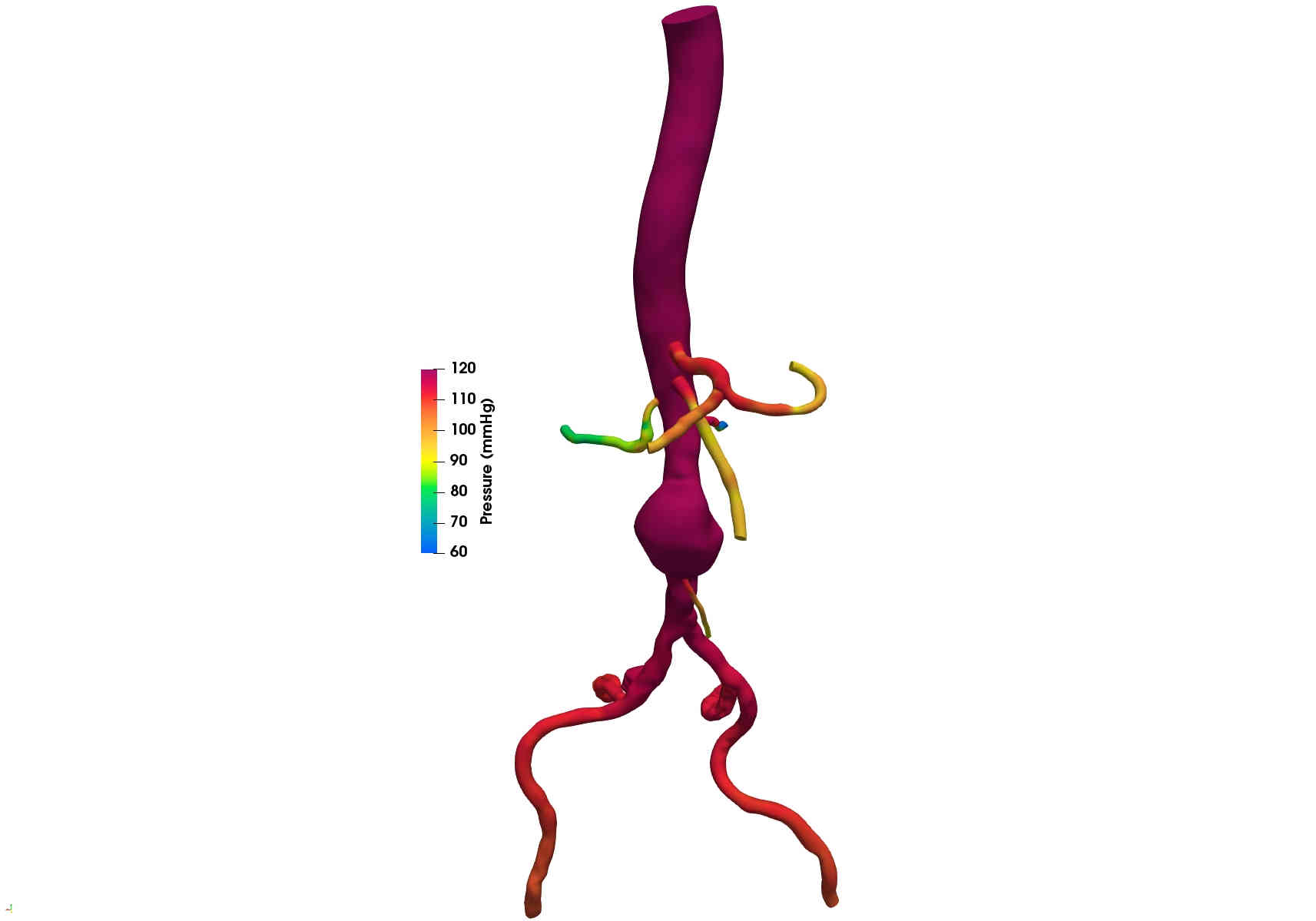} &
\includegraphics[angle=0, trim=600 0 450 0, clip=true, scale = 0.22]{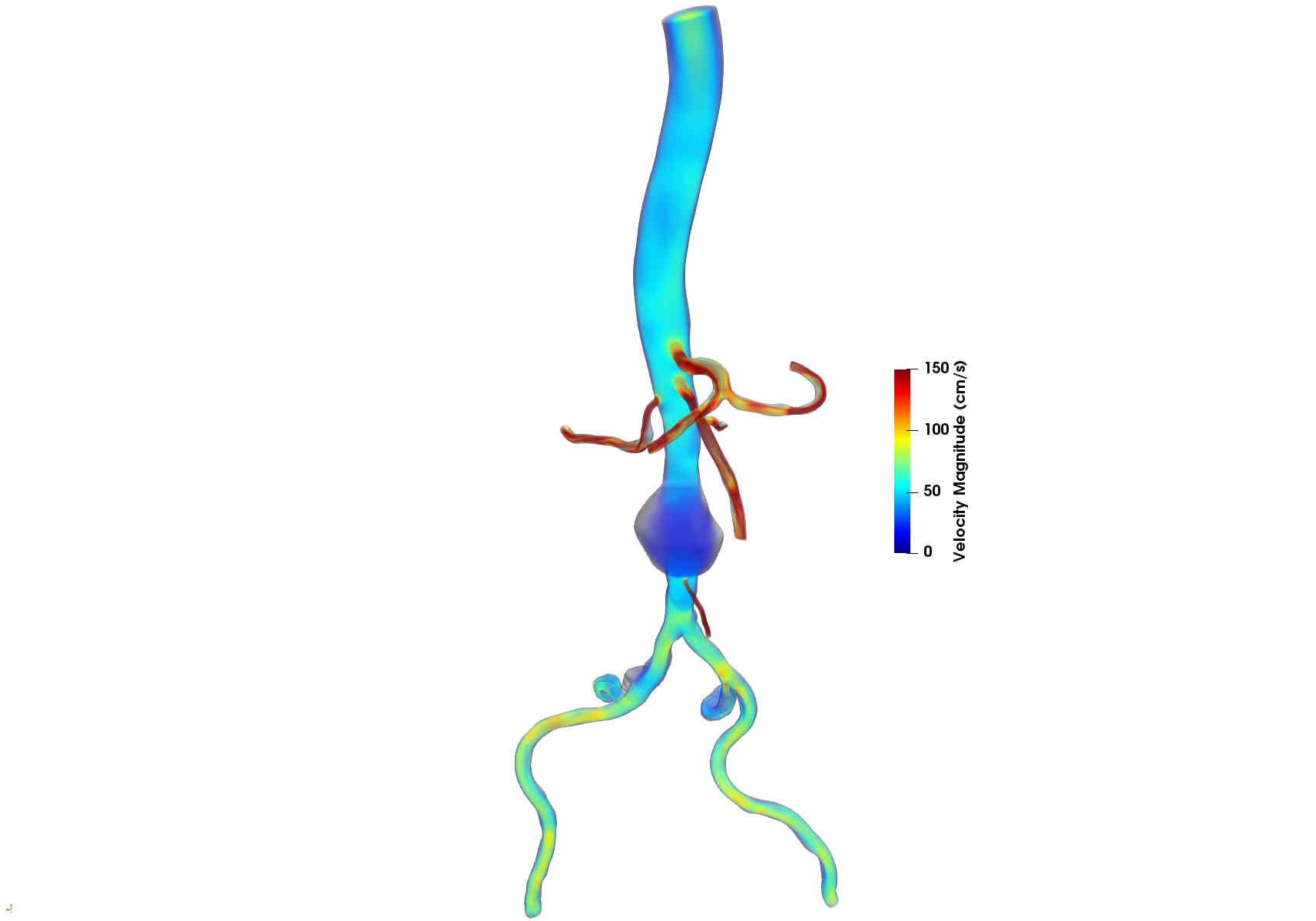} 
\end{tabular}
\caption{The blood pressure (left) and the velocity magnitude (right) at peak systole. } 
\label{fig:aaa_pressure_velocity}
\end{center}
\end{figure}

\begin{figure}
	\begin{center}
	\begin{tabular}{cc}
\includegraphics[angle=0, trim=500 0 600 0, clip=true, scale = 0.22]{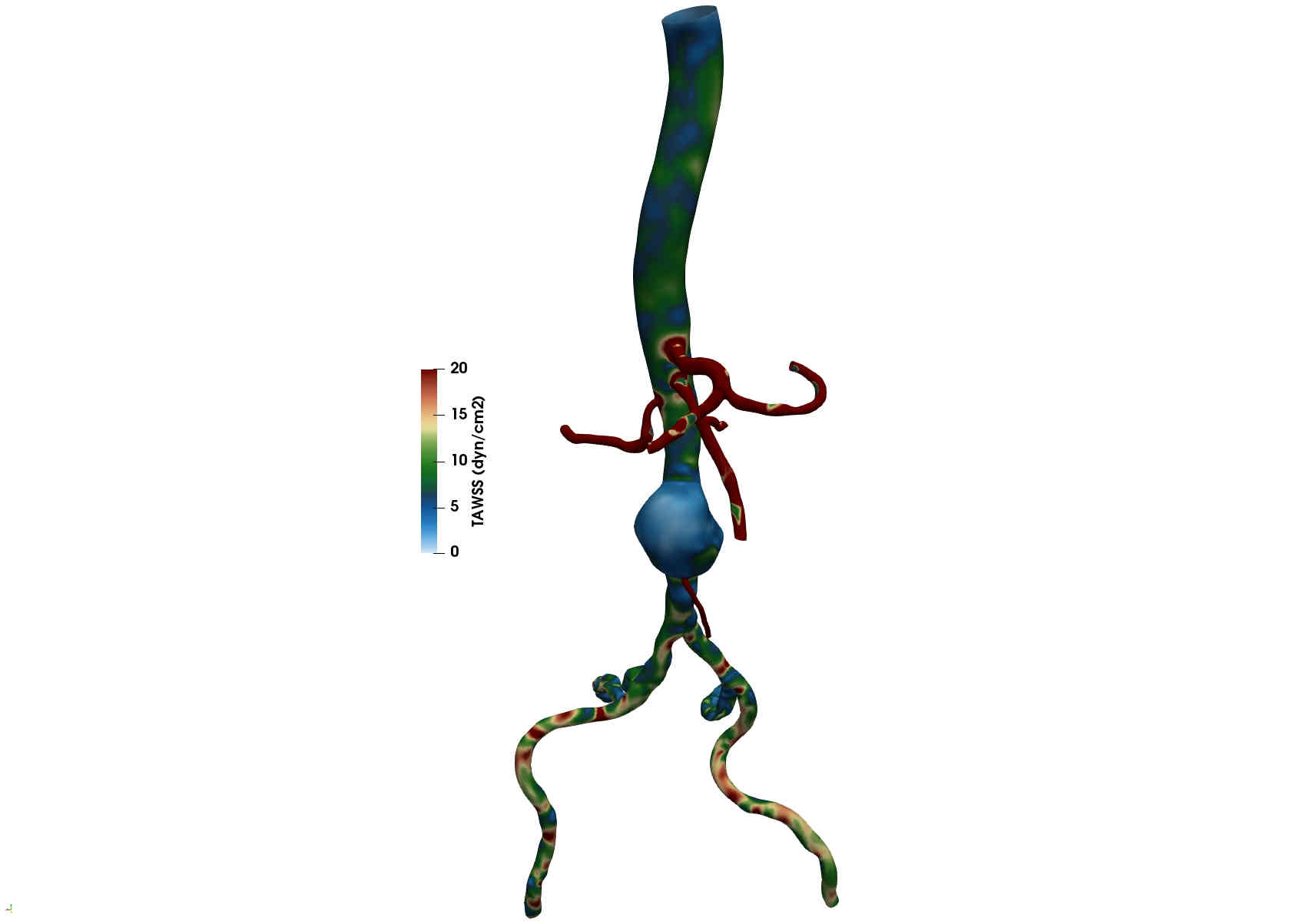} &
\includegraphics[angle=0, trim=600 0 450 0, clip=true, scale = 0.22]{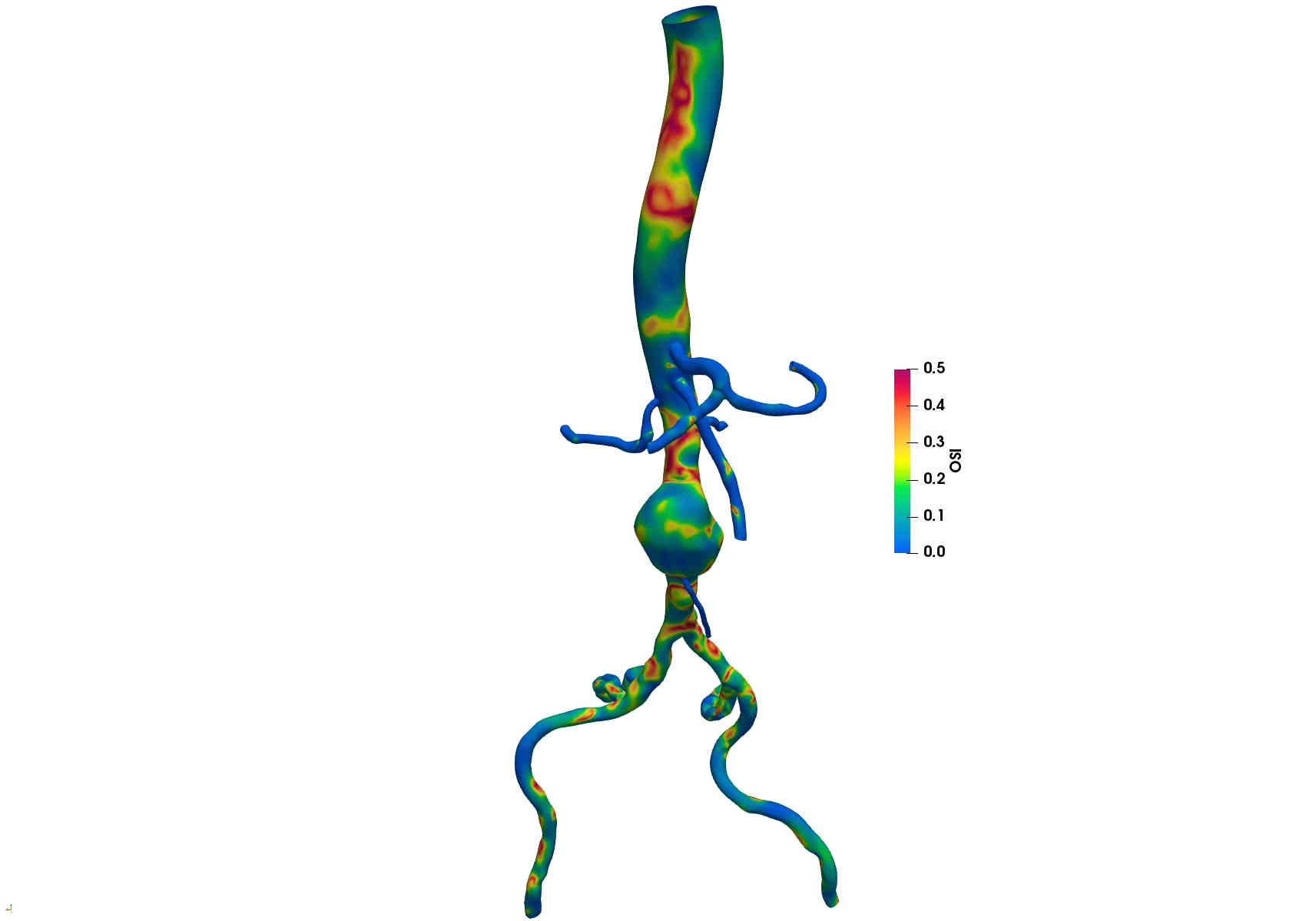} 
\end{tabular}
\caption{The time-averaged wall shear stress and the oscillatory shear index of the aortofemoral model.} 
\label{fig:aaa_tawss_osi}
\end{center}
\end{figure}

\begin{figure}
\begin{center}
	\begin{tabular}{cc}
\includegraphics[angle=0, trim=30 20 180 120, clip=true, scale = 0.083]{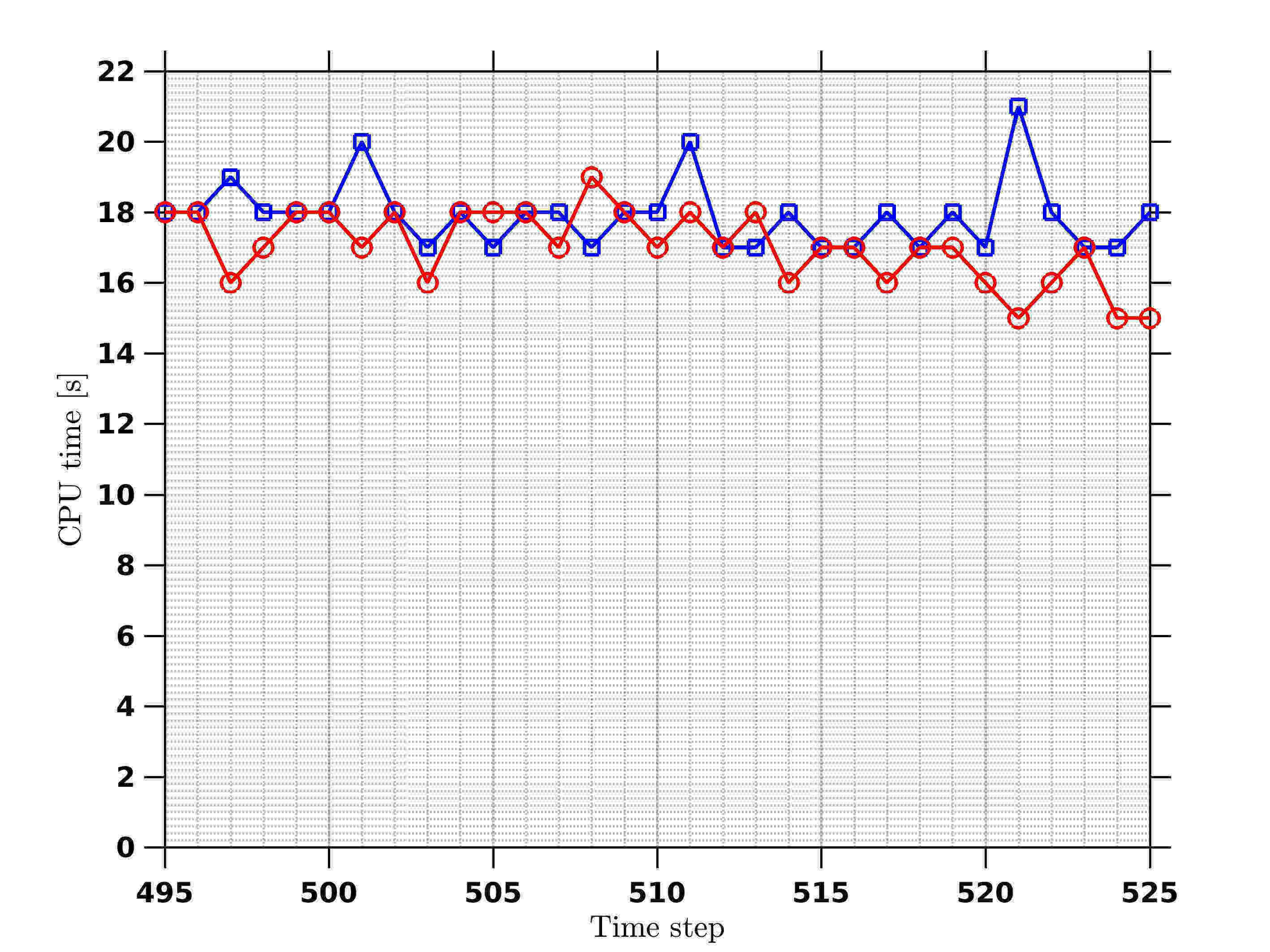} & 
\includegraphics[angle=0, trim=30 20 180 120, clip=true, scale = 0.083]{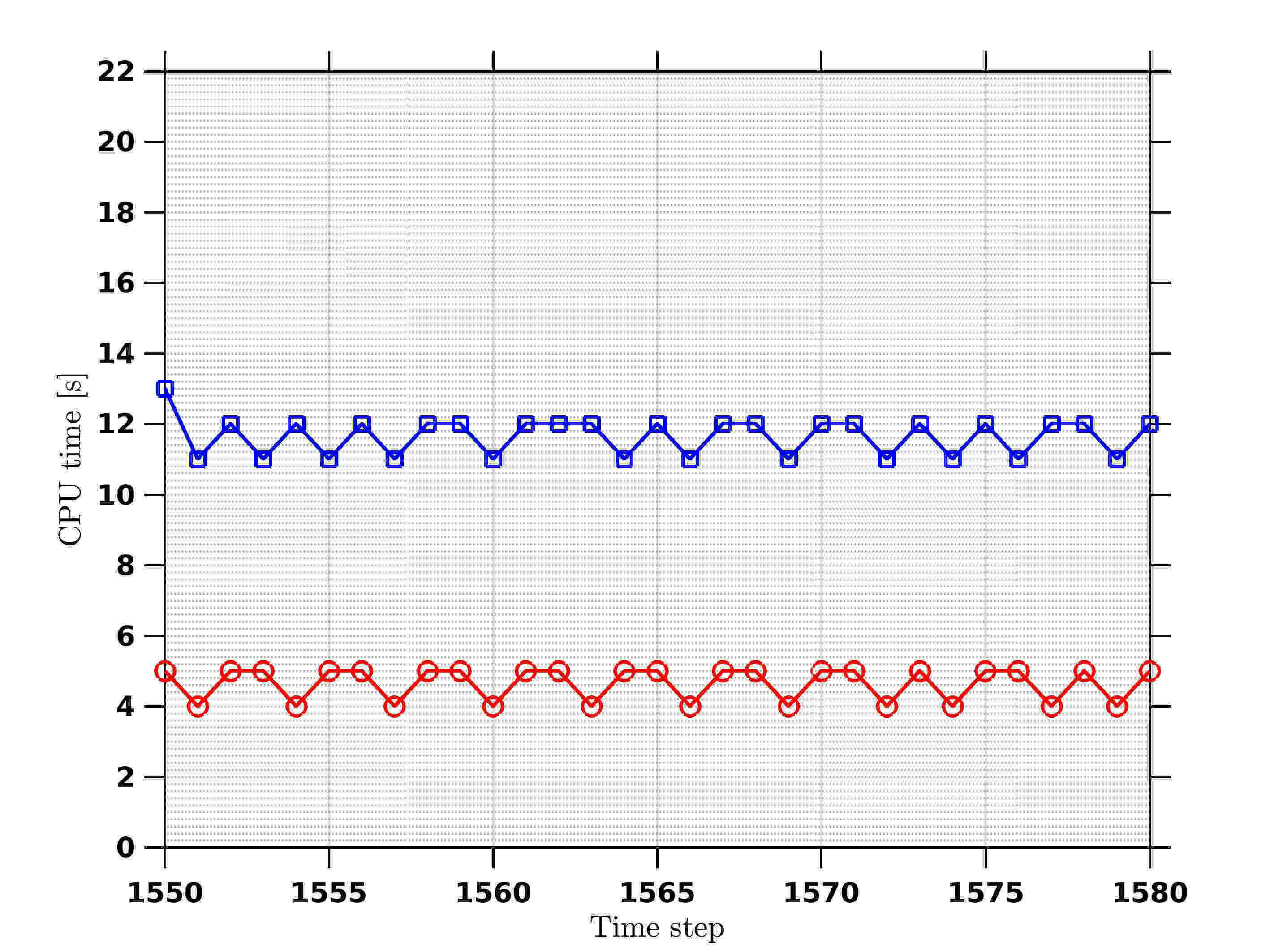}
\end{tabular}
\caption{Comparison of the CPU time for $30$ time steps near the peak systole (left) and mid diastole (right). The curves represents the CPU time spent in each time step for the proposed method (red) and \texttt{svSolver} (blue). During the systole stage, the proposed solution method requires $4$ Newton-Raphson iterations, while the \texttt{svSolver} reaches the maximum number of the Newton-Raphson iteration. During diastole, the proposed solution method converges in $2$ Newton-Raphson iterations, and the \texttt{svSolver} needs $12$ iterations.} 
\label{fig:sv0149_sv_pe_compare}
\end{center}
\end{figure}

\subsection{Patient-specific hemodynamic simulation}
\label{sec:patient_specific_examples}
In this section, we further examine the proposed method under physiologically realistic settings by considering two patient-specific examples: a diseased aortofemoral model and a healthy pulmonary model. The models and the clinical data are publicly available through the Cardiovascular and Pulmonary Model Repository \cite{AAA-model-url}. The two cases are selected to cover the systemic and pulmonary circulations. The meshes for both examples are generated by Meshsim \cite{simmetrix}, with boundary layer meshes near the arterial wall. We simulate the problems using the nested block preconditioner with $\delta^r = 10^{-2}$, $\delta^r_A = 10^{-3}$, $\delta^r_S = \delta^r_I =  10^{-2}$, and $\mathrm m_{\mathrm A} = \mathrm m_{\mathrm S}= \mathrm n^{\textup{max}}_{\mathrm A} = \mathrm n^{\textup{max}}_{\mathrm S} = 20$. Similar to the previous example, we use the Jacobi preconditioner for the inner and intermediate solvers associated with $\boldsymbol{\mathrm A}$. The same problems are also simulated with \texttt{svSolver}. In \texttt{svSolver}, the parameters adopt the default values except that the tolerances on momentum equations and the continuity equation are chosen to be $10^{-3}$ and $10^{-2}$, respectively. The above parameters are chosen to achieve the fastest performances for the two solvers separately. Since the linear systems in the two solvers are not identical, we only aim at comparing the overall algorithm efficiency, with the same accuracy for solving the nonlinear equations. 

\subsubsection{Aortofemoral model}
We first consider an aortofemoral model with an abdominal aortic aneurysm. The medical image and the volume rendering of the model are illustrated in Fig. \ref{fig:aortofemoral_image_and_mesh} (a). The volumetric flow rate used for prescribing the velocity on the inlet surface is illustrated in Fig. \ref{fig:aortofemoral_image_and_mesh} (b). The generated mesh consists of $8.80\times 10^6$ tetrahedral elements and $1.61\times 10^6$ vertices (Fig. \ref{fig:aortofemoral_image_and_mesh} (c)). The minimum element size is $\Delta x_{\textup{min}} = 9.21\times 10^{-3}$ cm, and the time step size is chosen as $\Delta t = 4.16\times 10^{-4}$ s. Notice that one cardiac cycle takes $0.832$ s, and thence it requires $2000$ time steps for simulating one cardiac cycle. The simulations are performed with $288$ CPUs and simulated for $12$ cardiac cycles. The instantaneous pressure and velocity at the peak systole are depicted in Fig. \ref{fig:aaa_pressure_velocity}; the time-averaged wall shear stress and the oscillatory shear index are illustrated in Fig. \ref{fig:aaa_tawss_osi}. The performances of the proposed solution method and the \texttt{svSolver} are monitored for $30$ time steps around peak systole and mid diastole (see Fig. \ref{fig:sv0149_sv_pe_compare}). For one cardiac cycle, the proposed solution method takes $14803$ s and \texttt{svSolver} takes $28200$ s.

\begin{figure}
	\begin{center}
	\begin{tabular}{cc}
\includegraphics[angle=0, trim=50 30 50 30, clip=true, scale = 0.16]{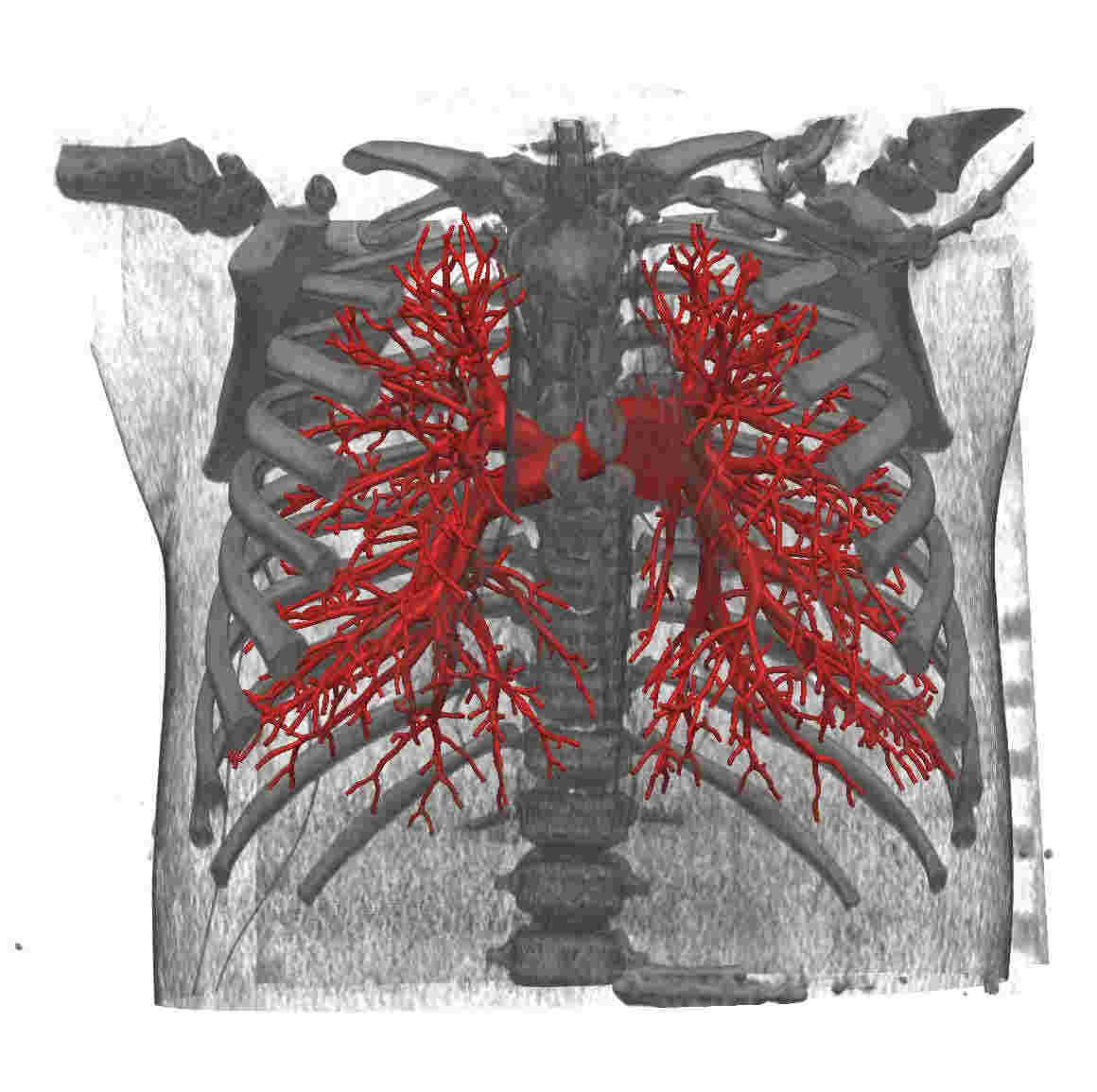} & 
\includegraphics[angle=0, trim=130 0 230 0, clip=true, scale = 0.08]{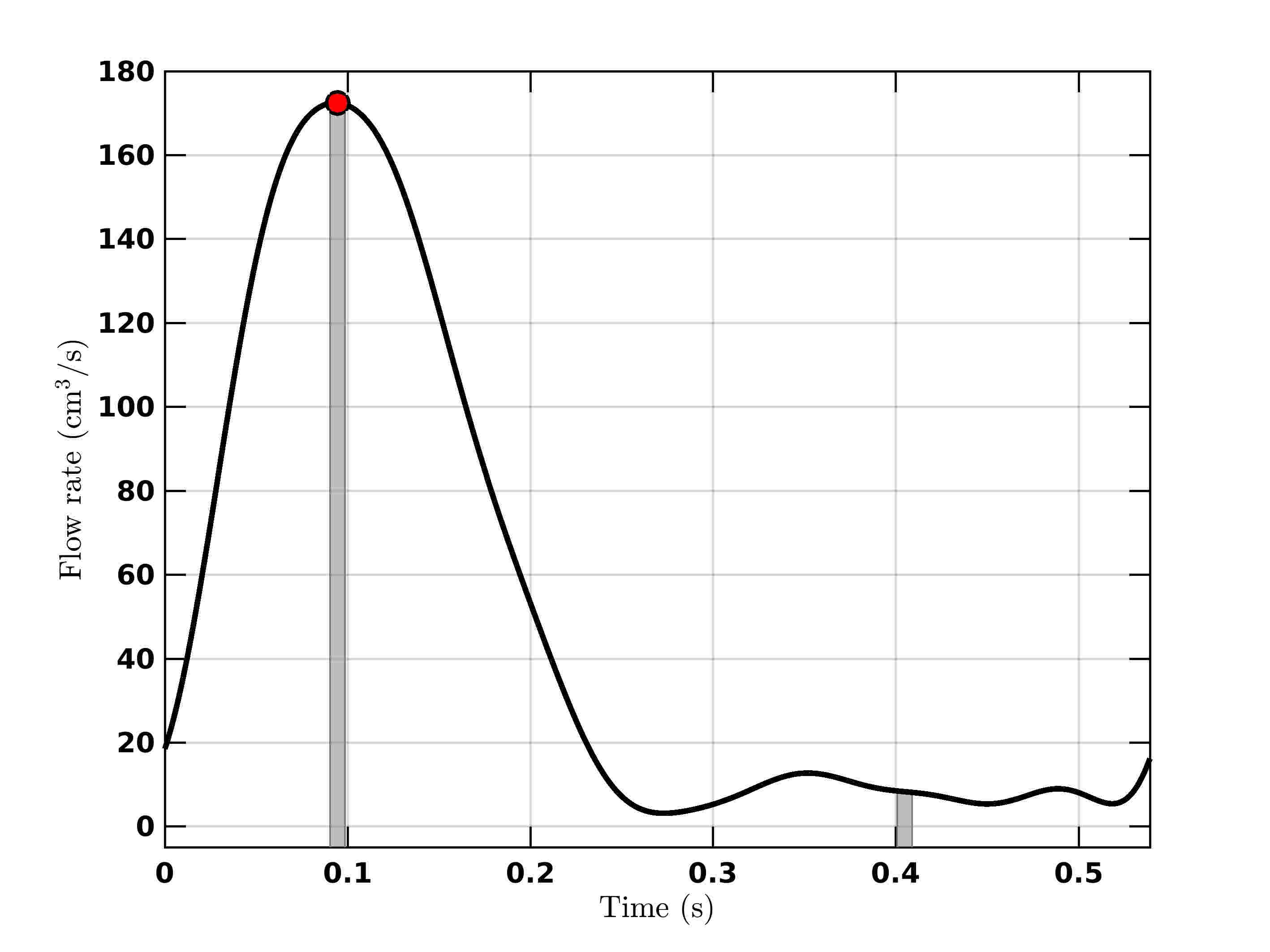}
\\
\multicolumn{2}{c}{ \includegraphics[angle=0, trim=500 500 150 500, clip=true, scale = 0.025]{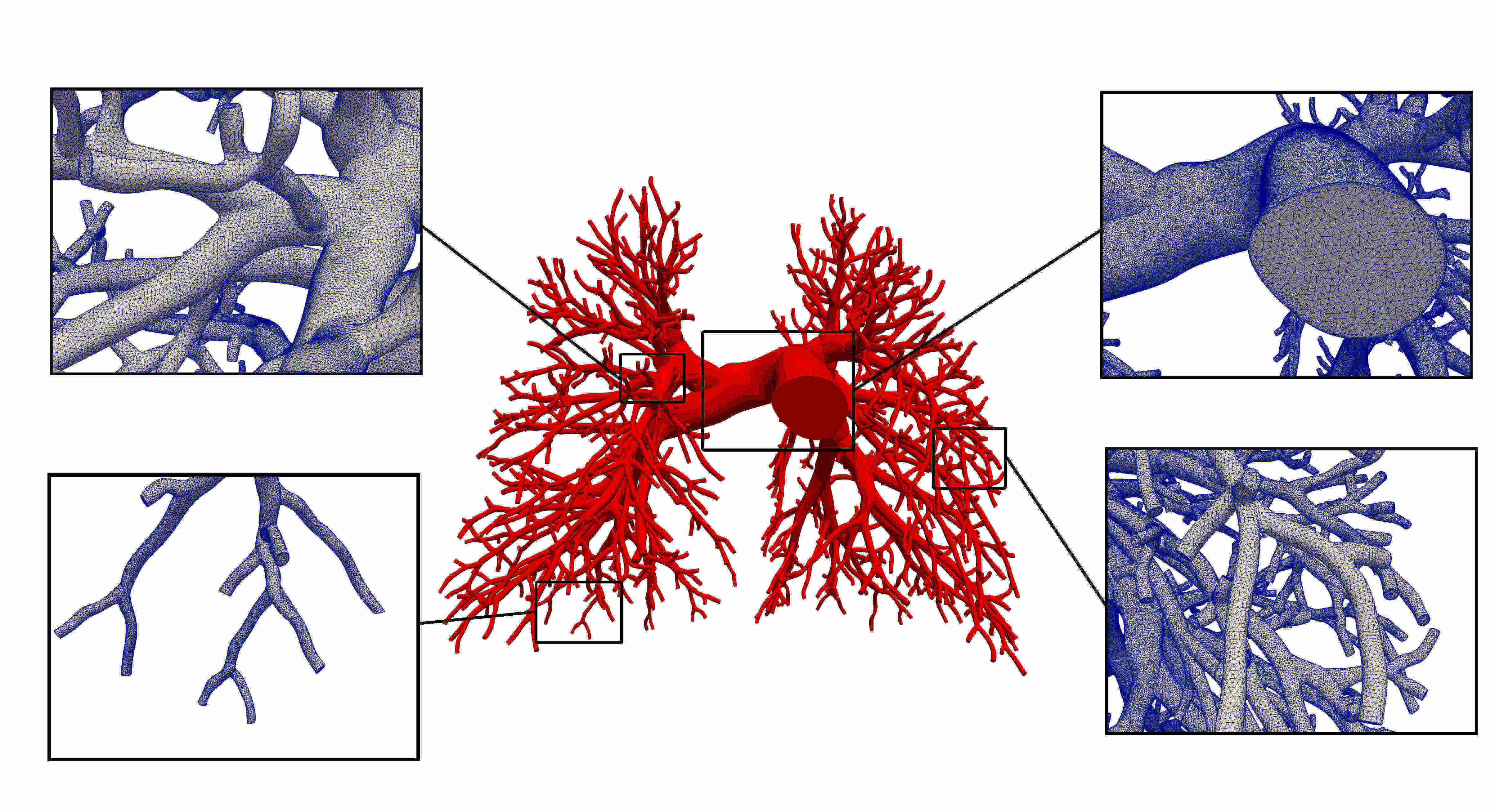} }
\end{tabular}
\caption{(a) Volume rendering of the medical data. (b) Volumetric flow waveform used to prescribe the velocity on the inlet surface. The red dot represents peak systole. The grey regions represent the periods of time where we monitor the CPU time for the two different solvers in Fig. \ref{fig:su0273_sv_pe_compare}. (c) The mesh of the pulmonary model.} 
\label{fig:pulmonary_image_and_mesh}
\end{center}
\end{figure}

\begin{figure}
	\begin{center}
	\begin{tabular}{c}
\includegraphics[angle=0, trim=50 0 430 0, clip=true, scale = 0.24]{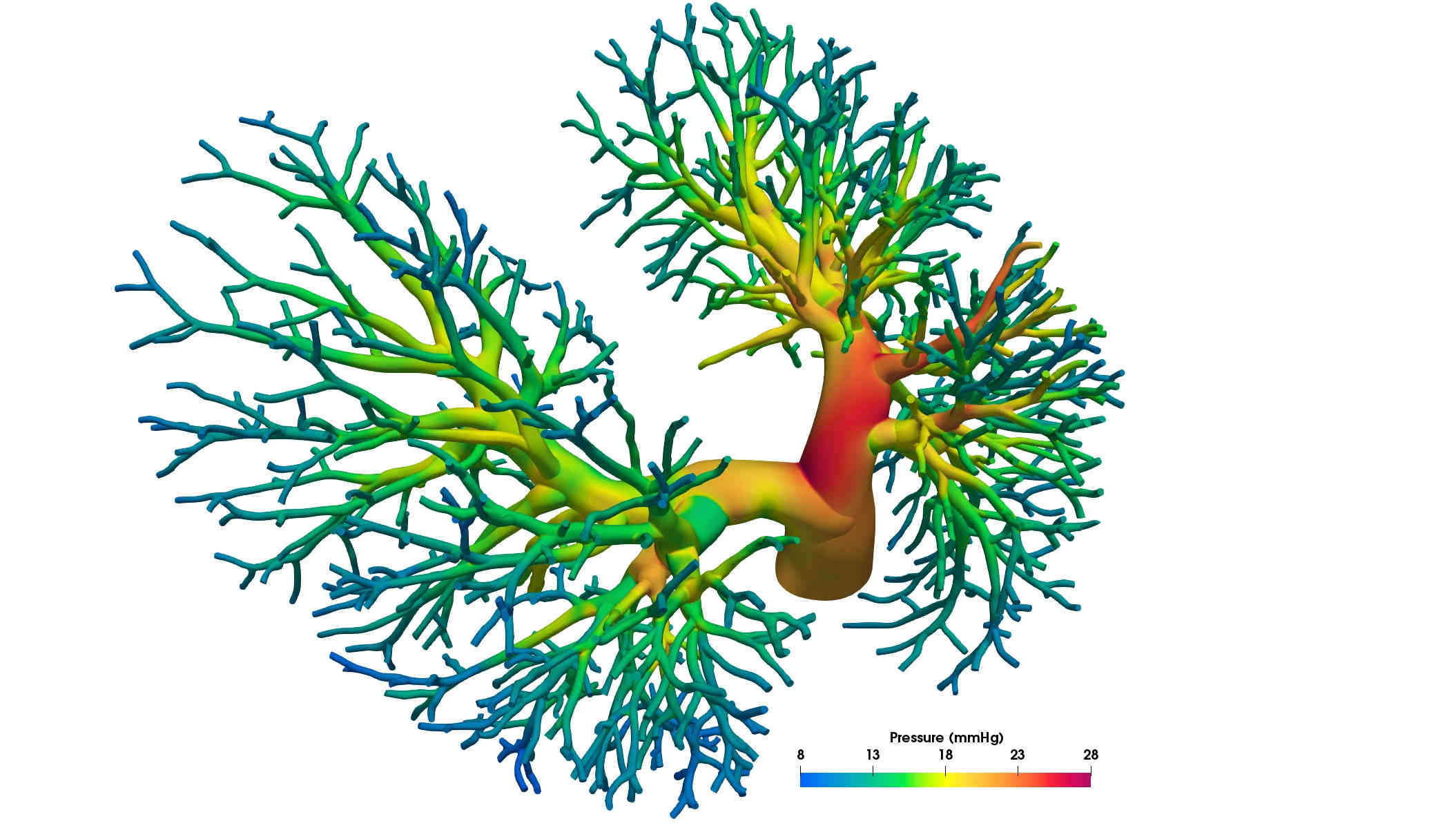} \\
\includegraphics[angle=0, trim=50 0 430 0, clip=true, scale = 0.24]{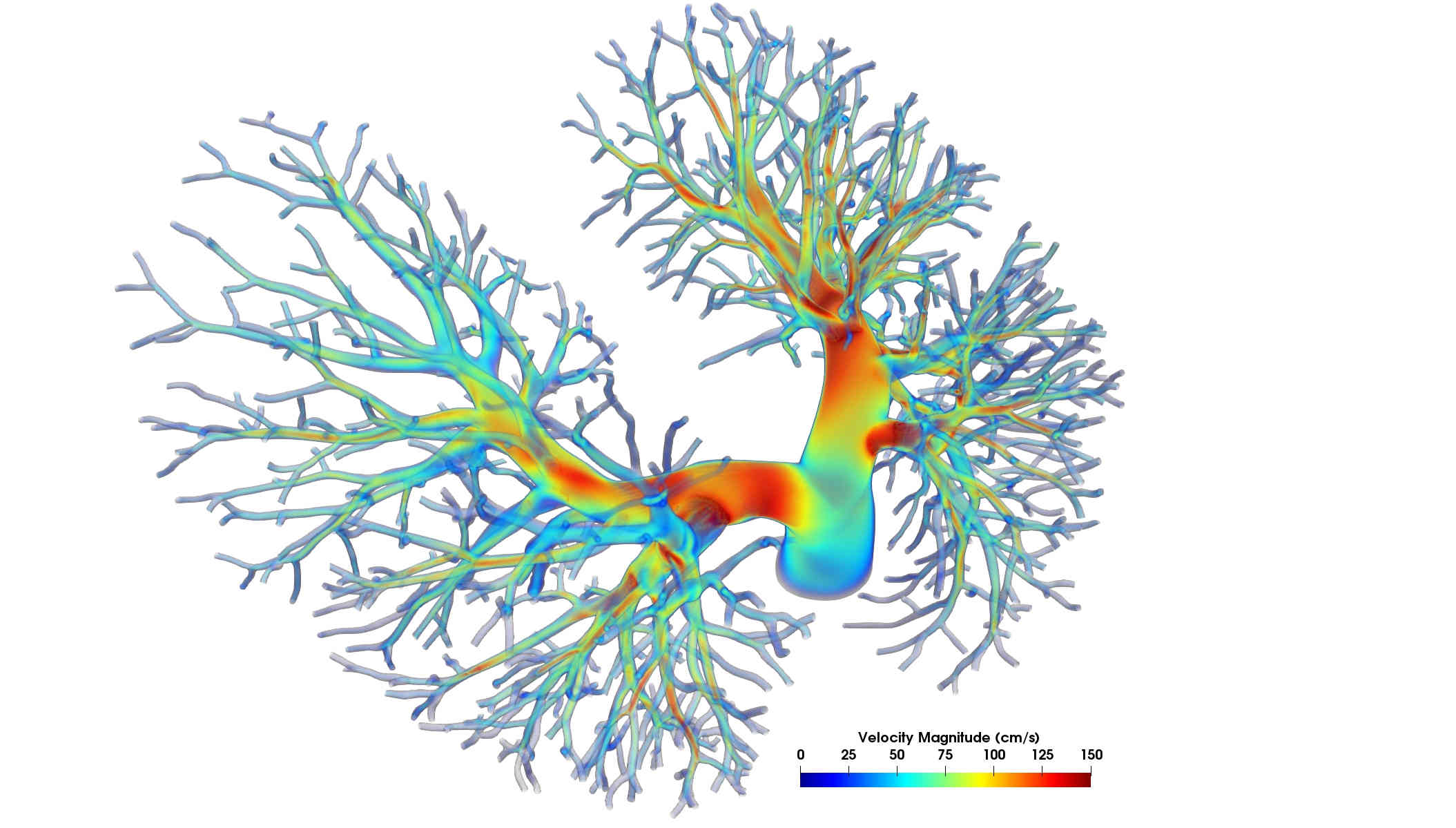}
\end{tabular}
\caption{The blood pressure and the volume rendering of the velocity magnitude at peak systole.} 
\label{fig:pulmonary_t35}
\end{center}
\end{figure}

\begin{figure}
	\begin{center}
	\begin{tabular}{c}
\includegraphics[angle=0, trim=50 0 430 0, clip=true, scale = 0.24]{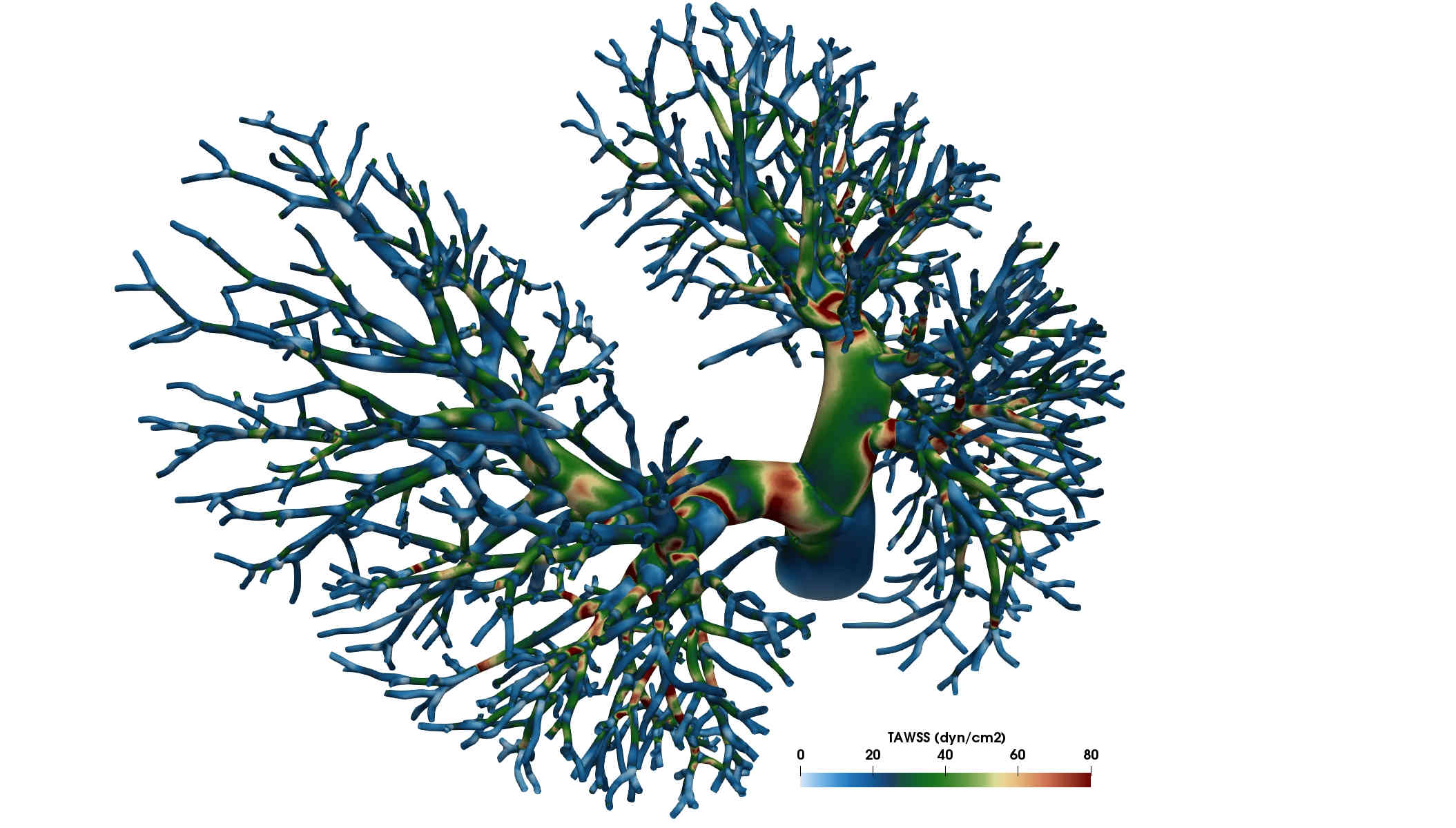} \\
\includegraphics[angle=0, trim=50 0 430 0, clip=true, scale = 0.24]{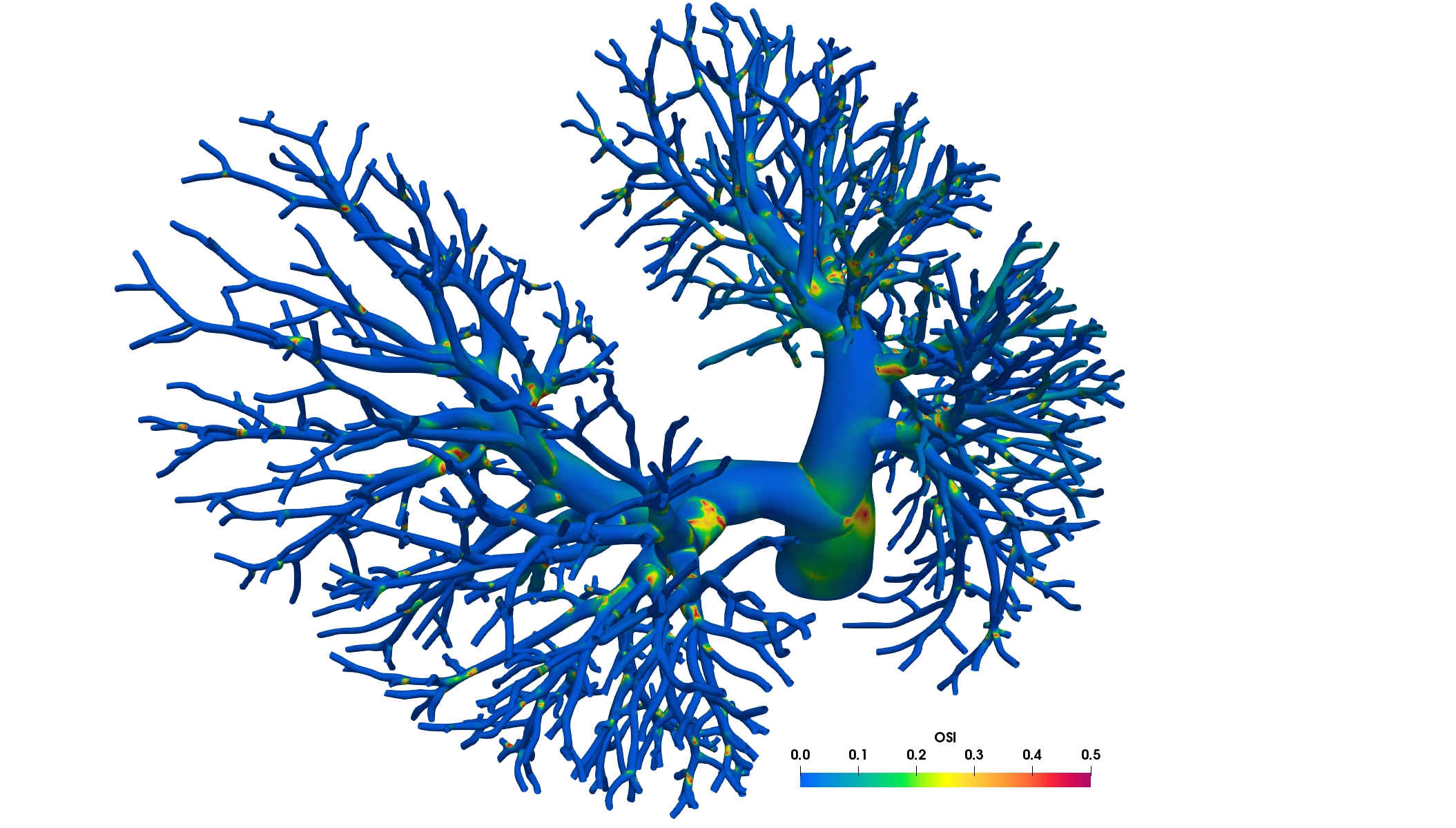}
\end{tabular}
\caption{The time-averaged wall shear stress and the oscillatory shear index of the pulmonary model.} 
\label{fig:pulmonary_tawss_osi}
\end{center}
\end{figure}

\begin{figure}
\begin{center}
	\begin{tabular}{cc}
\includegraphics[angle=0, trim=30 20 180 120, clip=true, scale = 0.083]{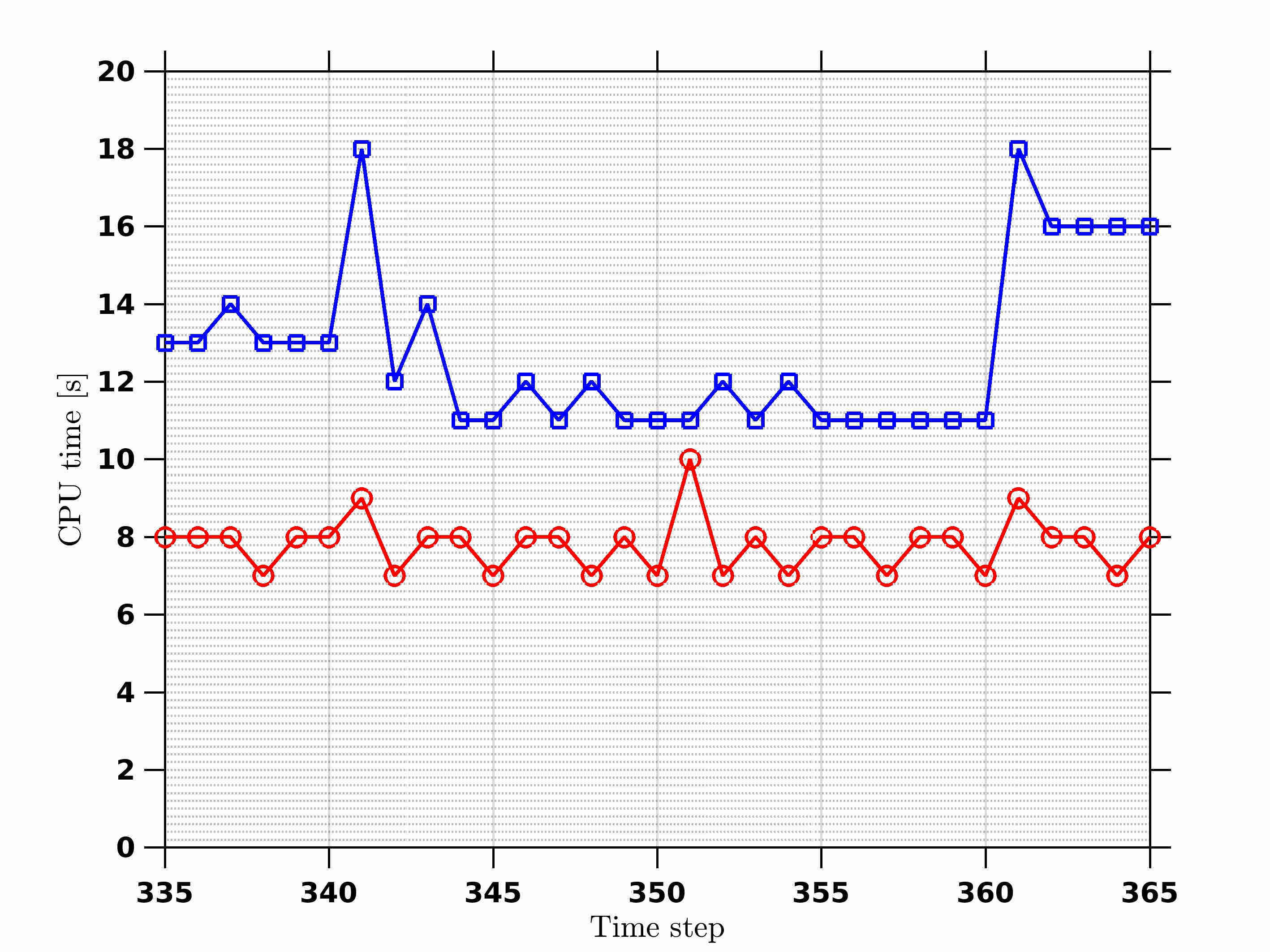} & 
\includegraphics[angle=0, trim=30 20 180 120, clip=true, scale = 0.083]{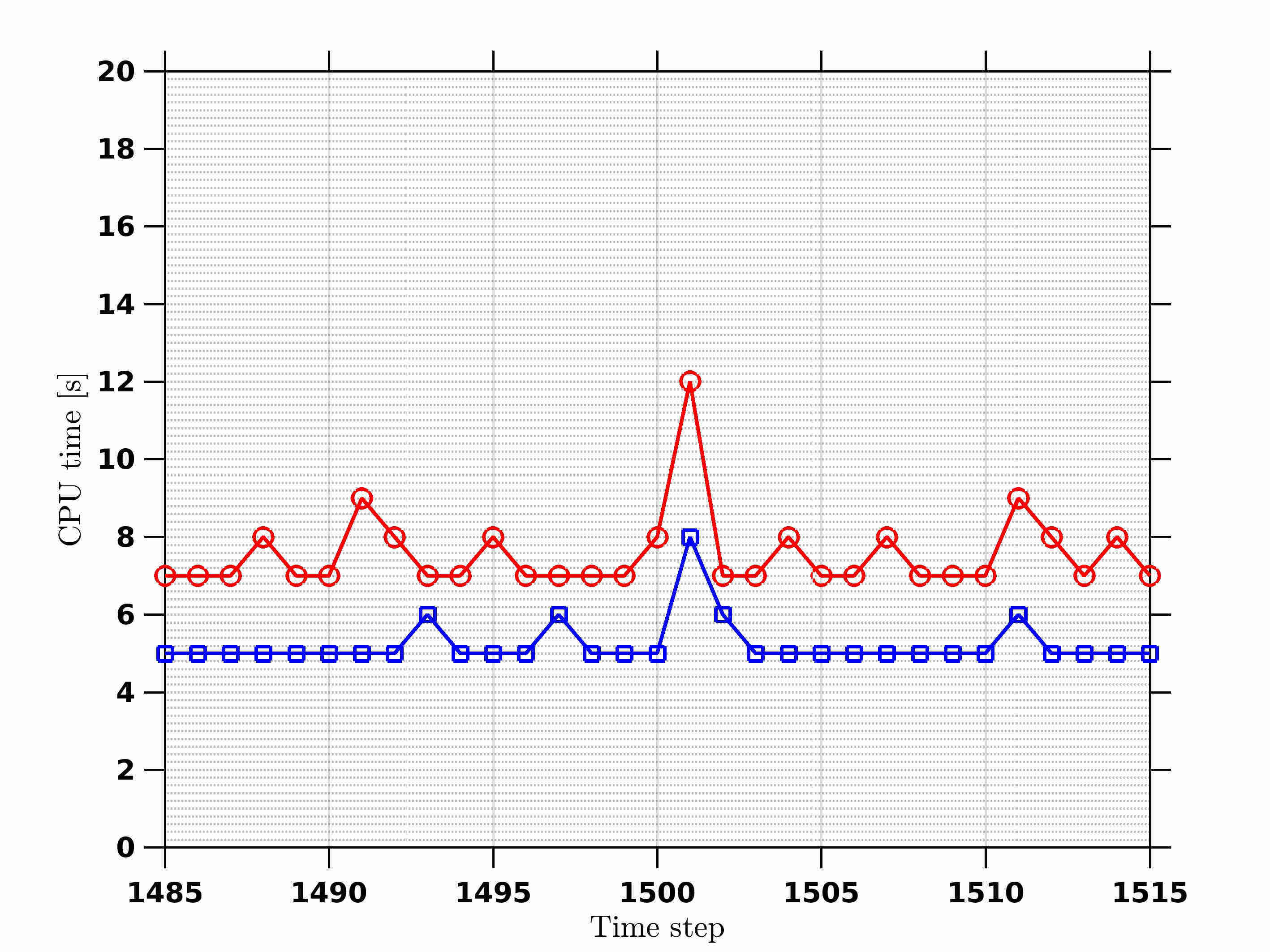}
\end{tabular}
\caption{Comparison of the CPU time for $30$ time steps near the peak systole (left) and mid diastole (right). The curves represent the CPU time spent in each time step for the proposed method (red) and the \texttt{svSolver} (blue). During systole, the proposed solution method requires $3$ Newton-Raphson iterations, while the \texttt{svSolver} requires $14$ iterations. During diastole, the proposed solution method converges in $2$ Newton-Raphson iterations, and the \texttt{svSolver} needs $4$ iterations.} 
\label{fig:su0273_sv_pe_compare}
\end{center}
\end{figure}

\subsubsection{Pulmonary model}
\label{sec:pulmonary_model}
The pulmonary circulation has lower pressures than the systemic circulation, but with the same amount of flow going through an extensive and concentrated tree of branching arteries. To demonstrate the performance of the proposed method in the pulmonary circulation, a model was built from a healthy 20-month-old male, consisting of 772 outlets. This mesh consists of $2.61\times 10^7$ tetrahedral elements and $4.95\times 10^6$ vertices. Correspondingly, there are about $20$ million unknowns in the resulting linear system of equations. The minimum element size is $\Delta x_{\textup{min}} = 5.16\times 10^{-3}$ cm. The time step size is fixed to be $\Delta t = 2.6995\times 10^{-4}$ s. The medical image, volumetric flow rate, and the detailed views of the mesh are illustrated in Fig. \ref{fig:pulmonary_image_and_mesh}. Notice that one cardiac cycle takes $0.5399$ s and it requires $2000$ steps in the time integration. The simulations are performed with $720$ CPUs and simulated for $3$ cardiac cycles. The pressure and velocity magnitude at the peak systole in the last cycle are depicted in Fig. \ref{fig:pulmonary_t35}; the time-averaged wall shear stress and the oscillatory shear index are illustrated in Fig. \ref{fig:pulmonary_tawss_osi}. Performance of the proposed solution method and the \texttt{svSolver} is monitored for $30$ time steps around peak systole and mid diastole (see Fig. \ref{fig:su0273_sv_pe_compare}). We notice that \texttt{svSolver} is faster than the proposed algorithm during the diastolic phase in this example. For one full cardiac cycle, the proposed method takes $14500$ s and \texttt{svSolver} takes $19500$ s.

\begin{figure}
\begin{center}
	\begin{tabular}{cc}
\includegraphics[angle=0, trim=30 20 200 120, clip=true, scale = 0.08]{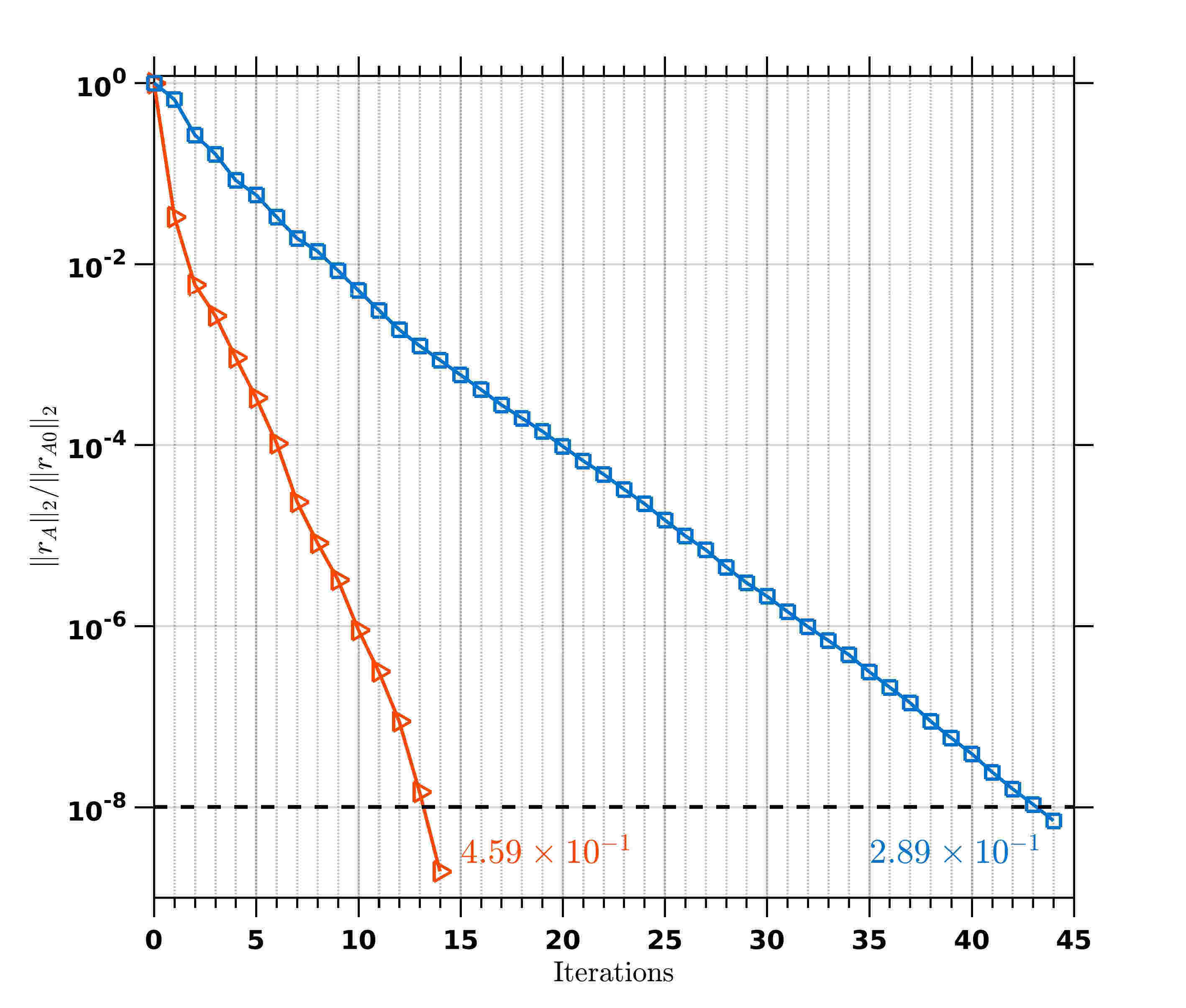} & 
\includegraphics[angle=0, trim=30 20 200 120, clip=true, scale = 0.08]{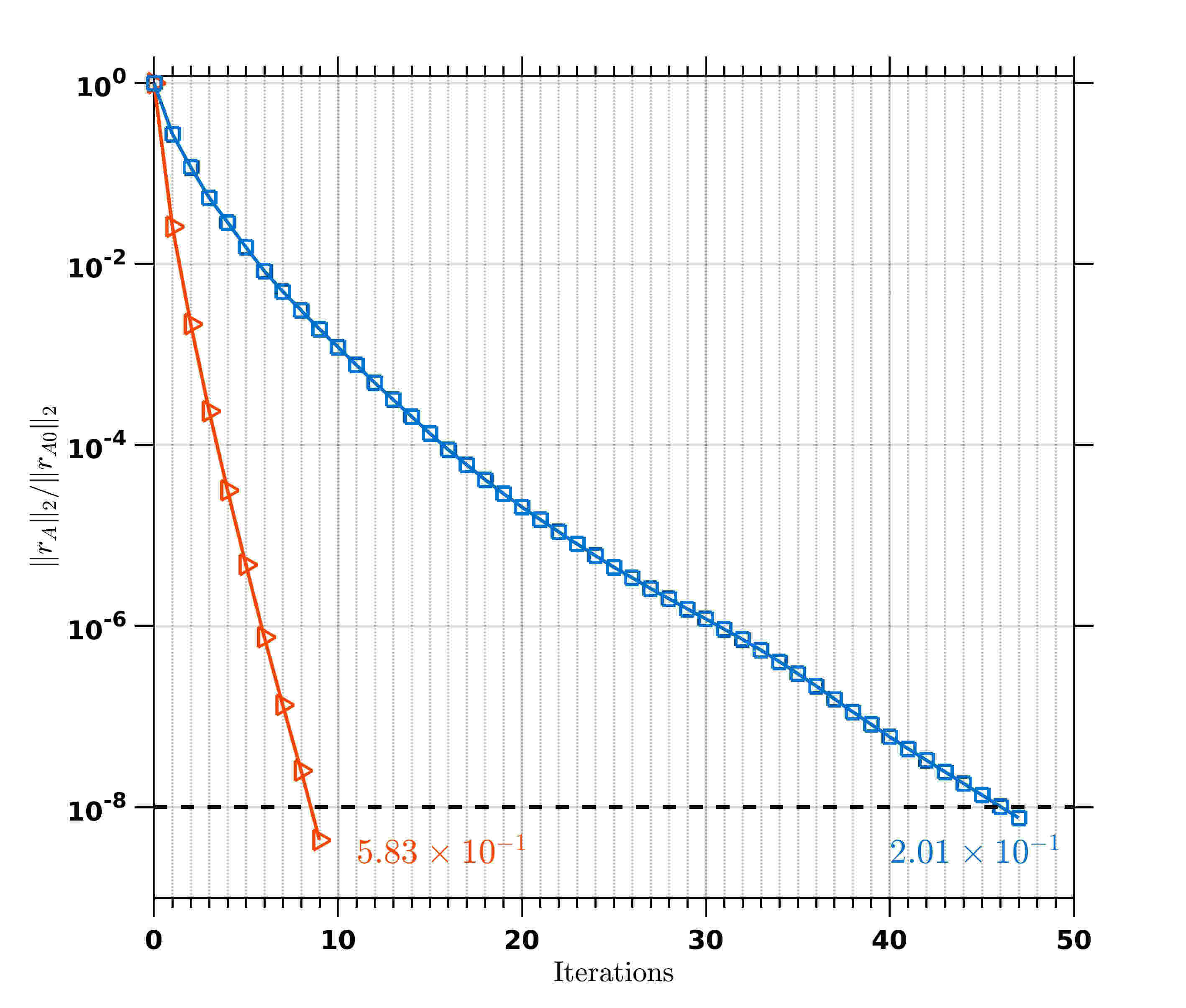} \\
\multicolumn{2}{c}{ \includegraphics[angle=0, trim=2360 100 2180 310, clip=true, scale = 0.12]{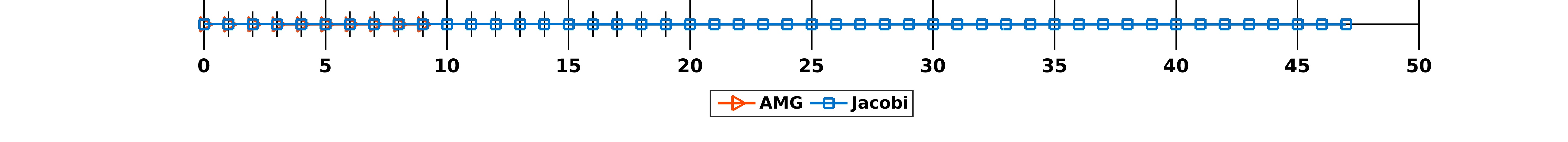} }
\end{tabular}
\caption{Convergence history for solving with $\boldsymbol{\mathrm A}$ with the AMG and Jacobi preconditioner for the aortofemoral model (left) and the pulmonary model (right). The numbers indicate the CPU time for solving the equations in the unit of seconds.} 
\label{fig:sv0149_amg_jacobi_compare}
\end{center}
\end{figure}

\subsubsection{Additional discussion}
\label{sec:additional_discussion}
As highlighted above, the inner solver plays a critical role in balancing the overall algorithm robustness and efficiency. The inner solver and the two solvers for \eqref{eq:seg_sol_int_disp} and \eqref{eq:seg_sol_disp} in the intermediate level are all associated with the block matrix $\boldsymbol{\mathrm A}$, and in this work, we naturally apply the same solution strategy for these equations. It has been suggested to use AMG in the design of traditional block preconditioners \cite{Elman2003,Burstedde2009}. Indeed, using AMG gives the most robust and scalable performance. However, the expensive setup phase of the AMG preconditioner results in a trade-off between scalability and efficiency. This can be particularly relevant when one simulates medium-sized problems when a scalable algorithm is not necessarily the most efficient one. For the problems considered in this work, we found that the AMG preconditioner is a very expensive option for $\boldsymbol{\mathrm A}$, and using the Jacobi preconditioner is often quite efficient, even for the pulmonary model with around $20$ million unknowns. In Fig. \ref{fig:sv0149_amg_jacobi_compare}, we tested the AMG and Jacobi preconditioners for solving  $\boldsymbol{\mathrm A}$ with $\delta^r_A = 10^{-8}$ using both patient-specific models. Although the AMG preconditioner requires fewer iterations, it costs almost twice the amount of CPU time compared to the Jacobi preconditioner. This convergence pattern is typical over iterations. Notice that the value of $\delta^r_A$ rarely reaches $10^{-8}$ in practice. For looser tolerances, the advantage of the Jacobi preconditioner is more pronounced. This is the reason we set $\delta^r_A=10^{-3}$ and $\mathrm m_{\mathrm A} = \mathrm n^{\textup{max}}_{\mathrm A} = 20$ for the patient-specific simulations. 

On the other hand we solve \eqref{eq:seg_sol_pres} with AMG because this step directly associates with the number of inner solver calls. We need a fast reduction rate when solving with the Schur complement. Based on our experience, AMG outperforms other options in this regard. We also set $\mathrm m_{\mathrm S} = \mathrm n^{\textup{max}}_{\mathrm S} = 20$ as a strategy to avoid spending too much time on solving with the Schur complement. In fact, we observe that the outer solver for $\boldsymbol{\mathcal A}$ typically converges in $1$ or $2$ FGMRES iterations in the numerical simulations presented in this section. This suggests that the settings in the intermediate and inner levels may be further released.

\section{Conclusion}
In this work, we designed and applied the three-level nested block preconditioning technique to the incompressible Navier-Stokes equations with reduced models coupled on the boundaries. The proposed preconditioning technique combines the merits of the SCR procedure and the conventional physics-based block preconditioners. In particular, it uses an algorithm, rather than an explicit algebraic form, to approximate the behavior of the Schur complement. In comparison with conventional block preconditioners, this approach is well-suited to problems involving multiscale and multiphysics coupling. With a proper setting of the intermediate and inner solvers, the proposed method is demonstrated to attain more robust and more efficient performance than the existing technologies, especially for hemodynamic simulations. One attribute of the nested block preconditioner is the proliferation of options for the intermediate and inner solvers. The choices of the preconditioners and stopping criteria at the intermediate and inner levels critically determine the actual overall performance of the algorithm. For a specific problem, it is often advisable to determine the sweet spot of the options to balance robustness and efficiency. In our study, we noticed that using the light-weight Jacobi preconditioner for $\boldsymbol{\mathrm A}$ always leads to highly efficient performances. This is somewhat in contrast to the common practices in the traditional block preconditioner, where the AMG preconditioner is recommended for $\boldsymbol{\mathrm A}$. In future work, the options defining the intermediate and inner solvers need to be extensively explored and compared for performance optimization purposes. Also, a crucial next step is to apply the iterative solution method for FSI problems within the unified continuum modeling framework \cite{Liu2018}.

\section*{Acknowledgements}
We thank A.W. Bergersen for helpful discussions on the FDA idealized medical device benchmark. This work is supported by the National Institutes of Health under the award numbers 1R01HL121754, 1R01HL123689, and R01EB01830204, the Department of Defense award W81XWH1810518, and the computational resources from the Stanford Research Computing Center and the Extreme Science and Engineering Discovery Environment (XSEDE) supported by the NSF grant ACI-1053575.

\bibliographystyle{elsarticle-num} 

\bibliography{solver_ns}

\begin{thebibliography}{10}
\expandafter\ifx\csname url\endcsname\relax
  \def\url#1{\texttt{#1}}\fi
\expandafter\ifx\csname urlprefix\endcsname\relax\def\urlprefix{URL }\fi
\expandafter\ifx\csname href\endcsname\relax
  \def\href#1#2{#2} \def\path#1{#1}\fi

\bibitem{Bao2014}
G.~Bao, Y.~Bazilevs, J.~Chung, P.~Decuzzi, H.~Espinosa, M.~Ferrari, H.~Gao,
  S.~Hossain, T.~Hughes, R.~Kamm, W.~Liu, A.~Marsden, B.~Schrefler, {USNCTAM}
  perspectives on mechanics in medicine, Journal of The Royal Society Interface
  11~(97).

\bibitem{Figueroa2017}
C.~Figueroa, C.~Taylor, A.~Marsden, Blood flow, Encyclopedia of Computational
  Mechanics (2017) 1--31.

\bibitem{Marsden2013a}
A.~Marsden, Optimization in cardiovascular modeling, Annual Review of Fluid
  Mechanics 46 (2013) 519--546.

\bibitem{top-500-list-site}
Top500 lists, \url{https://www.top500.org/}.

\bibitem{Taylor2013}
C.~Taylor, T.~Fonte, J.~Min, Computational fluid dynamics applied to cardiac
  computed tomography for noninvasive quantification of fractional flow
  reserve: scientific basis, Journal of the American College of Cardiology 61
  (2013) 2233--2241.

\bibitem{Zarins2013}
C.~Zarins, C.~Taylor, J.~Min, Computed fractional flow reserve ({FFT CT})
  derived from coronary ct angiography, Journal of cardiovascular translational
  research 6 (2013) 708--714.

\bibitem{Hwang2005}
F.~Hwang, X.~Cai, A parallel nonlinear additive {S}chwarz preconditioned
  inexact {N}ewton algorithm for incompressible {N}avier-{S}tokes equations,
  Journal of Computational Physics 204 (2005) 666--691.

\bibitem{Kong2019}
F.~Kong, V.~Kheyfets, E.~Finol, X.~Cai, Simulation of unsteady blood flows in a
  patient-specific compliant pulmonary artery with a highly parallel
  monolithically coupled fluid-structure interaction algorithm, International
  Journal for Numerical Methods in Biomedical Engineering 35 (2019) e3208.

\bibitem{Wu2014}
Y.~Wu, X.~Cai, A fully implicit domain decomposition based {ALE} framework for
  three-dimensional fluid–structure interaction with application in blood
  flow computation, Journal of Computational Physics 258 (2014) 524--537.

\bibitem{Elman2008}
H.~Elman, V.~Howle, J.~Shadid, R.~Shuttleworth, R.~Tuminaro, A taxonomy and
  comparison of parallel block multi-level preconditioners for the
  incompressible {N}avier-{S}tokes equations, Journal of Computational Physics
  227 (2008) 1790--1808.

\bibitem{Cyr2012}
E.~Cyr, J.~Shadid, R.~Tuminaro, Stabilization and scalable block
  preconditioning for the {N}avier-{S}tokes equations, Journal of Computational
  Physics 231 (2012) 345--363.

\bibitem{Washio2005}
T.~Washio, T.~Hisada, H.~Watanabe, T.~Tezduyar, A robust preconditioner for
  fluid-structure interaction problems, Computer Methods in Applied Mechanics
  and Engineering 194 (2005) 4027--4047.

\bibitem{Tuminaro2002}
R.~Tuminaro, C.~Tong, J.~Shadid, K.~Devine, D.~Day, On a multilevel
  preconditioning module for unstructured mesh {K}rylov solvers: two-level
  {S}chwarz, Communications in numerical methods in engineering 18 (2002)
  383--389.

\bibitem{Shadid2005}
J.~Shadid, R.~Tuminaro, K.~Devine, G.~Hennigan, P.~Lin, Performance of fully
  coupled domain decomposition preconditioners for finite element
  transport/reaction simulations, Journal of Computational Physics 205 (2005)
  24--47.

\bibitem{Benzi2005}
M.~Benzi, G.~Golub, J.~Liesen, Numerical solution of saddle point problems,
  Acta Numerica 14 (2005) 1--137.

\bibitem{Chorin1968}
A.~Chorin, Numerical solution of the {N}avier-{S}tokes equations, Mathematics
  of computation 22 (1968) 745--762.

\bibitem{Kim1985}
J.~Kim, P.~Moin, Application of a fractional-step method to incompressible
  {N}avier-{S}tokes equations, Journal of Computational Physics 59 (1985)
  308--323.

\bibitem{Gresho1998}
P.~Gresho, R.~Sani, Incompressible flow and the finite element method. {V}olume
  1: {A}dvection-diffusion and isothermal laminar flow, John Wiley and Sons,
  Inc., New York, NY (United States), 1998.

\bibitem{Jansen2000}
K.~Jansen, C.~Whiting, G.~Hulbert, A generalized-$\alpha$ method for
  integrating the filtered {N}avier-{S}tokes equations with a stabilized finite
  element method, Computer Methods in Applied Mechanics and Engineering 190
  (2000) 305--319.

\bibitem{Guermond2006}
J.~Guermond, P.~Minev, J.~Shen, An overview of projection methods for
  incompressible flows, Computer Methods in Applied Mechanics and Engineering
  195 (2006) 6011--6045.

\bibitem{Elman2003}
H.~Elman, V.~Howle, J.~Shadid, R.~Tuminaro, A parallel block multi-level
  preconditioner for the 3{D} incompressible {N}avier-{S}tokes equations,
  Journal of Computational Physics 187 (2003) 504--523.

\bibitem{Deparis2014}
S.~Deparis, G.~Grandperrin, A.~Quarteroni, Parallel preconditioners for the
  unsteady {N}avier-{S}tokes equations and applications to hemodynamics
  simulations, Computer \& Fluids 92 (2014) 253--273.

\bibitem{Elman2006}
H.~Elman, V.~Howle, J.~Shadid, R.~Shuttleworth, R.~Tuminaro, Block
  preconditioners based on approximate commutators, SIAM Journal on Scientific
  Computing 27 (2006) 1651--1668.

\bibitem{Kay2002}
D.~Kay, D.~Loghin, A.~Wathen, A preconditioner for the steady-state
  {N}avier-{S}tokes equations, SIAM Journal on Scientific Computing 24 (2002)
  237--256.

\bibitem{Shadid2016}
J.~Shadid, P.~Pawlowski, E.~Cyr, R.~Tuminaro, L.~Chacon, P.~Weber, Scalable
  implicit incompressible resistive {MHD} with stabilized {FE} and
  fully-coupled {N}ewton-{K}rylov-{AMG}, Computer Methods in Applied Mechanics
  and Engineering 304 (2016) 1--25.

\bibitem{Silvester2001}
D.~Silvester, H.~Elman, D.~K.~A. Wathen, Efficient preconditioning of the
  linearized {N}avier-{S}tokes equations for incompressible flow, Journal of
  Computational and Applied Mathematics 128 (2001) 261--279.

\bibitem{Turek1999}
S.~Turek, Efficient Solvers for Incompressible Flow Problems: {A}n Algorithmic
  and Computational Approach, Springer Science \& Business Media, 1999.

\bibitem{Peiro2009}
J.~Peir\'{o}, A.~Veneziani, Reduced models of the cardiovascular system,
  Springer, 2009, pp. 347--394.

\bibitem{Quarteroni2016}
A.~Quarteroni, A.~Veneziani, C.~Vergara, Geometric multiscale modeling of the
  cardiovascular system, between theory and practice, Computer Methods in
  Applied Mechanics and Engineering 302 (2016) 193--252.

\bibitem{Figueroa2006}
C.~Figueroa, I.~Vignon-Clementel, K.~Jansen, T.~Hughes, C.~Taylor, A coupled
  momentum method for modeling blood flow in three-dimensional deformable
  arteries, Computer Methods in Applied Mechanics and Engineering 195 (2006)
  5685--5706.

\bibitem{Moghadam2013a}
M.~Moghadam, Y.~Bazilevs, A.~Marsden, A new preconditioning technique for
  implicitly coupled multidomain simulations with applications to hemodynamics,
  Computational Mechanics 52 (2013) 1141--1152.

\bibitem{Moghadam2015}
M.~Moghadam, Y.~Bazilevs, A.~Marsden, A bi-partitioned iterative algorithm for
  solving linear systems obtained from incompressible flow problems, Computer
  Methods in Applied Mechanics and Engineering 286 (2015) 40--62.

\bibitem{simvascular-simulation-guide}
{S}im{V}ascular {S}imulation {G}uide,
  http://simvascular.github.io/docsFlowSolver.html, accessed: 2019-08-05.

\bibitem{Updegrove2017}
A.~Updegrove, N.~Wilson, J.~Merkow, H.~Lan, A.~Marsden, S.~Shadden,
  Sim{V}ascular: An open source pipeline for cardiovascular simulation, Annals
  of Biomedical Engineering 45 (2017) 525--541.

\bibitem{May2008}
D.~May, L.~Moresi, Preconditioned iterative methods for {S}tokes flow problems
  arising in computational geodynamics, Physics of the Earch and Planetary
  Interiors 171 (2008) 33--47.

\bibitem{Cyr2016}
E.~Cyr, J.~Shadid, R.~Tuminaro, Teko: {A} block preconditioning capability with
  concrete example applications in {N}avier-{S}tokes and {MHD}, SIAM Journal on
  Scientific Computing 38 (2016) S307--S331.

\bibitem{Manguoglu2008}
M.~Manguoglu, A.~Sameh, T.~Tezduyar, S.~Sathe, A nested iterative scheme for
  computation of incompressible flows in long domains, Computational Mechanics
  43 (2008) 73--80.

\bibitem{Manguoglu2011}
M.~Manguoglu, K.~Takizawa, A.~Sameh, T.~Tezduyar, Nested and parallel sparse
  algorithms for arterial fluid mechanics computations with boundary layer mesh
  refinement, International Journal for Numerical Methods in Fluids 65 (2011)
  135--149.

\bibitem{Manguoglu2009}
M.~Manguoglu, A.~Sameh, F.~Saied, T.~Tezduyar, S.~Sathe, Preconditioning
  techniques for nonsymmetric linear systems in the computation of
  incompressible flows, Journal of Applied Mechanics 76 (2009) 021204.

\bibitem{Saad1993}
Y.~Saad, A flexible inner-outer preconditioned {GMRES} algorithm, SIAM Journal
  on Scientific Computing 14 (1993) 461--469.

\bibitem{Liu2019}
J.~Liu, A.~Marsden, A robust iterative method for finite elastodynamics with
  nested block preconditioning, Journal of Computational Physics 383 (2019)
  72--93.

\bibitem{Liu2018}
J.~Liu, A.~Marsden, A unified continuum and variational multiscale formulation
  for fluids, solids, and fluid-structure interaction, Computer Methods in
  Applied Mechanics and Engineering 337 (2018) 549--597.

\bibitem{Liu2019a}
J.~Liu, A.~Marsden, Z.~Tao, An energy-stable mixed formulation for isogeometric
  analysis of incompressible hyper-elastodynamics, International Journal for
  Numerical Methods in Engineering 120 (2019) 937--963.

\bibitem{Taylor1998}
C.~Taylor, T.~Hughes, C.~Zarins, Finite element modeling of blood flow in
  arteries, Computer Methods in Applied Mechanics and Engineering 158 (1998)
  155--196.

\bibitem{Bazilevs2007a}
Y.~Bazilevs, V.~Calo, J.~Cottrell, T.~Hughes, A.~Reali, G.~Scovazzi,
  Variational multiscale residual-based turbulence modeling for large eddy
  simulation of incompressible flows, Computer Methods in Applied Mechanics and
  Engineering 197 (2007) 173--201.

\bibitem{Liu2020a}
J.~Liu, I.~Lan, O.~Tikenogullari, A.~Marsden, On the generalized-$\alpha$
  scheme for the incompressible {N}avier-{S}tokes equations, Computer Methods
  in Applied Mechanics and Engineering In preparation.

\bibitem{Moghadam2013}
M.~Moghadam, I.~Vignon-Clementel, R.~Figliola, A.~Marsden, {MOCHA}, A modular
  numerical method for implicit 0{D}/3{D} coupling in cardiovascular finite
  element simulations, Journal of Computational Physics 244 (2013) 63--79.

\bibitem{Pauli2017}
L.~Pauli, M.~Behr, On stabilized space-time {FEM} for anisotropic meshes:
  {I}ncompressible {N}avier-{S}tokes equations and applications to blood flow
  in medical devices, International Journal for Numerical Methods in Fluids 85
  (2017) 189--209.

\bibitem{Danwitz2019}
M.~von Danwitz, V.~Karyofylli, N.~Hosters, M.~Behr, Simplex space-time meshes
  in compressible flow simulations, International Journal for Numerical Methods
  in Fluids 91 (2019) 29--48.

\bibitem{Takizawa2010}
K.~Takizawa, J.~Christopher, T.~Tezduyar, S.~Sathe, Space-time finite element
  computation of arterial fluid-structure interactions with patient-specific
  data, International Journal for Numerical Methods in Biomedical Engineering
  26 (2010) 101--116.

\bibitem{Franca1992}
L.~Franca, S.~Frey, Stabilized finite element methods: {II}. {T}he
  incompressible {N}avier-{S}tokes equations, Computer Methods in Applied
  Mechanics and Engineering 99 (1992) 209--233.

\bibitem{Figuero2006}
C.~Figueroa, A coupled-momentum method to model blood flow and vessel
  deformation in human arteries: applications in disease research and
  simulation-based medical planning, Ph.D. thesis, Stanford university (2006).

\bibitem{Bazilevs2009a}
Y.~Bazilevs, J.~Gohean, T.~Hughes, R.~Moser, Y.~Zhang, Patient-specific
  isogeometric fluid–structure interaction analysis of thoracic aortic blood
  flow due to implantation of the {J}arvik 2000 left ventricular assist device,
  Computer Methods in Applied Mechanics and Engineering 198 (2009) 3534--3550.

\bibitem{Esmaily-Moghadam2011}
M.~Esmaily-Moghadam, Y.~Bazilevs, T.~Hsia, I.~Vignon-Clementel, A.~Marsden,
  {MOCHA}, A comparison of outlet boundary treatments for prevention of
  backflow divergence with relevance to blood flow simulations, Computational
  Mechanics 48 (2011) 277--291.

\bibitem{Bertoglio2016}
C.~Bertoglio, A.~Caiazzo, A {S}tokes-residual backflow stabilization method
  applied to physiological flows, Journal of Computational Physics 313 (2016)
  260--278.

\bibitem{Bertoglio2018}
C.~Bertoglio, A.~Caiazzo, Y.~Bazilevs, M.~Braack, M.~Esmaily, V.~Gravemeier,
  A.~Marsden, O.~Pironneau, I.~Vignon-Clementel, W.~Wall, Benchmark problems
  for numerical treatment of backflow at open boundaries, International journal
  for numerical methods in biomedical engineering 34 (2018) e2918.

\bibitem{Chung1993}
J.~Chung, G.~Hulbert, A time integration algorithm for structural dynamics with
  improved numerical dissipation: the generalized-$\alpha$ method, Journal of
  applied mechanics 60 (1993) 371--375.

\bibitem{Saad1986}
Y.~Saad, M.~Schultz, {GMRES}: {A} generalized minimal residual algorithm for
  solving nonsymmetric linear systems, SIAM Journal on scientific and
  statistical computing 7 (1986) 856--869.

\bibitem{Falgout2002}
R.~Falgout, U.~Yang, hypre: {A} library of high performance preconditioners,
  in: International Conference on Computational Science, Springer, 2002, pp.
  632--641.

\bibitem{hypre-boomeramg-moose-intro}
{H}ypre/{B}oomer{AMG} preconditioning,
  \url{https://mooseframework.org/application_development/hypre.html},
  accessed: 2019-08-19.

\bibitem{Anderson2018}
R.~Anderson, A.~Barker, J.~Bramwell, J.~Cerveny, J.~Dahm, V.~Dobrev,
  Y.~Dudouit, A.~Fisher, T.~Kolev, M.~Stowell, V.~Tomov, {MFEM}: {A} {M}odular
  {F}inite {E}lement {M}ethods {L}ibrary (2018).

\bibitem{Knoll2004}
D.~Knoll, D.~Keyes, Jacobian-free {N}ewton-{K}rylov methods: a survey of
  approaches and applications, Journal of Computational Physics 193 (2004)
  357--397.

\bibitem{Pernice2001}
M.~Pernice, M.~Tocci, A {M}ultigrid-{P}reconditioned {N}ewton-{K}rylov {M}ethod
  for the {I}ncompressible {N}avier-{S}tokes {E}quations, SIAM Journal on
  Scientific Computing 23 (2001) 398--418.

\bibitem{Patankar1983}
S.~Patankar, D.~Spalding, A calculation procedure for heat, mass and momentum
  transfer in three-dimensional parabolic flows, in: Numerical Prediction of
  Flow, Heat Transfer, Turbulence and Combustion, Elsevier, 1983, pp. 54--73.

\bibitem{Deparis2016}
S.~Deparis, D.~Forti, G.~Grandperrin, A.~Quarteroni, {FaSCI}: A block parallel
  preconditioner for fluid–structure interaction in hemodynamics, Journal of
  Computational Physics 327 (2016) 700--718.

\bibitem{White2011}
J.~White, R.~Borja, Block-preconditioned {N}ewton-{K}rylov solvers for fully
  coupled flow and geomechanics, Computational Geosciences 15 (2011) 647.

\bibitem{Verdugo2016}
F.~Verdugo, W.~Wall, Unified computational framework for the efficient solution
  of n-field coupled problems with monolithic schemes, Computer Methods in
  Applied Mechanics and Engineering 310 (2016) 335--366.

\bibitem{Quarteroni2000}
A.~Quarteroni, F.~Saleri, A.~Veneziani, Factorization methods for the numerical
  approximation of {N}avier-{S}tokes equations, Computer Methods in Applied
  Mechanics and Engineering 188.

\bibitem{Esmaily-Moghadam2015}
M.~Esmaily-Moghadam, Y.~Bazilevs, A.~Marsden, A bi-partitioned iterative
  algorithm for solving linear systems arising from incompressible flow
  problems, Computer Methods in Applied Mechanics and Engineering 286 (2015)
  40--62.

\bibitem{Esmaily-Moghadam2013}
M.~Esmaily-Moghadam, Y.~Bazilevs, A.~Marsden, A new preconditioning technique
  for implicitly coupled multidomain simulations with applications to
  hemodynamics, Computational Mechanics 52 (2013) 1141--1152.

\bibitem{Akguen2001}
M.~Akg{\"u}n, J.~Garcelon, R.~Haftka, Fast exact linear and non-linear
  structural reanalysis and the sherman--morrison--woodbury formulas,
  International Journal for Numerical Methods in Engineering 50 (2001)
  1587--1606.

\bibitem{Stampede2-user-guide}
{S}tampede2 {U}ser {G}uide,
  \url{https://portal.tacc.utexas.edu/user-guides/stampede2}, accessed:
  2019-06-23.

\bibitem{Bergersen2019}
A.~Bergersen, M.~Mortensen, K.~Valen-Sendstad, The {FDA} nozzle
  benchmark:``{I}n theory there is no difference between theory and practice,
  but in practice there is”, International Journal for Numerical Methods in
  Biomedical Engineering 35 (2019) e3150.

\bibitem{Geuzaine2009}
C.~Geuzaine, J.~Remacle, Gmsh: {A} 3-{D} finite element mesh generator with
  built-in pre- and post-processing facilities, International Journal for
  Numerical Methods in Biomedical Engineering 79 (2009) 1309--1331.

\bibitem{Zmijanovic2017}
V.~Zmijanovic, S.~Mendez, V.~Moureau, F.~Nicoud, About the numerical robustness
  of biomedical benchmark cases: interlaboratory {FDA}'s idealized medical
  device, International journal for numerical methods in biomedical engineering
  33 (2017) e02789.

\bibitem{FDA_Critical_Path_Project_webpage}
The {FDA's} ``{C}ritical {P}ath" {C}omputational {F}luid {D}ynamics
  ({CFD})/{B}lood {D}amage {P}roject, \url{https://nciphub.org/wiki/FDA_CFD},
  accessed: 2019-03-28.

\bibitem{petsc-user-ref}
S.~Balay, S.~Abhyankar, M.~Adams, J.~Brown, P.~Brune, K.~Buschelman, L.~Dalcin,
  V.~Eijkhout, W.~Gropp, D.~Kaushik, M.~Knepley, D.~May, L.~McInnes, K.~Rupp,
  P.~Sanan, B.~Smith, S.~Zampini, H.~Zhang, {PETS}c users manual, Tech. Rep.
  ANL-95/11 - Revision 3.8, Argonne National Laboratory (2017).

\bibitem{Si2015}
H.~Si, Tet{G}en, a {D}elaunay-{B}ased {Q}uality {T}etrahedral {M}esh
  {G}enerator, ACM Transactions on Mathematical Software 41 (2015) 11.

\bibitem{AAA-model-url}
Cardiovascular and pulmonary model repository,
  \url{http://www.vascularmodel.com}.

\bibitem{simmetrix}
Simmetrix, \url{http://www.simmetrix.com/}.

\bibitem{Burstedde2009}
C.~Burstedde, O.~Ghattas, G.~Stadler, T.~Tu, L.~Wilcox, Parallel scalable
  adjoint-based adaptive solution of variable-viscosity {S}tokes flow problems,
  Computer Methods in Applied Mechanics and Engineering 198 (2009) 1691--1700.

\end{thebibliography}
\end{document}